\documentclass{m2an}

\usepackage[utf8]{inputenc}
\usepackage{amsmath}
\usepackage{graphicx}
\usepackage{subcaption}
\usepackage{mathSymbols}
\usepackage{geometry}
\usepackage{url}
\usepackage{graphicx}
\usepackage{color}
\usepackage[final]{changes}

\newcommand{\xadded}[1]{#1}
\newcommand{\xdeleted}[1]{}
\newcommand{\xreplaced}[2]{#1}

\usepackage[colorlinks=true, allcolors = blue]{hyperref}

\definechangesauthor[name=R1, color=orange]{R1}
\definechangesauthor[name=R2, color=red]{R2}
\definechangesauthor[name=additional, color=magenta]{additional}

\definechangesauthor[name=R1, color=orange]{R1-new}
\definechangesauthor[name=R2, color=red]{R2-new}

\numberwithin{equation}{section}

\renewcommand{\ie}{\emph{i.e.}, }

\newcommand{\SE}[2]{SE#1#2}

\newcommand{\lambdaLk}{\lambda_{\vecL_k}}

\newcommand{\xiLk}{\xi_{\vecL_k}}

\begin{document}

\title[Analysis and improvement of semi-Lagrangian exponential schemes for the SWE]{Analysis and improvement of a semi-Lagrangian exponential scheme for the shallow-water equations on the rotating sphere}
% \thanks{This work was supported by the São Paulo Research Foundation (FAPESP) grants 2021/03777-2 and 2021/06176-0, as well as the Brazilian National Council for Scientific and Technological Development (CNPq), Grant 303436/2022-0.}\thanks{This project also received funding from the Federal Ministry of Education and Research and the European High-Performance Computing Joint Undertaking (JU) under grant agreement No 955701, Time-X. The JU receives support from the European Union’s Horizon 2020 research and innovation programme and Belgium, France, Germany, Switzerland.}

\runningtitle{Analysis and improvement of semi-Lagrangian exponential schemes for the SWE}

\author{João Guilherme Caldas Steinstraesser}\address{Applied Mathematics, Universidade de Sao Paulo, São Paulo, São Paulo, CEP 05508-090, Brazil}
\author{Martin Schreiber}\address{Université Grenoble Alpes / Laboratoire Jean Kuntzmann, Saint-Martin-d'Hères, 38400, France; Inria AIRSEA team, Grenoble, 38058, France}\secondaddress{Technical University of Munich, Garching, 85748, Germany}
\author{Pedro S. Peixoto}\sameaddress{1}

\date{}



\begin{abstract}

   In this work, we study and extend a class of semi-Lagrangian exponential methods, which combine exponential time integration techniques, suitable for integrating stiff linear terms, with a semi-Lagrangian treatment of nonlinear advection terms. Partial differential equations involving both processes arise for instance in atmospheric circulation models. Through a truncation error analysis, we show that previously formulated semi-Lagrangian exponential schemes are limited to first-order accuracy due to \xadded{the approximation of the integration factor acting on }the discretization of the linear term; we then formulate a new discretization leading to second-order accuracy. Also, a detailed stability study is conducted to compare several Eulerian and semi-Lagrangian exponential schemes, as well as a well-established semi-Lagrangian semi-implicit method, which is used in operational atmospheric models. Numerical simulations of the shallow-water equations on the rotating sphere are performed to assess the orders of convergence, stability properties, and computational cost of each method. The proposed second-order semi-Lagrangian exponential method was shown to be more stable and accurate than the previously formulated schemes of the same class at the expense of larger wall-clock times; however, the method is more stable and has a similar cost compared to the well-established semi-Lagrangian semi-implicit method; therefore, it is a competitive candidate for potential operational applications in atmospheric circulation modeling.
\end{abstract}

\maketitle

\section{Introduction}

\indent Exponential integration methods are a large class of time-stepping schemes that have raised increasing interest for the integration of time-dependent partial differential equations containing stiff linear operators, since they allow integrating exactly (or at least very accurately) the linear terms of the governing equations, thus overcoming stability constraints related to their stiffness \cite{minchev_wright:2005, hochbruck_ostermann:2010}. These methods require the evaluation of matrix exponentials, which is a main challenge for their efficient implementation. The recent development of accurate and efficient linear algebra methods has directly contributed to the development and popularization of exponential integration schemes \cite{minchev_wright:2005}, which have been increasingly seen as viable alternatives to overcome Courant-Friedrichs-Lewy (CFL)-related time step limitations usually found in complex fluid dynamics problems such as in atmospheric modeling \cite{gaudreault_pudykiewicz:2016}. 

\indent The development of exponential integration methods goes back to 1960, by \cite{certaine:1960}, in which the variation-of-constants formula was used to derive the scheme, leading to the family of Exponential Time Differencing (ETD) methods. An alternative approach, proposed by \cite{lawson:1967}, consists of defining a change of variables leading to a differential equation in which the linear operator appears only as arguments of exponential functions, a class of exponential schemes known as Integrating Factor (IF) methods. Both classes are derived considering a partition of the governing equations into linear and nonlinear terms. A third well-known class of exponential schemes is composed of the exponential Rosenbrock methods \cite{hochbruck_ostermann:2006}, which are derived from a linearization of the full time evolution operator, thus requiring the computation of its Jacobian. Since then, several ETD, IF and Rosenbrock exponential methods have been proposed, as well as novel classes of exponential schemes. Detailed reviews on exponential integration methods can be found, \eg in \cite{minchev_wright:2005, hochbruck_ostermann:2010}. These methods have been studied and applied in a very diverse range of contexts. In particular, several recent works have focused on geophysical fluid dynamics, in particular atmospheric circulation models, which is the application domain targeted in this work. \cite{bonaventura:2010, clancy_pudykiewicz:2013} assessed good stability properties and improved accuracy of Rosenbrock-type exponential schemes, compared to more traditional semi-implicit methods, when applied to the shallow-water equations (SWE) on the rotating sphere, but highlighted the need for more efficient approaches for computing the exponential of large matrices. Similar conclusions have been drawn by \cite{brachet_al:2022} comparing various time integration schemes in the context of the one-dimensional SWE. The challenge of improving the exponential computation has been tackled by \cite{gaudreault_pudykiewicz:2016, thai_luan_al:2019,vasylkevych_nedjeljka:2021}. Applications to more complex, three-dimensional atmospheric models have been investigated by \cite{rainwater_al:2023,pudykiewicz_clancy:2022}, with stable and accurate integration of both slow and fast atmospheric dynamics. Finally, methods considering rational approximations to the exponential functions, which lead to parallel-in-time methods, have been investigated for the SWE on the rotating sphere by \cite{schreiber_loft:2018, schreiber_al:2019}.

\indent More recently, \cite{peixoto_schreiber:2019} proposed a class of semi-Lagrangian Runge-Kutta-type ETD schemes (SL-ETDRK) to tackle advective PDEs, inspired by the Eulerian ETDRK methods proposed by \cite{cox_matthews:2002}. Semi-Lagrangian schemes, which follow the trajectories of particles along each time step instead of relying on fixed spatial grids as in Eulerian approaches, are popular methods in atmospheric modeling due to their improved stability properties. This improved stability was indeed verified for the proposed SL-ETDRK schemes in \cite{peixoto_schreiber:2019}, through the numerical simulation of the SWE on the bi-periodic plane with a spectral discretization in space. A semi-Lagrangian exponential method has also been developed by \cite{shashkin_goyman:2020}, based on the Rosenbrock exponential approach, leading to a second-order scheme whose accuracy was confirmed on numerical tests of the SWE on the rotating sphere considering a \xreplaced{cubed-sphere}{cubed grid} discretization.

\indent This work proposes a better analytical understanding of convergence and stability properties of the SL-ETDRK schemes, as well as improvements for these methods. Indeed, although these methods have been derived analogously to their respective Eulerian counterparts ETD1RK and ETD2RK (respectively first- and second-order accurate) \cite{cox_matthews:2002}, the expected second-order convergence of SL-ETD2RK is not verified in practice. The numerical simulations presented by \cite{peixoto_schreiber:2019} reveal that both SL-ETD1RK and SL-ETD2RK are first-order accurate\xreplaced{. Despite some indications that it could be related to the approximation for the integration factor arising on the derivation of the schemes, its influence on the overall accuracy of SL-ETD2RK was not completely clear}{, but the reasons for that were not clear.} Here, we prove that the observed first-order accuracy of SL-ETD2RK is \xreplaced{indeed related to the approximation of the integration factor, more specifically in the }{due to the }discretization of the linear term, which involves both an exponential operator and a spatial interpolation coming from the semi-Lagrangian approach. This issue is naturally absent in Eulerian schemes, in which the exponential integrates exactly the linear term. Then, we propose an alternative discretization, consisting of a splitting of the exponential operator, and prove that it leads indeed to second-order accuracy. We also conduct a detailed analytical stability study of both Eulerian and semi-Lagrangian ETDRK schemes, as well as the popular semi-Lagrangian semi-implicit SETTLS method \cite{hortal:2002} used in operational atmospheric models; we show in particular that standard linear stability analyses are not able to provide concluding insights on the stability of semi-Lagrangian exponential schemes due to a lack of commutation between exponential and interpolation operators; therefore, we conduct an empirical simulation-based stability study following a procedure described by \cite{peixoto_al:2017}. Finally, we perform numerical experiments for the assessment of convergence and stability of the proposed improved schemes, which demonstrate its advantages compared to the previous Eulerian and semi-Lagrangian exponential methods and the semi-Lagrangian semi-implicit SETTLS.

\indent In terms of application, our focus here is on the use of semi-Lagrangian exponential methods in the context of atmospheric modeling. Atmospheric circulation models are characterized by a vast range of waves propagating at different speeds, with the fast ones being responsible for often severe stability constraints, which has motivated the development of various discretization methods searching to establish a good trade-off between stability, accuracy, and computational cost  \cite{williamson:2007}; therefore, combining the semi-Lagrangian and exponential approaches could better handle both advection processes and the stiffness of the linear terms. In particular, we expect semi-Lagrangian exponential methods to be more stable than Eulerian exponential ones, and more accurate than traditional semi-Lagrangian schemes (at least when stiff, linear processes have relatively important magnitude), such that good compromises between computational cost and accuracy could be obtained. We consider, for the numerical experiments and stability analyses conducted here, the SWE on the rotating sphere, using challenging benchmark test cases commonly used in atmospheric modeling research; thus, this work not only analyses and improves the semi-Lagrangian exponential methods proposed by \cite{peixoto_schreiber:2019}, but it also extends their application domain since \cite{peixoto_schreiber:2019} focused on the SWE on the plane. Finally, we also compare the exponential semi-Lagrangian scheme to the well-established SL-SI-SETTLS, in order to assess the competitiveness of semi-Lagrangian exponential methods.

\indent This paper is organized as follows: the shallow-water equations on the rotating sphere are introduced in Section \ref{sec:SWE}; the time integrators considered in this work (the Eulerian ETDRK schemes by \cite{cox_matthews:2002}, the SL-ETDRK schemes by \cite{peixoto_schreiber:2019} and SL-SI-SETTLS by \cite{hortal:2002}) are presented in Section \ref{sec:time_integrators}; a truncation error analysis indicating the sources of the first-order accuracy of the SL-ETD2RK scheme is developed in Section \ref{sec:source_first_order}; an effectively second-order accurate semi-Lagrangian exponential method is proposed in Section \ref{sec:second_order_SLETDRK}, as well as its truncation error and computational complexity analyses; a detailed linear stability study comparing the various schemes considered in this work is conducted in Section \ref{subsec:serial_stability}; numerical simulations of the SWE on the rotating sphere are performed in Section \ref{sec:numerical_simulations} in order to assess and compare stability, convergence and computational cost properties; and conclusions are presented in Section \ref{sec:conclusion}.

\section{The shallow-water equations on the rotating sphere}
\label{sec:SWE}

\indent The shallow-water equations on the rotating sphere are a popular two-dimensional model in the early stages of atmospheric circulation research since it is simplified \wrt more complete, three-dimensional equations but still contains most of the challenges related to the spatial discretization of atmospheric models. Therefore, it allows an easier but insightful study of spatial and temporal discretization schemes for atmospheric models \cite{williamson_al:1992}. The equations read

\begin{equation}
    \label{eq:swe_sphere}
    \begin{aligned}
    \pdertt{\vecU} &= \vecL_{\vecG}(\vecU) + \vecL_{\vecC}(\vecU) + \vecN_{\vecA}(\vecU) + \vecN_{\vecR}(\vecU) + \vecL_{\nu}(\vecU) +
    \vecb,
    \end{aligned}
\end{equation}

\noindent where the solution vector $\vecU = \vecthreeT{\Phi}{\xi}{\delta}$ contains the geopotential $\Phi = \ol{\Phi} + \Phi' = gh$, vorticity  $\xi := \vecz \cdot (\nabla \times \vecV)$ and divergence $\delta := \nabla \cdot \vecV$ fields. $\Phi$, $\xi$ and $\delta$ are functions of $(\lambda, \mu, t)$, where $\lambda$ is the longitude, $\mu = \sin(\theta)$, $\theta$ is the latitude and $t$ is the time. $g$ denotes the gravitational acceleration, $h$ is the fluid depth, $\ol{\Phi}$ and $\Phi'$ are respectively the mean geopotential and the geopotential perturbation, $\vecV := \vectwoT{u}{v}$ is the horizontal velocity field and $\vecz$ is the unit vector in the vertical direction. The RHS of \eqref{eq:swe_sphere} is decomposed as a sum of linear and nonlinear terms accounting for different physical processes: $\vecL_{\vecG}$ accounts for linear processes related to the gravitation; $\vecL_{\vecC}$ contains linear rotation-related phenomena, described by the Coriolis parameter $f = 2\Omega \sin \theta$, where $\Omega$ is the Earth's angular velocity; $\vecN_{\vecA}$ is the nonlinear advection term, $\vecN_{\vecR}$ contains the remaining nonlinear terms, $\vecL_{\nu}(\vecU)$ is a linear term accounting for an artificial (hyper-)viscosity of order $q$ and coefficient $\nu \geq 0$, with $q$ even, and $\vecb$ is related to the topography field $b$. These terms are given by

\begin{equation}
    \label{eq:swe_sphere_terms}
    \begin{gathered}
        \vecL_{\vecG}(\vecU) = \matthree{0}{0}{-\ol{\Phi}}{0}{0}{0}{-\nabla^2}{0}{0}\vecU,\qquad
        \vecL_{\vecC}(\vecU) = \vecthree{0}{-\nabla \cdot(f\vecV)}{\vecz \cdot \nabla \times (f\vecV)}\\
        \vecN_{\vecA}(\vecU) = \vecthree{-\vecV \cdot \nabla \Phi}{-\nabla \cdot (\xi \vecV)}{-\nabla^2 \left( \frac{\vecV \cdot \vecV}{2} \right) + \vecz \cdot \nabla \times (\xi \vecV)},\qquad
        \vecN_{\vecR}(\vecU) = \vecthree{-\Phi'\delta}{0}{0},\\
        \vecL_{\nu}(\vecU) = (-1)^{\frac{q}{2} + 1}\nu 
        \vecthree{\nabla^q \Phi'}{ \nabla^q \xi}{\nabla^q \delta }, \qquad \vecb = \vecthree{0}{0}{-g\nabla^2 b}.
    \end{gathered}
\end{equation}

% \indent A detailed description of how each time stepping scheme considered in this work treats each term is provided in Section \ref{sec:time_integrators}.

 % \indent In all the time stepping schemes described below, the governing equations are solved considering the terms $\vecL_{\vecG}$, $\vecL_{\vecC}$, $\vecN_{\vecA}$ and $\vecN_{\vecR}$. The remaining term, $\vecL_{\nu}$, accounting for the viscosity, is solved as a post-processing step, with the application of an implicit solver at the end of each time step of the main integrator.

 \subsection{Spectral discretization}

\indent In this work, the spatial dimension is discretized considering a spectral approach using spherical harmonics, which is briefly described below. Details can be found in \xadded{\cite{swarztrauber:1996}}. Let $Y_{m,n}$ denote the spherical harmonic function of zonal and total wavenumbers $m$ and $n$, respectively. The spherical harmonics expansion of a given smooth field $\psi$ reads

\begin{equation}
    \label{eq:spherical_harmonics}
    \psi(\lambda, \mu, t) = \sum_{m=-\infty}^{\infty} \sum_{n  = |m|}^{\infty} \psi_{m,n}(t) Y_{m,n}(\lambda,\mu)
\end{equation}

\noindent where $\psi_{m,n}$ is the expansion coefficient corresponding to the wavenumber pair $(m,n)$. We consider a triangular truncation with resolution $M$:

\begin{equation}
    \label{eq:spherical_harmonics_truncated}
    \psi(\lambda, \mu, t) = \sum_{m=-M}^{M} \sum_{n  = |m|}^{M} \psi_{m,n}(t) Y_{m,n}(\lambda,\mu)
\end{equation}

\indent As described in the next section, the use of a spectral discretization based on spherical harmonics, combined with a proper arrangement of the terms of the governing equations, is essential for ensuring an efficient computation of matrix exponentials in both Eulerian and semi-Lagrangian exponential methods.

%\section{Eulerian and Semi-Lagrangian exponential integration methods}

\section{Time integrators}
\label{sec:time_integrators}

\indent We present in this section the time stepping schemes considered in the work. We begin by briefly presenting the Eulerian exponential integration methods ETD1RK and ETD2RK \cite{cox_matthews:2002}, on which their semi-Lagrangian versions are based. Then, as an introduction to semi-Lagrangian schemes, we present the semi-Lagrangian semi-implicit SETTLS \cite{hortal:2002}, which we will also compare in the numerical experiments. Finally, we present the semi-Lagrangian exponential schemes proposed by \cite{peixoto_schreiber:2019}. As explained below, we propose a modification to the name of this last class of methods in order to take into account their actual order of convergence.

\subsection{Eulerian Runge-Kutta-type exponential integration methods (ETDRK)}
\label{subsec:ETDRK}

\indent Exponential Time Differencing (ETD) schemes are a class of exponential integration methods to solve time-dependent problems of the form

\begin{equation}
    \label{eq:base_equation_exp}
    \dert{\vecU} = \vecL \vecU + \vecN(\vecU), \qquad \vecU(t = 0) = \vecU_0,
\end{equation}

\noindent which are obtained, \eg from a spatial discretization of a PDE, where $\vecU \in \reals^N$, $\vecL \in \reals^{N \times N}$ is a linear operator and $\vecN: \reals^N \rightarrow \reals^N$ is a nonlinear function. As a common feature, exponential methods integrate the linear term $\vecL \vecU$ very accurately (or exactly, down to machine precision, if the exponential terms can be evaluated exactly), by writing its evolution using the exponential of the linear operator, whereas an approximation is proposed to compute the nonlinear term. In the following paragraphs, we briefly describe the main ideas behind exponential methods, more specifically concerning the schemes considered in this work, namely ET1DRK and ETD2RK \cite{cox_matthews:2002}. We refer the reader to \cite{cox_matthews:2002, hochbruck_ostermann:2010} for further details on exponential integration methods and\xreplaced{, in particular, to \cite{peixoto_schreiber:2019} for their derivation through the solution of an integration factor problem, which is also the starting point for the derivation of semi-Lagrangian exponential schemes in Section \ref{subsec:sletdrk}.}{their derivation}.

\indent The exact solution of \eqref{eq:base_equation_exp} at time $t_{n+1} = t_n + \Dt$, integrated from $\vecU(t_n)$, reads

\xadded
{
\begin{equation}
    \label{eq:exact_exp}
    \begin{aligned}
    \vecU(t_{n+1}) & = P_n^{-1}(t_{n+1}) \vecU(t_n) + P_n^{-1}(t_{n+1}) \int_{t_n}^{t_{n+1}}P_n(s)\vecN(\vecU(s))ds \\
    & = e^{\Dt \vecL} \vecU(t_n) + \int_{t_n}^{t_{n+1}} e^{-(s-t_n)\vecL}\vecN(\vecU(s))ds,
    \end{aligned}
\end{equation}
}

\noindent which is the well-known variation-of-constants formula, \xreplaced{where, assuming that $\vecL$ does not depend on time, the integration factor}{and can be obtained by solving an integration factor problem}

\xadded{
\begin{equation}
    \label{eq:integration_factor_eulerian}
    P_n(t) = e^{-(t-t_n)\vecL} \in \reals^{N \times N}
\end{equation}
}

\noindent \xadded{is the solution of the integration factor problem}

\xadded
{
\begin{equation}
    \label{eq:integration_factor_problem}
    \dert{P_n}(t) = -P_n(t)L, \qquad P_n(t_n) = I.
\end{equation}
}

% \noindent \xadded{see \cite{peixoto_schreiber:2019} for details on its derivation. }

\indent The first term in \eqref{eq:exact_exp} is the exact integration of the linear term. Different exponential methods can be formulated proposing different approximations to the integral containing the evolution of the nonlinear term. In particular, ETDRK methods are a class of Runge-Kutta type exponential schemes proposed by \cite{cox_matthews:2002}, with improved stability properties \wrt previous exponential methods. The first- and second-order ETDRK schemes, named respectively ETD1RK and ETD2RK, read

\begin{subequations}
    \begin{align}
    \vecU^{n+1}_{\text{ETD1RK}} &= \varphi_0(\Dt \vecL) \vecU^n + \Dt \varphi_1(\Dt L)\vecN(\vecU^n) \label{eq:ETD1RK},\\
    \vecU^{n+1}_{\text{EDT2RK}} &= \vecU^{n+1}_{\text{ETD1RK}} + \Dt \varphi_2(\Dt \vecL)\left[ \vecN(\vecU^{n+1}_{\text{ETD1RK}}) - \vecN(\vecU^n)     \label{eq:ETD2RK}\right],
    \end{align}
\end{subequations}
\noindent where 

\begin{equation}
    \label{eq:phin}
    \begin{aligned}
    \varphi_0(z) & = e^{z}, \\ \varphi_k(z) & = z^{-1}\left( \varphi_{k-1}(z) - \varphi_{k-1}(0) \right), \qquad k \geq 1
    \end{aligned}
\end{equation}

\noindent with singularities at $z = 0$ being computed via series expansions in its neighborhood. \xadded{The functions $\varphi_k$ are defined both when $z$ is a scalar or a matrix (in particular, we are interested in the case $z = \Dt \vecL$).} The integrals in the first- and second-order methods are approximated by considering respectively constant and linear approximations to the nonlinear function along $[t_n, t_{n+1}]$:

\begin{equation}
    \label{eq:nonlinear_terms_sletdrk}
    \begin{gathered}
        \vecN(\vecU(s)) \approx \vecN(\vecU(t_n))\\
        \vecN(\vecU(s)) \approx \vecN(\vecU(t_n)) + \frac{s-t_n}{\Dt}\left[ \vecN(\vecU_{\text{ETD1RK}}(t_{n+1})) - \vecN(\vecU(t_n)) \right]
    \end{gathered}
\end{equation}

\indent As mentioned above, these schemes solve the linear term exactly, if the exponential terms can be computed exactly (down to machine precision), which is the case, \eg of scalar equations. However, in the case of systems of equations, this is possible only if the matrix exponentials $e^z = \sum_{k=0}^{\infty} z^k/k!$ can be computed exactly. In general, the matrix exponential must be evaluated approximately, and the efficiency of traditional algorithms limits their practical implementation in geophysical contexts \cite{saad:1980, gaudreault_pudykiewicz:2016} thus requiring the development of more efficient algorithms. Some examples, that have been applied to atmospheric modeling, are Krylov subspace methods with incomplete orthogonalization \cite{gaudreault_pudykiewicz:2016} and rational approximation to exponential integrators \cite{haut_al:2015, schreiber_al:2019}. Nevertheless, if $z$ is diagonalizable, then the matrix exponential can be computed exactly \cite{peixoto_schreiber:2019}: if $z = Q\Lambda Q^{-1}$, where the columns of $Q$ are the eigenvectors of $z$ and $\Lambda$ is a diagonal matrix containing the respective eigenvalues, then

\begin{equation}
    \label{eq:exp_diag}
    \varphi_k(z) = Q \varphi_k(\Lambda) Q^{-1}, \qquad k \geq 0,
\end{equation}

\noindent where $\varphi_k(\Lambda)$ can be computed elementwise since $\Lambda$ is diagonal.

\indent We make use of this property to evaluate the matrix exponentials in a relatively efficient manner. We remark that, in the case of the SWE on the rotating sphere, the linear term $\vecL_{\vecG} + \vecL_{\vecC}$ is not diagonalizable, since the Coriolis parameter is spatially dependent. In the case where $f$ is constant (\eg in the so-called $f$-plane approximation), $\vecL_{\vecG} + \vecL_{\vecC}$ is diagonalizable, which is taken into account by \cite{peixoto_schreiber:2019} to apply exponential schemes to the SWE on the plane. In this work, we overcome this issue by including the Coriolis effects into the nonlinear term, \ie we consider as linear and nonlinear terms respectively $\vecL := \vecL_{\vecG}$ and $\vecN := \vecL_{\vecC} + \vecN_{\vecA} + \vecN_{\vecR} + \vecb$. Alternative splitting between linear and nonlinear terms is also done in operational models, \eg the Integrated Forecast System of the European Centre for Medium-Range Weather Forecast (IFS-ECMWF) \cite{ECMWF:2003}, in which the linear Coriolis term is incorporated into the nonlinear advection.

\indent Since we consider a spectral discretization using spherical harmonics, the time integration \eqref{eq:base_equation_exp} is performed for the spectral coefficients $\vecU_{m,n} := \vecthreeT{\Phi_{m,n}}{\xi_{m,n}}{\delta_{m,n}} $ with zonal and total wavenumbers $m$ and $n$, respectively. In the spectral space, $\vecL\vecU$ becomes, for each wavenumber pair $(m,n)$, $ \mc{L} \vecU_{n,m}$, where $\mc{L}$ is the symbol matrix of $\vecL = \vecL_{\vecG}$, given by:

\begin{equation}
    \label{eq:linear_operator_spectral}
    \mc{L} = \matthree{0}{0}{-\ol{\Phi}}{0}{0}{0}{-D_n}{0}{0}
\end{equation}

\noindent where $D_n := -n(n+1)/a^2$ and $a$ is the Earth's radius. The eigenvalues of the symbol matrix, which are the diagonal elements of $\Lambda$ in \eqref{eq:exp_diag}, read $-\sqrt{-D_n \ol{\Phi}}$, $0$ and $\sqrt{-D_n \ol{\Phi}}$.

\indent An alternative approach to obtain a diagonalizable linear operator would be to use spectral methods considering Hough harmonics \cite{longuet-higgins:1968}, which are eigenmodes of $\vecL_{\vecG} + \vecL_{\vecC}$, as basis functions. This is done by \cite{vasylkevych_nedjeljka:2021} to solve the SWE on the rotating sphere using exponential methods. However, implementing the Hough decomposition is much less straightforward and computationally more expensive than the spherical harmonics one, for which there are efficient available open-source libraries, \eg SHTns \cite{schaeffer:2013}; therefore, we consider the latter in this work.

\indent In terms of computational efficiency, it is important to notice that, for a fixed problem (\ie for a fixed linear operator $\vecL$), the $\varphi_k$ functions evaluated at $\Dt \vecL$, as used in the ETDRK schemes, depend only on the time step. Therefore, if the time step is constant during a numerical simulation, these functions can be precomputed. This may lead to a significant reduction in terms of computing time, since the $\varphi_k$ functions need to be computed for every pair of wavenumbers $(m,n)$. 

\subsection{SL-SI-SETTLS}
\label{subsec:settls}

\indent Semi-Lagrangian (SL) schemes are spatio-temporal discretization methods for advective PDEs that combine the Eulerian and Lagrangian approaches. The former is based on fixed spatial grids and usually has restrictions on the time step size imposed by Courant-Friedrichs-Lewy (CFL) stability conditions, and the latter consists of following fluid particles along their trajectories, which allows the use of larger time steps but is more complex in terms of implementation compared to fixed grid approaches. The SL approach follows the Lagrangian trajectories only along individual time steps: the trajectories arriving at each point of the fixed spatial grid at time $t_{n+1}$ are estimated, and a spatial interpolation procedure is performed to retrieve the solution at the departure points at time $t_n$. This simplified Lagrangian approach, which still allows the use of relatively large time step sizes \cite{staniforth_cote:1991}, explains the popularity of SL methods, which are used by several weather and climate agencies \cite{mengaldo_al:2018}.

\indent An important example of SL methods is the semi-Lagrangian semi-implicit SETTLS (SL-SI-SETTLS), a second-order scheme used in the IFS-ECMWF \cite{ECMWF:2003}. This scheme uses a semi-implicit (Crank-Nicolson) approach for the linear terms and the Stable Extrapolation Two-Time-Level Scheme (SETTLS), proposed by \cite{hortal:2002}, for estimating the Lagrangian trajectories $(t, \vecx(t))$ at each time step $[t_n, t_{n+1}]$. The trajectories between the departure point $(t_n, \vecx_d)$ and the arrival (grid) point $(t_{n_1}, \vecx_j)$ are considered to be linear and determined by the velocity field $\vecv$ at its midpoint $(t_{n+1/2}, \vecx_m)$, where $t_{n+1/2} := t_n + \Dt / 2$. This velocity field, on its turn, is approximated by a linear interpolation between $\vecv(t_n, \vecx_d)$ and $\vecv(t_{n+1}, x_j)$, with this last one being estimated via a linear extrapolation using the two previous time steps:

\begin{equation}
    \label{eq:SETTLS}
    \vecx_j - \vecx_d = \Dt \vecv(t_{n+1/2}, \vecx_m) = \frac{\Dt}{2} \left[ 2\vecv(t_n,\vecx_d) - \vecv(t_{n-1},\vecx_d) + \vecv(t_n,\vecx_j) \right].
\end{equation}

\indent For each grid point $\vecx_j$, Eq. \eqref{eq:SETTLS} is solved iteratively to obtain the departure point $\vecx_d$. Then, SL-SI-SETTLS applied to \eqref{eq:swe_sphere} reads

\begin{equation}
    \label{eq:SL_SI_SETTLS}
    \xadded
    {
    \frac{\vecU^{n+1} - \vecU^n_*}{\Dt} = \frac{1}{2}\left(\vecL \vecU^{n+1} + (\vecL \vecU^n)^n_*\right) + \frac{1}{2} \left([2\tilde{\vecN}(\vecU^n) - \tilde{\vecN}(\vecU^{n-1}) ]^n_* + \tilde{\vecN}(\vecU^{n})\right)
    }
\end{equation}

\noindent where quantities denoted by the subscript $*$ are interpolated to the departure point $\vecx_d$ at time instant $t_n$. As in the IFS-ECMWF model, we consider the linear term to include only gravity-related processes, \ie $\vecL := \vecL_{\vecG}$; also, we consider both the Coriolis- and the bathymetry-related linear terms to be treated as nonlinear ones, such that $\tilde{\vecN} := \vecN_{\vecR} + \vecL_{\vecC} + \vecb$.

\subsection{SL-ETDRK}
\label{subsec:sletdrk}

\indent A semi-Lagrangian version of the ETDRK schemes is proposed by \cite{peixoto_schreiber:2019}, aiming to tackle problems of the form

\begin{equation}
    \label{eq:base_equation_slexp}
    \frac{\xDif\vecU}{\xDif t} = \vecL \vecU + \tilde{\vecN}(\vecU), \qquad \vecU(t = 0) = \vecU_0,
\end{equation}

\noindent where $\xDif{}/\xDif{t} := \partial/\partial t + \vecV \cdot \nabla$ is the material derivative, $\vecV$ is a velocity field, $\tilde{\vecN}$ is a nonlinear term not containing the nonlinear advection $\vecV \cdot \nabla$ and $\vecU = \vecU(t, \vecx(t))$ is the solution defined on the Lagrangian trajectories. The schemes are derived analogously to the Eulerian case (Section \ref{subsec:ETDRK}): the exact solution of \eqref{eq:base_equation_slexp} reads

\begin{equation}
    \label{eq:exact_solution_slexp}
    \vecU(t_{n+1}, \vecx(t_{n+1})) = P_n^{-1}(t_{n+1})\vecU(t_n, \vecx(t_n)) + P_n^{-1}(t_{n+1}) \int_{t_n}^{t_{n+1}} P_n(s) \tilde{\vecN}(\vecU(s, \vecx(s)))ds,
\end{equation}

% \noindent where $P_n(t)$ is the solution of an integration factor problem \xadded{along the Lagrangian trajectory} in $[t_n, t_{n+1}]$, assumed to exist and be invertible for all times. \xadded{If $L$ depends on the space, then it depends on the specific Lagrangian trajectory $(t,\vecx(t))$. Note that this is the case of the SWE on the rotating sphere \eqref{eq:swe_sphere}, even if the Coriolis effects are not included in the linear operator (as considered in this work), since the spherical Laplacian is spatial dependent in usual spherical coordinates.} Let $\vecL_{\vecx(t)}$ denote the trajectory-dependent linear operator. If it satisfies $\vecL_{\vecx(t)}\vecL_{\vecx(s)} = \vecL_{\vecx(s)}\vecL_{\vecx(t)}$ for every $t,s$, then $P_{n}$ reads

\noindent where $P_n(t)$ is the solution of an integration factor problem \xadded{along the Lagrangian trajectory} in $[t_n, t_{n+1}]$, assumed to exist and be invertible for all times. \xadded{If $\vecL$ depends on the space, then it depends on the specific Lagrangian trajectory $(t,\vecx(t))$. Note that this is the case of the SWE on the rotating sphere \eqref{eq:swe_sphere}, since $\vecL$ is a differential operator; if $\vecL$ were constant, then the derivative (\ie the application of $\vecL$) would be constant along the trajectory.} Let $\vecL_{\vecx(t)}$ denote the trajectory-dependent linear operator. If it satisfies $\vecL_{\vecx(t)}\vecL_{\vecx(s)} = \vecL_{\vecx(s)}\vecL_{\vecx(t)}$ for every $t,s$, then $P_{n}$ reads

\begin{equation}
    \label{eq:integration_factor}
    P_n(t) = e^{- \int_{t_n}^{t}\vecL_{\vecx(s)}ds}.
\end{equation}

\indent In \cite{peixoto_schreiber:2019}, the linear term is assumed to be constant along each trajectory and each time step, such that the integration factor assumes the same form as in the Eulerian exponential schemes presented in Section \ref{subsec:ETDRK}. The schemes analogous to ETD1RK and ETD2RK were named respectively SL-ETD1RK and SL-ETD2RK, in which first- and second-order approximations are proposed to the nonlinear term. However, as will be shown in Section \ref{sec:source_first_order}, the approximation of the linear term is actually a first-order one, such that both schemes are globally first-order. Therefore, we rename these methods respectively as \emph{\SE{1}{1}} and \emph{\SE{1}{2}}, where the first and second numbers refer to the orders of approximation of the linear and nonlinear terms, respectively. These schemes read

\begin{subequations}
    \begin{align}
    \vecU^{n+1}_{\text{\SE{1}{1}}} & = \varphi_0(\Dt \vecL) \left[ \vecU^n + \Dt \psi_1(\Dt \vecL) \tilde{\vecN}(\vecU^n) \right]_*^n \label{eq:SL_ETD1RK}, \\      
    \vecU^{n+1}_{\text{\SE{1}{2}}} & =  \vecU^{n+1}_{\text{\SE{1}{1}}} + \Dt \varphi_0(\Dt \vecL ) \left[\psi_2(\Dt \vecL) \tilde{\vecN}(\vecU^{n+1}_{\text{\SE{1}{1}}}) - \left( \psi_2(\Dt \vecL) \tilde{\vecN} (\vecU^n)        \right)^n_* \right], \label{eq:SL_ETD2RK}
    \end{align}
\end{subequations}

\noindent where the functions $\psi_k$ also satisfy the property \eqref{eq:exp_diag} and are defined by

\begin{equation}
    \label{eq:psin}
    \psi_k(z) := (-1)^{k+1}\varphi_k(-z) + \sum_{l=1}^{k-1}\varphi_l(-z), \qquad k \geq 1.
\end{equation}

\indent The semi-Lagrangian trajectories are estimated as in SL-SI-SETTLS, by solving \eqref{eq:SETTLS} iteratively. Finally, in \eqref{eq:SL_ETD1RK} and \eqref{eq:SL_ETD2RK}, the nonlinear term does not contain the advection effects, which are treated by the SL approach. Thus, the only nonlinear term included in $\tilde{\vecN}$ is $\vecN_{\vecR}$, such that $\tilde{\vecN} = \vecN_{\vecR} + \vecL_{\vecC} + \vecb$. As in the Eulerian exponential schemes, the treatment of the Coriolis effects as a nonlinear term allows us to compute the matrix exponentials accurately using \eqref{eq:exp_diag}, since $\vecL := \vecL_{\vecG}$. We remark that, in the SL framework, the Coriolis term could also be dealt with by considering it as an advected quantity \cite{temperton:1997}.

\section{Source of the first order accuracy of \SE{1}{2}}
\label{sec:source_first_order}

\indent We demonstrate in the following paragraphs that the term $\varphi_0(\Dt L)[\vecU^n]^n_*$ is responsible for the first-order accuracy of \SE{1}{2}, despite a second-order approximation to the nonlinear term. More specifically, it is due to the approximation of the integral in the solution \eqref{eq:integration_factor} of the integration factor problem by considering the integrand to be constant in $[t_n, t_{n+1}]$. For that, we compute the truncation error of the scheme

\begin{equation}
    \label{eq:exp_linear}
    U^{n+1} = \varphi_0(\Dt L)[U^n]^n_*
\end{equation}

\noindent approximating the linear PDE,

\begin{equation}
    \label{eq:PDE}
    \frac{\xDif{u}}{\xDif{t}} = Lu, \qquad \frac{\xDif{}}{\xDif{t}} := \pdert{} + v\pderx{},
\end{equation}

\noindent where $v = v(t,x)$ is a velocity field and $L = L(x)$ is a space-dependent linear operator. For the sake of conciseness, we restrict ourselves to the one-dimensional case in space, with a homogeneous spatial discretization with mesh size $\Dx$, and we assume that all quantities are scalars. In Section \ref{sec:second_order_SLETDRK}, in which we propose a second-order method, some remarks are necessary in the case of systems of PDEs.

% \indent We first show that the constant value approximating the integrand in \eqref{eq:integration_factor} is $L_{x(t_{n+1})}$, \ie the linear operator evaluated at the arrival point of the Lagrangian trajectory. 

\indent Along each Lagrangian trajectory $(t,x(t))$, \eqref{eq:PDE} reduces to the ODE

\begin{equation}
    \label{eq:ODE_scalar}
    \frac{du}{dt}(t,x(t)) = L(x(t))u(t,x(t)).
\end{equation}

\begin{figure}[!htbp]
    \centering
    \includegraphics{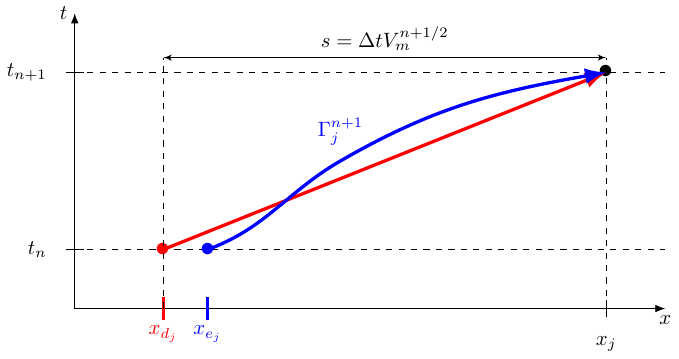}
    \caption{\xadded{Exact Lagrangian trajectory $\Gamma^{n+1}_j$ (in blue) in $[t_n, t_{n+1}]$, and estimated trajectory provided by SETTLS (in red), where $V^{n+1/2}_m$ is an estimation of the velocity field on its midpoint $(t_{n+1/2}, x_m)$ (see Eq. \eqref{eq:SETTLS}).}}
    \label{fig:lagrangian_trajectory}
\end{figure}

\indent Consider a spatial discretization of \eqref{eq:ODE_scalar}, with $P+1$ points $x_0, \dots, x_{P}$ and homogeneous step $\Dx$. Let $\Gamma_j^{n+1}(t) := (t,x^j(t))$ be the Lagrangian trajectory in $[t_n, t_{n+1}]$, arriving at the grid point $x_j$ at time $t_{n+1}$, such that $x^j(t_{n+1}) = x_j$. Under this notation, we also define $u^n_{e_j} := u(t_n, x_{e_j}) := u(t_n, x^j(t_n))$ the exact solution computed at the exact departure point $x_{e_j}$ of the trajectory $\Gamma_j^{n+1}$; and $u^n_{d_j} := u(t_n, x_{d_j})$ the exact solution computed at the approximated departure point $x_{d_j}$ of the trajectory $\Gamma_j^{n+1}$. The linear operator defined along $\Gamma^{n+1}_j$ is denoted by $L_{x^j(t)}$. Some of these definitions are presented in Figure \ref{fig:lagrangian_trajectory}.

\indent This discretization leads to a system of uncoupled ODEs which are approximated by \eqref{eq:exp_linear} as

\begin{equation}
    \label{eq:exp_linear_j}
    U^{n+1}_j = \varphi_0(\Dt L_j) U^n_{*_j}, \qquad j = 0, \dots, P
\end{equation}

\noindent where $w_{*_{j}}$, for a given function $w$ defined on the spatial grid, denotes interpolation to the estimated departure point $x_{d_j}$, and $L_j := L_{x^j(t_{n+1})} = L(x_j)$.

\indent Using cubic Lagrange interpolation to the approximated departure point $x_d$, it is well known that $w_{*_j} = w_{d_j} + \mc{O}(\Dx^4)$ if $w$ is smooth enough. It is also well known that SETTLS \eqref{eq:SETTLS} provides a second-order approximation to the solution computed at the true departure point $x_{e_j}$, more precisely, $u_{d_j} = u_{e_j} + \mc{O}(\Dt^3) + \mc{O}(\Dt\Dx^2)$ \cite[Chapter 7]{durran:2010}. Combining these results, we get, in the case of SETTLS using cubic Lagrange interpolation,

\begin{equation}
    \label{eq:error_cubic_interpolation}
    u_{*_j} = u_{d_j} + \mc{O}(\Dx^4) = u_{e_j} + \mc{O}(\Dt^3) + \mc{O}(\Dt\Dx^2) + \mc{O}(\Dx^4)
\end{equation}

% \indent Then, substituting the exact solution $u$ in \eqref{eq:exp_linear} applied to the $j$-th equation of \eqref{eq:PDE_scalar_semidiscrete} yields

% \begin{equation}
%     \label{eq:first_order_exact}
%     u^{n+1}_j = \varphi_0(\Dt L_j) u^n_{*_j} = \varphi_0(\Dt L_j) u_{e_j} + \mc{O}(\Dt^3) + \mc{O}(\Dx^3)
% \end{equation}

% \indent Then, by defining the matrix $B := (\beta_{j,k}) \in \reals^{P \times P}$, $a = a(x)$ a given spatial dependent function, and $A = \text{diag}(a_0, \dots, a_{P-1}) \in \reals^{P \times P}$, where $a_j = a(x_j)$, we can write

% \begin{equation}
%     \label{eq:linearoperator_interpolation}
%      A\vecu_* = A\left(B\vecu\right) \implies \left(A\vecu_*\right)_j = a_j \sum_{k=-t}^q\beta_{j,k} u_{k} = a_j \left(u_{d_j} + \mc{O}(\Dx^3)\right) = a_j u_{d_j} + \mc{O}(\Dx^3)
% \end{equation}

\noindent From \eqref{eq:exp_linear_j}, it is clear that the exponential function is computed at the arrival point $x_j = x^j(t_{n+1})$ of each \replaced[id=R1-new]{Lagrangian}{Larangian} trajectory. In other words, the integral in $P_n(t_{n+1})$ (see Eq. \eqref{eq:integration_factor}) is approximated by considering the integrand to be constant in $[t_n,t_{n+1}]$ and equal to $L_{x^j(t_{n+1})}$. Note that it is coherent with the interpretation of the order in which the operators are applied in \eqref{eq:exp_linear}: the solution is obtained at $t_{n+1}$ through interpolation, then the linear operator, also evaluated at $t_{n+1}$, is applied to it. This statement also allows us to evaluate the error in the approximation of the integral in \eqref{eq:integration_factor}. From the theory of numerical integration using Newton-Cotes formulae (see \eg \cite[chapter 3]{stoer_bulirsch:2002}), it is well known that 

\begin{equation}
    \label{eq:error_integration_first_order}
    \int_{t_n}^{t_{n+1}} L_{x^j(s)}ds = \Dt L_{x^j(t_{n+1})} + C_{1,j}\Dt^2 = \Dt L_j + C_{1,j}\Dt^2
\end{equation}

\noindent for a given constant $C_{1,j}$ depending on the first derivative of $L$ along the trajectory $\Gamma^{n+1}_j$. Replacing \eqref{eq:error_integration_first_order} in \eqref{eq:integration_factor} and denoting by $P_{n_j}(t)$ the integration factor \eqref{eq:integration_factor} defined along the Lagrangian trajectory $\Gamma^{n+1}_j$, we obtain

\begin{equation}
    \label{eq:integration_factor_first_order}
    P_{n_j}(t_{n+1}) = e^{-\Dt L_j}e^{-C_{1,j}\Dt^2}.
\end{equation}

\indent We are now able to obtain the local truncation error $\tau^{n+1}_j$ at $(t_{n+1},x_j)$ of \eqref{eq:exp_linear} applied to \eqref{eq:PDE}. The exact solution of \eqref{eq:ODE_scalar} along $\Gamma^{n+1}_j$ reads

\begin{equation}
    \label{eq:exact_solution_linear_j}
    u^{n+1}_j = u(t_{n+1}, x^j(t_{n+1})) = P_{n_j}^{-1}(t_{n+1})u(t_n,x^j(t_n)) = P_{n_j}^{-1}(t_{n+1})u^n_{e_j}.
\end{equation}

\indent Replacing the exact solution into the $j$-th equation of \eqref{eq:exp_linear_j}, and using \eqref{eq:error_cubic_interpolation}, \eqref{eq:integration_factor_first_order} and \eqref{eq:exact_solution_linear_j}, yields

\begin{equation}
    \label{eq:LTE_first_order}
    \begin{aligned}
        \tau^{n+1}_j & = \frac{1}{\Dt} \left[ u^{n+1}_j - \varphi_0(\Dt L_j) u^n_{*_j} \right]\\
        & = \frac{1}{\Dt}\left[ P^{-1}_{n_j}(t_{n+1}) u^n_{e_j} - e^{\Dt L_j} u^n_{*_j}\right]\\
        &= \frac{1}{\Dt} \left[e^{\Dt L_j }e^{C_{1,j}\Dt^2} u^n_{e_j} - e^{\Dt L_j}\left( u^{n}_{e_j} + \mc{O}(\Dt^3) + \mc{O}(\Dt\Dx^2) + \mc{O}(\Dx^4)\right)\right]\\
        & = \frac{1}{\Dt} \left[  e^{\Dt L_j} \left(1 + C_{1,j}\Dt^2 + \mc{O}(\Dt^4)\right)u^{n}_{e_j} - e^{\Dt L_j}\left( u^{n}_{e_j} + \mc{O}(\Dt^3) + \mc{O}(\Dt\Dx^2) + \mc{O}(\Dx^4)\right) \right]\\
        & = \mc{O}(\Dt) + \mc{O}(\Dx^2) + \mc{O}\left(\frac{\Dx^4}{\Dt}\right),
    \end{aligned}
\end{equation}

\noindent confirming that the scheme \eqref{eq:exp_linear} provides a first-order approximation to \eqref{eq:PDE}. Note that this first-order is due exclusively to the term $e^{C_{1,j}\Dt^2}$, which arises from the approximation to the integral in the integration factor $P_{n_j}(t_{n+1})$; all other approximations are second-order accurate.

\section{Second-order semi-Lagrangian exponential integration method}
\label{sec:second_order_SLETDRK}

\indent The truncation error analysis presented above indicates that the approximation of the integral in the integration factor \eqref{eq:integration_factor} is responsible for the first-order accuracy of the \SE{1}{2} scheme. Therefore, we provide an alternative formulation for this approximation which is effectively a second-order one. A natural approach is to consider an $\mc{O}(\Dt^3)$ approximation the integral; for instance, the trapezoidal rule. In the scalar case, using the notation defined in the previous section, it reads

\begin{equation}
    \label{eq:error_integration_second_order}
    \int_{t_n}^{t_{n+1}} L_{x^j(s)}ds = \frac{\Dt}{2} \left[L_{x^j(t_{n+1})} + L_{x^j(t_{n})} \right] + C_{2,j}\Dt^3 = \frac{\Dt}{2}\left[L_j + L_{e_j}\right] + C_{2,j}\Dt^3
\end{equation}

\noindent for a constant $C_{2,j}$ depending on the second derivative of $L$ along $\Gamma^{n+1}_j$, and $L_{e_j} := L_{x^j(t_n)} = L(x_{e_j})$. Then, the integration factor \eqref{eq:integration_factor} along the Lagrangian trajectory $\Gamma^{n+1}_j$ satisfies

\begin{equation}
    \label{eq:integration_factor_second_order}
    P_{n_j}(t_{n+1}) = e^{-\int_{t_n}^t L_{x^j(s)}ds} = e^{-\frac{\Dt}{2}\left[L_j + L_{e_j}\right] - C_{2,j}\Dt^3} = e^{-\frac{\Dt L_j}{2}} e^{-\frac{\Dt L_{e_j}}{2}} e^{-C_{2,j}\Dt^3},
\end{equation}

\noindent yielding the scheme

\begin{equation}
    \label{eq:exp_linear_second_order}
    U^{n+1} = \varphi_0\left(\frac{\Dt L}{2}\right)\left[\varphi_0\left(\frac{\Dt L}{2}\right)U^n\right]^n_*
\end{equation}

\indent Note that the last equality in \eqref{eq:integration_factor_second_order} holds because $L_j$ and $L_{e_j}$ are scalars. The order in which the operators are applied in \eqref{eq:exp_linear_second_order} (application of the exponential inside the brackets, followed by a spatial interpolation, followed by a second application of the exponential) arises from \eqref{eq:integration_factor_second_order} combined with a similar interpretation as in the first-order scheme: half time step with the linear operator, evaluated at $t_n$ (\ie at the departure point, $L_{e_j}$), is applied to the solution computed at $t_n$; then the solution is obtained at $t_{n+1}$ through interpolation; finally, half time step with the linear operator at $t_{n+1}$ (\ie at the arrival point, $L_{j}$) is applied to the solution. 

\indent In the vectorial case 

\begin{equation}
    \label{eq:base_equation_slexp_linear}
    \frac{\xDif{\vecU}}{\xDif{t}} = \vecL \vecU
\end{equation}

\noindent with $\vecU \in \reals^N$, $\vecL \in \reals^{N \times N}$, we propose an analogous scheme, constructed directly from \eqref{eq:exp_linear_second_order}:

\begin{equation}
    \label{eq:exp_linear_second_order_vectorial}
    \vecU^{n+1} = \varphi_0\left(\frac{\Dt \vecL}{2}\right)\left[\varphi_0\left(\frac{\Dt \vecL}{2}\right)\vecU^n\right]^n_*
\end{equation}

\indent Note, however, that an analogous derivation is not necessarily valid. Indeed, we have, in general,

\begin{equation}
    \label{eq:non_commutation}
    e^{-\frac{\Dt}{2}\left[\vecL_{j} + \vecL_{e_j}\right]} \neq e^{-\frac{\Dt}{2}\vecL_{j}}e^{-\frac{\Dt}{2}\vecL_{e_j}} 
\end{equation}

\noindent except if $\vecL_j$ and $\vecL_{e_j}$ commute.
In Section \ref{subsec:LTE_second_order} we demonstrate the second-order accuracy of \eqref{eq:exp_linear_second_order}, in the scalar case, and of \eqref{eq:exp_linear_second_order_vectorial} in the vectorial case when $\vecL_j$ and $\vecL_{e_j}$ commute; the numerical simulations presented in Section \ref{sec:numerical_simulations} indicate that the second-order accuracy is preserved even if this commutation hypothesis does not hold.

\indent Applying \eqref{eq:exp_linear_second_order_vectorial} to the full PDE \xadded{\eqref{eq:base_equation_exp}}, we obtain the methods:

\begin{subequations}
    \begin{align}
    \vecU^{n+1}_{\text{\SE{2}{1}}} &=  \varphi_0\left(\frac{\Dt L}{2}\right) \left[ \varphi_0\left(\frac{\Dt L}{2}\right) \vecU^n \right]^n_* +   \varphi_0(\Dt L) \left[\Dt \psi_1(\Dt L) \tilde{\vecN}(\vecU^n) \right]_*^n, \label{eq:SL-ETD1RK-new}\\
    \vecU^{n+1}_{\text{\SE{2}{2}}} &=  \vecU^{n+1}_{\text{\SE{2}{1}}} + \Dt \varphi_0(\Dt L ) \left[\psi_2(\Dt L) \tilde{\vecN}(\vecU^{n+1}_{\text{\SE{2}{1}}}) - \left( \psi_2(\Dt L) \tilde{\vecN} (\vecU^n)        \right)^n_* \right] \label{eq:SL-ETD2RK-new}
    \end{align}
\end{subequations}

\noindent \xadded{where the Lagrangian trajectories are estimated as in the \SE{1}{1} and \SE{1}{2} schemes, using SETTLS (Eq. \eqref{eq:SETTLS}).}

\indent As proved in Section \ref{subsec:LTE_second_order} and/or verified through numerical simulations in Section \ref{sec:numerical_simulations}, \eqref{eq:SL-ETD2RK-new} provides a second-order approximation to \eqref{eq:base_equation_exp}. Therefore, the names of the methods were chosen following the same reasoning as before. In \eqref{eq:SL-ETD1RK-new}, a second- and a first-order discretizations of the linear and nonlinear terms, respectively, resulting in a global first-order approximation; in \eqref{eq:SL-ETD2RK-new}, a global second-order accuracy is ensured by second-order discretizations of both linear and nonlinear terms. 

\indent Note that the proposed discretization, using the approximation \eqref{eq:integration_factor_second_order} to the integration factor, is applied only to the linear term, \ie the nonlinear one is discretized as in \SE{1}{1} (eq. \eqref{eq:SL_ETD1RK}) and \SE{1}{2} (eq. \eqref{eq:SL_ETD2RK}). As verified in the numerical experiments in Section \ref{sec:numerical_simulations}, the discretization of the nonlinear term proposed by \cite{peixoto_schreiber:2019} in \SE{1}{2} has second-order accuracy, such that \SE{2}{2} achieves second-order with a relatively simple approximation of the exponential functions; the application of \eqref{eq:integration_factor_second_order} to the nonlinear term would lead to a much more complex formulation, since the integration factor acts both inside and outside the integral in \eqref{eq:exact_solution_slexp}.

\indent In this section, we prove analytically the second-order convergence of the proposed discretization considering only the linear term. In order to provide a more complete comparative study with the original formulation of the SL-ETDRK schemes, as well as with SL-SI-SETTLS and the Eulerian ETDRK methods, we also present an estimate of the computing complexity and, in Section \ref{subsec:serial_stability}, a detailed linear stability analysis. In Section \ref{sec:numerical_simulations}, we provide a numerical verification of these aspects considering \xadded{both a one-dimensional advection equation and, with more details,} the full nonlinear SWE on the rotating sphere.

\indent It is worth noting that linear and interpolation operators do not commute in general (see Appendix of \cite{peixoto_schreiber:2019}), such that, even with $\varphi_0(\Dt L /2) \varphi_0(\Dt L /2) = \varphi_0(\Dt L)$, we may have

\begin{equation}
    \label{eq:commutation}
    \varphi_0(\Dt L)[\vecU^n]^n_* \neq \varphi_0\left(\frac{\Dt L}{2}\right)\left[\varphi_0\left(\frac{\Dt L}{2}\right) \vecU^n\right]^n_* .
\end{equation}

\indent Finally, we highlight the similarities and differences of the proposed method \wrt the semi-Lagrangian Rosenbrock-type exponential schemes proposed by \cite{shashkin_goyman:2020}. In both methods, a splitting of the time step integration is performed in order to achieve second-order accuracy. In \cite{shashkin_goyman:2020}, this is made by evolving the fluid particles along the Lagrangian trajectory from $t_n$ to $t_{n+1/2} := t_n + \Dt/2$, applying the exponential operator at $t_{n+1/2}$, and evolving the fluid particles from $t_{n+1/2}$ to $t_n$. Here, we consider a different splitting, with two applications of the exponential operator and a single application of the semi-Lagrangian approach; also, both applications of the exponential operator are identical, corresponding to a left multiplication by the constant matrix $\varphi_0(\Dt L/2)$. Moreover, as described in Section \ref{subsec:ETDRK}, we consider a spectral discretization in space, which, with a proper rearrangement of the terms of the governing equations, allows for a straightforward computation of the matrix exponentials, while these are approximated in \cite{shashkin_goyman:2020} using Krylov-subspace methods.

\subsection{Truncation error analysis}
\label{subsec:LTE_second_order}

\indent The \replaced[id=R2-new]{truncation error}{consistency} analysis of the proposed scheme is analogous to the first order scheme. We begin by considering the scalar case, \ie \eqref{eq:exp_linear_second_order} applied to \eqref{eq:exp_linear}. It leads to the system of uncoupled ODEs

\begin{equation}
    \label{eq:exp_linear_j_second_order}
    U^{n+1}_j = \varphi_0\left(\frac{\Dt L_j}{2}\right) \left[\varphi_0\left(\frac{\Dt L}{2}\right)U^n\right]_{*_j}, \qquad j = 0, \dots, P
\end{equation}

\indent We begin by inserting the exact solution $u$ into the $j-$th equation of \eqref{eq:exp_linear_j_second_order} and using the error estimate \eqref{eq:error_cubic_interpolation} for SETTLS using cubic Lagrange interpolation:

\begin{equation}
    \label{eq:second_order_exact}
    \begin{aligned}
        u^{n+1}_j &= \varphi_0\left(\frac{\Dt L_j}{2}\right)\left[ \varphi_0\left(\frac{\Dt L}{2}\right) u^n\right]^n_{*_j}\\
        &= \varphi_0\left(\frac{\Dt L_j}{2}\right) \varphi_0\left(\frac{\Dt L_{e_j}}{2}\right) u^n_{e_j} + \mc{O}(\Dt^2) + \mc{O}(\Dt\Dx^2) + \mc{O}(\Dx^4)
    \end{aligned}
\end{equation}

\indent Then, using the exact solution \eqref{eq:exact_solution_linear_j} and the expression \eqref{eq:integration_factor_second_order} for the integration factor along the Lagrangian trajectory $\Gamma^{n+1}_j$, we get the local truncation error

\begin{equation}
    \label{eq:LTE_second_order}
    \begin{aligned}
        \tau^{n+1}_j & = \frac{1}{\Dt} \left[ u^{n+1}_j - \varphi_0\left(\frac{\Dt L_j}{2}\right) \left(\frac{\Dt L_{e_j}}{2}\right) u^n_{e_j} + \mc{O}(\Dt^2) + \mc{O}(\Dt\Dx^2) + \mc{O}(\Dx^4)\right]\\
        &= \frac{1}{\Dt} \left[e^{C_{2,j}\Dt^3} e^{\frac{\Dt L_{e_j}}{2}} e^{\frac{\Dt L_{j}}{2}}  u^n_{e_j} - e^{\frac{\Dt L_{j}}{2}} e^{\frac{\Dt L_{e_j}}{2}} u^{n}_{e_j} + \mc{O}(\Dt^3) + \mc{O}(\Dt\Dx^2) + \mc{O}(\Dx^4) \right]\\
        &= \frac{1}{\Dt}\left[(C_{2,j}\Dt^3 + \mc{O}(\Dt^6))e^{\frac{\Dt L_{e_j}}{2}} e^{\frac{\Dt L_{j}}{2}}  u^n_{e_j} + \mc{O}(\Dt^3) + \mc{O}(\Dt\Dx^2) + \mc{O}(\Dx^4) \right]\\
        & = \mc{O}(\Dt^2) + \mc{O}(\Dx^2) + \mc{O}\left(\frac{\Dx^4}{\Dt}\right)
    \end{aligned}
\end{equation}

\noindent confirming that the proposed method is second-order accurate.

\indent In the vectorial case (Eq. \eqref{eq:exp_linear_second_order_vectorial} applied to \eqref{eq:base_equation_slexp_linear}), this result is also true if $\vecL_{j}$ and $\vecL_{e_j}$ commute. Indeed, in this case, the last equality in \eqref{eq:integration_factor_second_order} would hold and $\varphi_0(\Dt \vecL_j /2)$ and $\varphi_0(\Dt \vecL_{e_j} /2)$ would commute, such that the third equality in the truncation error analysis \eqref{eq:LTE_second_order} would be true. In the general case of a non-commutative linear operator $\vecL$ along each Lagrangian trajectory, a more careful analysis would be required. For non-commuting square matrices $A,B$ of same size, it holds

\begin{equation}
    \label{eq:non_commuting}
    e^{A+B} - e^{A}e^{B} = \frac{1}{2}\left[BA - AB\right] + \text{ third-order terms} 
\end{equation}

\noindent such that, with $A = -\Dt \vecL_j/2$ and $B = -\Dt \vecL_{e_j}/2$, the proposed approximation \eqref{eq:integration_factor_second_order} to the integration factor would have a leading error term of order $\mc{O}(\Dt^2)$, which would contribute with $\mc{O}(\Dt)$ to the local truncation error, thus ensuring only first-order accuracy. Note, however, that if $\vecL$ is smooth enough, we can expand $\vecL_j$ around $\vecL_{e_j}$ along $\Gamma^{n+1}_j$ to get

\begin{equation}
    \label{eq:non_commuting_expansion2}
    \begin{aligned}
        \vecL_{e_j}\vecL_j - \vecL_j\vecL_{e_j} && = & \ \vecL_{e_j}\left(\vecL_{e_j} + \Dt\dert{\vecL}|_{e_j} + \mc{O}(\Dt^2)\right) - \left(\vecL_{e_j} + \Dt\dert{\vecL}|_{e_j} + \mc{O}(\Dt^2)\right)\vecL_{e_j}\\
        && = & \Dt \left( \vecL_{e_j}\dert{\vecL}|_{e_j} - \dert{\vecL}|_{e_j} \vecL_{e_j}\right) + \mc{O}(\Dt^2) = \mc{O}(\Dt)
    \end{aligned}
\end{equation}

\noindent such that

\begin{equation}
    \label{eq:non_commuting_expansion3}
    \begin{aligned}
    e^{-\frac{\Dt}{2}\left[\vecL_j + \vecL_{e_j}\right]} - e^{-\frac{\Dt}{2}\vecL_j}e^{-\frac{\Dt}{2}\vecL_{e_j}} & = \frac{\Dt^2}{8}\left[\vecL_{e_j}\vecL_j - \vecL_j\vecL_{e_j}\right] + \mc{O}(\Dt^3) = \mc{O}(\Dt^3)
    \end{aligned}
\end{equation}

\indent This implies that, if the integration factor can be written in the matrix exponential form \eqref{eq:integration_factor_second_order}, then it actually satisfies

\begin{equation}
    \label{eq:integration_factor_second_order_noncommut}
    P_{n_j}(t_{n+1}) = e^{-\int_{t_n}^{t_{n+1}} \vecL_{x^j(s)}ds} = e^{-\frac{\Dt}{2}\left[\vecL_j + \vecL_{e_j}\right] - C_{2,j}\Dt^3} = \left(e^{-\frac{\Dt \vecL_j}{2}} e^{-\frac{\Dt \vecL_{e_j}}{2}} + C_j\Dt^3\right) e^{-C_{2,j}\Dt^3}
\end{equation}

\noindent for a given constant $C_j$ accounting for the non-commutation of $\vecL_j$ and $\vecL_{e_j}$. In other words, the splitting of the exponential would introduce an $\mc{O}(\Dt^3)$ term that would contribute with $\mc{O}(\Dt^2)$ to the local truncation error, thus preserving the second-order accuracy. Nevertheless, it is important to notice that, in the case where $\vecL_{\vecx(t)}$ and $\vecL_{\vecx(s)}$ do not commute, we would have to additionally assume that the integration factor exists and can be written in the matrix exponential form \eqref{eq:integration_factor}, otherwise, the development of the local truncation error \eqref{eq:LTE_second_order} would not be valid, since we would not be able to write an explicit form to the exact solution $u^{n+1}_j$ in \eqref{eq:exact_solution_linear_j}. However, the numerical simulations presented in Section \ref{sec:numerical_simulations} indicate that second-order accuracy is achieved by \eqref{eq:exp_linear_second_order_vectorial} even if this commutation is not asserted. Finally, note that commutation-related problems are known issues in semi-Lagrangian schemes (\eg the non-commutation between Laplace transforms and linear differential operators along Lagrangian trajectories, as studied in \cite{lynch:2022}).

% \indent As a last comment on the accuracy of semi-Lagrangian exponential schemes, the truncation error analysis presented above suggests some possible ways to derive higher-order discretizations of the linear term, by considering higher-order approximations of the integral, \eg the Simpson's rule; also, an alternative second-order scheme could possibly be obtained by using the midpoint rule. In both cases, additional difficulties concerning the evaluation of $\vecL_{\vecx(s)}$ at intermediate points of the Lagrangian trajectories would need to be dealt with. However, the $\mc{O}(\Dt^2)$ error due to the non-commutation of the linear operator evaluated at different points of the Lagrangian trajectory could prevent higher-order accuracy, and more careful studies should be conducted.

% \input{alternative_SLETD2RK/truncation_error}

\subsection{Computational complexity}

\indent We now provide an estimate of the computational complexity of the proposed modified semi-Lagrangian exponential methods and the other time integration schemes considered here. In each time step, the following computations are performed:

\begin{itemize}
    \item In the case of exponential methods, the evaluation of the $\varphi_k$ and/or the $\psi_k$ functions. Since these functions are defined recursively, their cost increases for larger $k$. These functions need to be computed for each pair of wavenumber $(m,n)$, and we recall that the matrix exponentials are computed in this work considering the diagonalization of the matrix (eq. \eqref{eq:exp_diag}), such that we need only to compute the exponential of the eigenvalues of $\vecL = \vecL_{\vecG}$ (two per wavenumber pair);
    \item In the case of semi-Lagrangian schemes, the trajectories need to be estimated, and spatial interpolations to the departure points need to be performed;
    \item In the case of SL-SI-SETTLS \eqref{eq:SL_SI_SETTLS}, the semi-implicit discretization for the linear term requires both the application of $\vecL$ and $\vecL^{-1}$, the latter through the solution of a Helmholtz equation \cite{temperton:1997};
    \item Finally, all schemes require the evaluation of at least one nonlinear term. \xadded{We take into account reutilizations of the computed terms between consecutive time steps.}
\end{itemize}

\indent Table \ref{tab:complexity} summarizes the number of each operation listed above performed by each time integration scheme. In this table, each $\varphi_k / \psi_k$ function is depicted individually; $\vecx_d$ and $()_*$ denote respectively the evaluation of the Lagrangian trajectories and the interpolation to the departure points; $\vecL$ and $\vecL^{-1}$ correspond to the semi-implicit scheme described above; and $\vecN_{\vecA}$ and $\vecN_{\vecR}$ indicate the evaluation of the (respectively advection and remaining) nonlinear terms.

\begin{table}[!htbp]
    \centering
    \begin{tabular}{|c|c|c|c|c|c||c|c||c|c||c|c|}
        \hline
         Scheme & $\varphi_0$ & $\varphi_1$ & $\varphi_2$ & $\psi_1$ & $\psi_2$ & $\vecx_d$ & $()_*$ &  $\vecL$ & $\vecL^{-1}$ & $\vecN_{\vecA}$ & $\vecN_{\vecR}$  \\
         \hline
         ETD1RK & 1 & 1 & 0 & 0 & 0 & 0 & 0 & 0 & 0  & 1 & 1 \\
         \hline
         ETD2RK & 1 & 1 & 1 & 0 & 0 & 0 & 0 & 0 & 0 & 2 & 2 \\
         \hline 
         \SE{1}{1} & 1 & 0 & 0 & 1 & 0 & 1 & 1 & 0 & 0  & 0  & 1 \\
         \hline
         \SE{1}{2} & 2 & 0 & 0 & 1 & 2 & 1 & 2 & 0 & 0 & 0 & 2 \\
         \hline
         \SE{2}{1} & 3 & 0 & 0 & 1 & 0 & 1 & 2 & 0 & 0 & 0 & 1 \\
         \hline
         \SE{2}{2} & 4 & 0 & 0 & 1 & 2 & 1 & 3 & 0 & 0 &  0 & 2 \\
         \hline
         SL-SI-SETTLS & 0 & 0 & 0 & 0 & 0 & 1 & 2 & 1 & 1 & 0 & \xreplaced{1}{2} \\
         \hline
    \end{tabular}
    \caption{Number of operations performed per time step by each time integration scheme. $\vecx_d$ and $()_*$ denote respectively the evaluation of the Lagrangian trajectories and the interpolation to the departure points; $\vecL$ and $\vecL^{-1}$ correspond to the semi-implicit scheme.}
    \label{tab:complexity}
\end{table}

\indent It is clear from this table that the proposed modifications of the semi-Lagrangian exponential schemes introduce additional complexity per time step, with two more exponential evaluations and one more spatial interpolation \wrt their respective original formulations. However, as assessed in the numerical simulations presented in Section \ref{sec:numerical_simulations}, the improved accuracy and stability of \SE{2}{2} allow us to choose larger time step sizes. Moreover, it has been identified in the numerical simulations that the computation of the exponential functions is responsible for an important fraction of the computing time in our implementation, such that important savings are possible by precomputing them; additional savings could also be obtained with a more efficient implementation (\eg by using vectorization).

\section{Linear stability analysis}
\label{subsec:serial_stability}

\indent In the following paragraphs, we conduct a stability study comparing exponential and semi-Lagrangian exponential schemes, as well as SL-SI-SETTLS. Linear stability analyses of Eulerian exponential schemes have been conducted, \eg by \cite{cox_matthews:2002, crouseilles_al:2020, buvoli_minion:2022}. The last two cited works show notably that these methods are in general unstable when applied to non-diffusive problems, but depending on configurations of the problem and the exponential scheme, the amplifications are small such that instabilities may remain unnoticed within the simulation time. In this section, we focus on comparing stability regions of both Eulerian and semi-Lagrangian schemes.

\indent The stability study proposed here is inspired by the procedure presented by \cite{cox_matthews:2002}. Consider a generic nonlinear ODE in the form 

\begin{equation}
    \label{eq:base_ODE}
    \dot{u} = \lambdaL u + \vecN(u)
\end{equation}

\indent Substituting $u = \ol{u} + u'$ in \eqref{eq:base_ODE}, where $\ol{u}$ is a steady state (\ie $\lambdaL \ol{u} + \vecN(\ol{u}) = 0$) and $u'$ is a perturbation, expanding $N(\ol{u} + u')$ around $\ol{u}$, truncating at first order and defining $\lambdaN := \vecN'(\ol{u})$, we obtain the linearized equation

\begin{equation}
        \label{eq:base_ODE_linearized}
        \dot{u}' = \lambdaL u' + \lambdaN u'
\end{equation}

\noindent which we consider hereafter, with the primes omitted for conciseness. However, we need to take into account that we are comparing both Eulerian (ETDRK) and semi-Lagrangian (SL-SI-SETTLS and SE**) methods, and this study should be able to capture the influence of advective effects and interpolation procedures. In the case of SL schemes, the advection is incorporated into the material derivative, and we can directly use the same linearized ODE as above:

\begin{equation}
    \label{eq:pde_SL}
    \dert{u} =  \lambdaL u + \lambdaN u
\end{equation} 

\noindent where $u, \lambdaL, \lambdaN \in \complex$. This equation arises, \eg when \eqref{eq:base_equation_exp} is solved using a spectral method, in which case $u = \hat{u}_k(t)$ is a spectral coefficient of the solution, and $\lambdaL$, $\lambdaN$ are modes associated respectively to $\vecL$ and the linearized term $\vecN$. If $\vecL$ and $\vecN$ are diagonalizable in the spectral space with the same set of eigenvectors, then \eqref{eq:base_equation_exp} leads to a system of uncoupled ODEs of the form \eqref{eq:pde_SL}, each one corresponding to a wavenumber $k$.

\indent On the other hand, in the case of Eulerian schemes, we need to consider the advection explicitly, leading to a linearized one-dimensional scalar advection equation with source term:

\begin{equation}
    \label{eq:pde_eulerian}
    \pdert{u} + v\pderx{u} = \lambdaL u + \vecN(u) \implies \pdert{u} + v\pderx{u} =  \lambdaL u + \lambdaN u
\end{equation}

\noindent in which the spatial derivative becomes a scalar multiplication when solved in the spectral space, and $v$ is a scalar velocity field assumed to be constant for simplicity.

\indent The stability region of a given numerical scheme with stability function \replaced[id=R2-new]{$A := \left|u^{n+1}\right|/\left|u^n\right|$}{$A := |u^{n+1}/u^n|$} is defined as the region in $\complex^2$ in which $|A(\Dt\lambdaL, \Dt\lambdaN)| \leq 1$, where $u^n = u(t_n, x) = A^n e^{i \kappa x}$ is the computed approximation to $u(t_n)$, $\kappa = 2 \pi k / L$ is the normalized wavenumber and $L$ is the length of the spatial domain. A complete stability study should be performed in a four-dimensional space, since $A$ depends on the real and imaginary parts of $\lambdaL$ and $\lambdaN$; however, since a hyperbolic problem such as the SWE propagates purely hyperbolic modes, we simplify this study by \xadded{considering both $\xiL := \Dt \lambdaL$ and $\xiN:= \Dt \lambdaN$ to be purely imaginary (Figure \ref{fig:stability_schemes_2D}), and, in particular, by }fixing values of $\xiL \in i\reals$ \xdeleted{(since the eigenvalues of $\vecL$ are purely imaginary)} and plotting the amplification factors as functions of $\Im(\xiN)$ \xadded{(Figure \ref{fig:stability_schemes_1D}). We present below the stability functions of each time integration scheme considered in this work}.

\subsection{Stability of ETD1RK and ETD2RK}

\indent When applying \eqref{eq:ETD1RK} and \eqref{eq:ETD2RK} to \eqref{eq:pde_eulerian}, the advection term could be treated either as a linear or nonlinear term, yielding different stability functions. We consider only the latter case, but the overall conclusions on the stability of ETD1RK and ETD2RK compared to the other schemes considered in this work are similar. We have

\begin{subequations}
    \begin{align}
        u^{n+1}_{\text{ETD1RK}} & = \varphi_0(\xiL) u^n + \Dt  \varphi_1(\xiL) \lambdaN u^n - \Dt \varphi_1(\xiL)  i v \kappa u^n  \label{eq:ETD1RK_ode} \\
        & = \left[ \varphi_0(\xiL) + \varphi_1(\xiL) (\xiN - i \kappa s)\right] u^n  \notag \\
        u^{n+1}_{\text{EDT2RK}} & = u^{n+1}_{\text{ETD1RK}} + \Dt \varphi_2(\xiL)\left[ (\lambdaN - iv\kappa) u^{n+1}_{\text{ETD1RK}} - (\lambdaN - v i \kappa) u^n \right] \label{eq:ETD2RK_ode} \\
        & = u^{n+1}_{\text{ETD1RK}} +  \varphi_2(\xiL) \left( \xiN - i \kappa s \right) \left(u^{n+1}_{\text{ETD1RK}} - u^n \right) \notag
    \end{align}
\end{subequations}

\noindent where $s := v \Dt = x_j - x_{d_j}$ is the spatial displacement of each particle from $t_n$ to $t_{n+1}$ (see Figure \ref{fig:lagrangian_trajectory}). Note that we have considered the spatial derivative to be computed exactly (which can be achieved, \eg by using a spectral method, as considered in this work), \ie $(u^n)_x = i \kappa u^n$, since we want to analyze exclusively the stability properties due to the temporal discretization.

\indent Thus, the stability functions of the Eulerian ETD1RK and ETD2RK schemes read

\begin{subequations}
    \begin{align}
    A_{\text{ETD1RK}} & = \varphi_0(\xiL) + \varphi_1(\xiL) (\xiN - i \kappa s) \label{eq:A_ETD1RK} \\
     A_{\text{ETD2RK}} & =  A_{\text{ETD1RK}} + \varphi_2(\xiL)\left(\xiN - i \kappa s \right) \left[A_{\text{ETD1RK}} - 1 \right]  \label{eq:A_ETD2RK}
    \end{align}
\end{subequations}

\indent Note that the presence of the term $-i \kappa s$ makes the stability regions dependent on the spatial information. This implies that the schemes may not be simultaneously stable for a large range of wavenumbers, since the intersection of their stability regions may be small or even empty. Following \cite{hortal:2002}, we compute stability functions $A_{\kappa s} = A_{\kappa s}(\xiL, \xiN)$ for $\kappa s \in \mc{K} := \{0, \pi/10, 2\pi/10, \dots, 2\pi\}$; then, for each $(\xiL, \xiN)$, the amplification factor is $A(\xiL, \xiN) = \max_{\kappa s \in \mc{K}} A_{\kappa s}(\xiN, \xiL)$.

\subsection{Stability of SL-SI-SETTLS}

\indent We now develop a linear stability analysis for the semi-implicit semi-Lagrangian SETTLS \cite{hortal:2002} described in Section \ref{subsec:settls}. In the stability analysis of semi-Lagrangian schemes, we also need to take into account the spatial dependence of the solution. Using $u(t_n, x) = A^n e^{i\kappa x}$ in \eqref{eq:SL_SI_SETTLS} applied to \eqref{eq:pde_SL}, we obtain 

\begin{equation}
    \label{eq:SETTLS_ode}
    \added[id=R2-new]
    {A^{n+1}e^{i\kappa x_j} - A^{n}e^{i\kappa x_{d_j}} = \frac{\xiL}{2}\left( A^{n+1}e^{i\kappa x_j} + A^{n}e^{i\kappa x_{d_j}} \right)  + \frac{\xiN}{2} \left( \left[2 A^{n} - A^{n-1} \right]e^{i\kappa x_{d_j}} + A^n e^{i\kappa x_j}\right)}
\end{equation}

\noindent which, after division by $A^{n-1} e^{i \kappa x_j}$ and recalling that $s: = x_j - x_{d_j}$, leads to the quadratic equation

\begin{equation}
    \label{eq:stability_SL_SI_SETTLS}
    A^2 \left( 1 - \frac{\xiL}{2} \right) - A \left( e^{-i\kappa s} \left(1 +  \frac{\xiL}{2} + \xiN \right) + \frac{\xiN}{2} \right) + \frac{\xiN}{2}e^{-i \kappa s} = 0
\end{equation}

\noindent whose roots $A  = A_{\text{SL-SI-SETTLS}}$ define the stability function of SL-SI-SETTLS; its stability region is the intersection between the stability regions defined by each of the two roots of \eqref{eq:stability_SL_SI_SETTLS}. Moreover, as in the case of ETD1RK and ETD2RK considering the advection term, these roots depend on the spatial wavenumber through the term $e^{-i \kappa s}$, and the stability region is defined as the intersection of the individual stability regions for each $\kappa s$.

\subsection{Stability of \SE{1}{1} and \SE{1}{2}}
\label{subsec:stability_SE11_SE12}

\indent We now study the stability of the semi-Lagrangian exponential schemes proposed by \cite{peixoto_schreiber:2019} and described in Section \ref{subsec:sletdrk}. Using $u(x,t_n) = A^n e^{i\kappa x}$ in \eqref{eq:SL_ETD1RK} and \eqref{eq:SL_ETD2RK} applied to \eqref{eq:pde_SL} and the fact that $\varphi_k(z) = \varphi_0(z) \psi_k(z) \ \forall k > 0$, we get

\begin{subequations}
    \begin{align}
    A_{\text{\SE{1}{1}}} & = \varphi_0(\xiL) e^{-i\kappa s} + \xiN \varphi_1(\xiL) e^{-i \kappa s} = e^{-i \kappa s} \varphi_0(\xiL) [ 1 + \xiN \psi_1(\xiL)] \label{eq:A_SL_ETD1RK} \\
    A_{\text{\SE{1}{2}}} & =  A_{\text{\SE{1}{1}}} + \xiN \varphi_2(\xiL) \left[A_{\text{\SE{1}{1}}} - e^{-i\kappa s} \right] \label{eq:A_SL_ETD2RK}\\
    & = e^{-i\kappa s} \left[ \left(1 + \xiN \varphi_2(\xiL)\right) \varphi_0(\xiL) \left( 1 + \xiN \psi_1(\xiL)\right) - \xiN \varphi_2(\xiL)\right] \notag
    \end{align}
\end{subequations}

\indent Since $|e^{-iks}| = 1$, the stability regions of both \SE{1}{1} and \SE{1}{2} do not depend on $\kappa s$. This is a remarkable feature of these schemes, since, contrary to what happens to SL-SI-SETTLS and ETDRK, there is no reduction of the global stability region when intersecting the regions corresponding to each wavenumber.

\begin{figure}[!htbp]
    \begin{subfigure}{.3\linewidth}
        \centering
        \includegraphics[clip=true,trim=0 0 35 0, scale=.375]{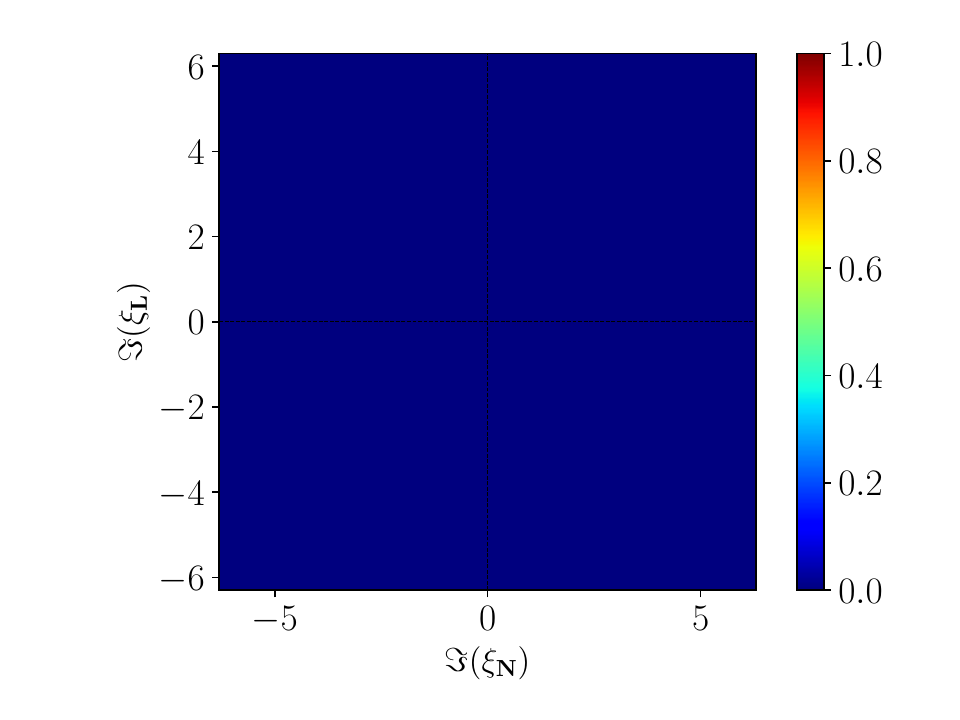}
        \caption{ETD1RK}
    \end{subfigure}
    \begin{subfigure}{.3\linewidth}
        \centering
        \includegraphics[clip=true,trim=0 0 35 0, scale=.375]{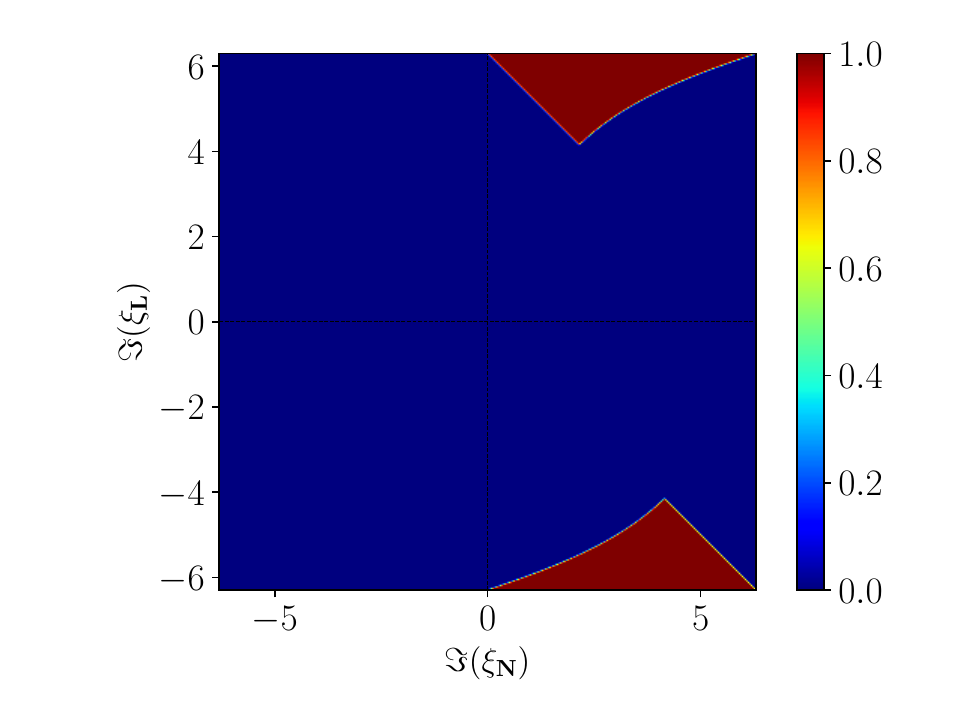}
        \caption{ETD2RK}
    \end{subfigure}
    \begin{subfigure}{.3\linewidth}
        \centering
        \includegraphics[clip=true,trim=0 0 35 0, scale=.375]{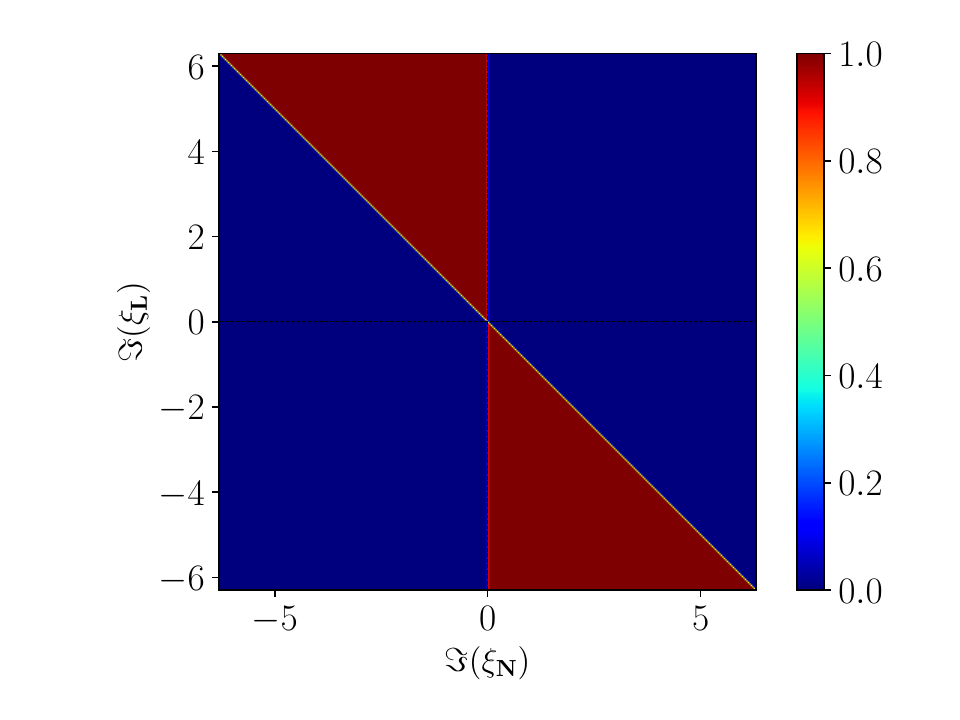}
        \caption{SL-SI-SETTLS}
    \end{subfigure}
    \begin{subfigure}{.3\linewidth}
        \centering
        \includegraphics[clip=true,trim=0 0 35 0, scale=.375]{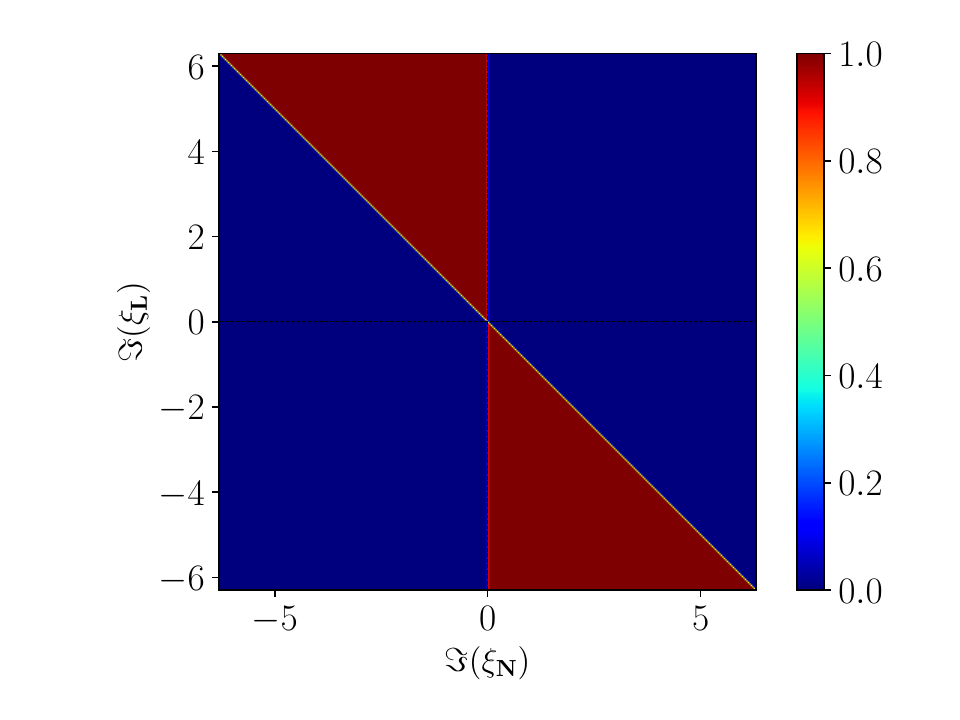}
        \caption{\SE{1}{1}}
    \end{subfigure}
    \begin{subfigure}{.3\linewidth}
        \centering
        \includegraphics[clip=true,trim=0 0 35 0, scale=.375]{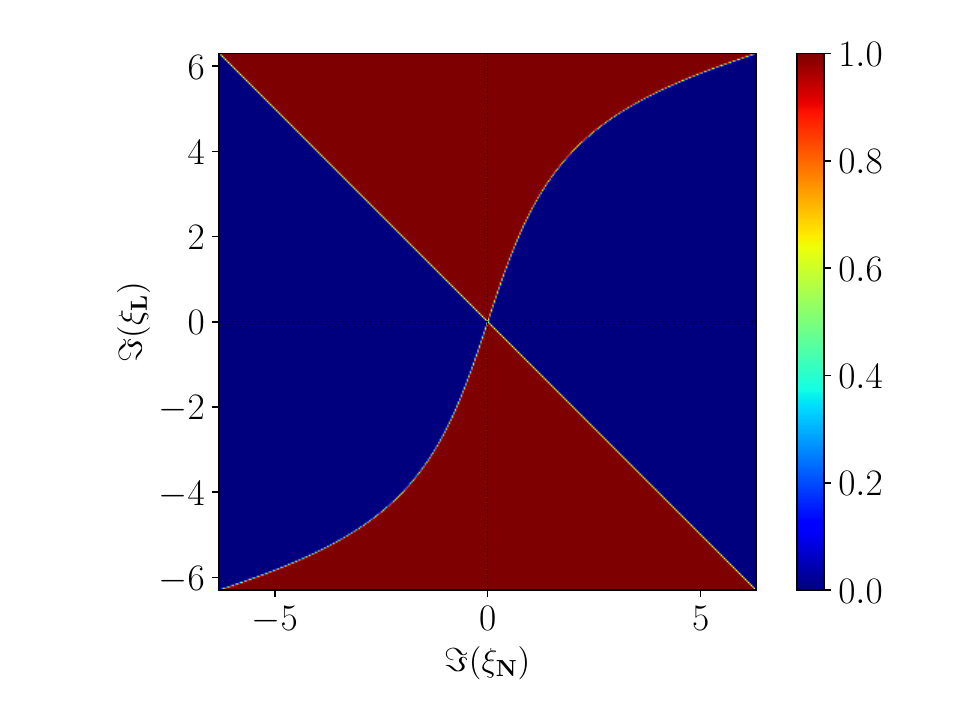}
        \caption{\SE{1}{2}}
    \end{subfigure}
    \caption{\xadded{Stability regions (in red) of ETD1RK and ETD2RK applied to \eqref{eq:pde_eulerian} (considering a nonlinear treatment of the advection term), and SL-SI-SETTLS, \SE{1}{1} and \SE{1}{2} applied to \eqref{eq:pde_SL}, in the \replaced[id=R2-new]{($\Im(\xiN),\Im(\xiL))$}{$\Im(\xiN)-\Im(\xiL)$} plane, considering both $\xiN$ and $\xiL$ to be purely imaginary.}}
    \label{fig:stability_schemes_2D}
\end{figure}

\begin{figure}[!htbp]
    \begin{subfigure}{.475\linewidth}
        \centering
        \includegraphics[clip=true,trim=0 0 0 0, scale=.425]{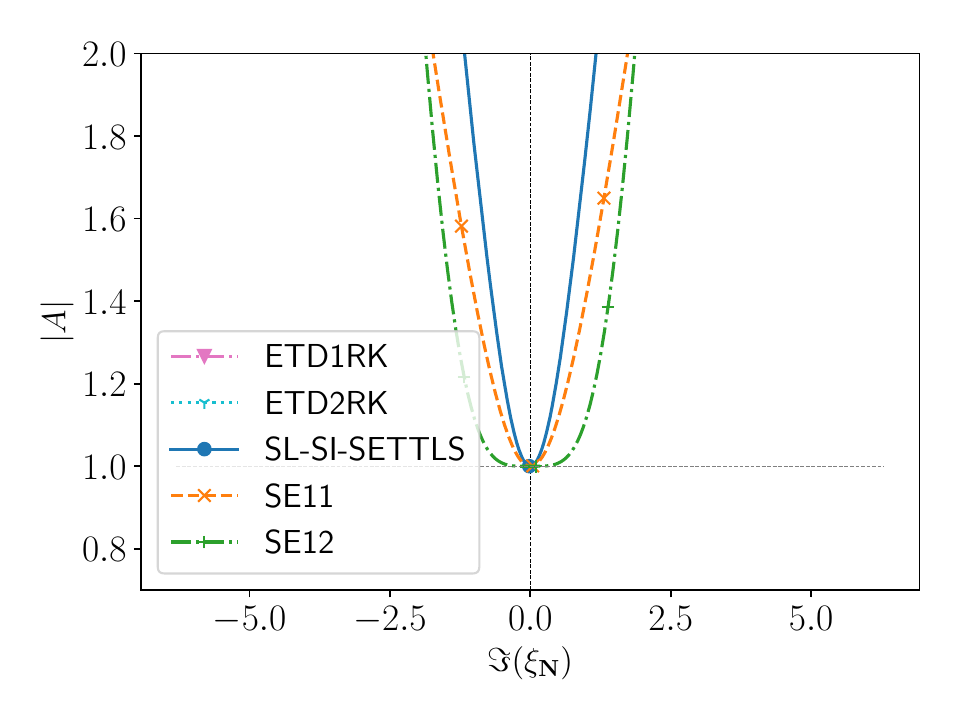}
        \caption{$\xiL = 0$}
    \end{subfigure}
    \begin{subfigure}{.475\linewidth}
        \centering
        \includegraphics[clip=true,trim=0 0 0 0, scale=.425]{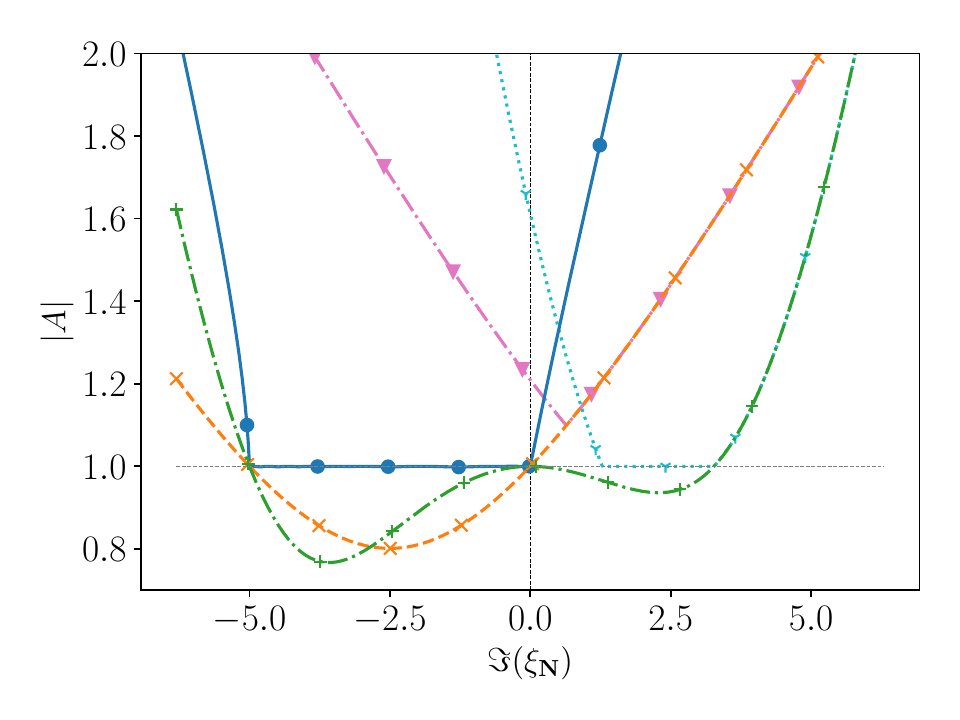}
        \caption{$\xiL = 20000\tilde{\xi}_{\vecL}$}
    \end{subfigure}
    \caption{Stability functions of ETD1RK and ETD2RK applied to \eqref{eq:pde_eulerian} (considering a nonlinear treatment of the advection term), and SL-SI-SETTLS, \SE{1}{1} and \SE{1}{2} applied to \eqref{eq:pde_SL}, as a function of $\Im(\xiN)$ for fixed values of $\xiL$ and $\Re(\xiN) = 0$. The stability regions of ETD1RK and ETD2RK are empty in the case $\xiL = 0$. The horizontal, dashed line indicates the boundary of the stability region requiring $|A| \leq 1$. The same legend holds for both plots.}
    \label{fig:stability_schemes_1D}
\end{figure}

\indent We compare in Figure \ref{fig:stability_schemes_2D} the stability functions of ETD1RK and ETD2RK applied to \eqref{eq:pde_eulerian}, and SL-SI-SETTLS, \SE{1}{1} and \SE{1}{2} applied to \eqref{eq:pde_SL}, \xadded{in the plane \replaced[id=R2-new]{$(\Im(\xiN),\Im(\xiL))$}{$\Im(\xiN)-\Im(\xiN)$}. In Figure \ref{fig:stability_schemes_1D}, we present the amplification factor along two horizontal slices of these planes, with $\xiL$ chosen as} \xdeleted{We present the regions for} integer multiples of $\tilde{\xi}_{\vecL} = i \sqrt{\ol{\Phi}/a^2} \approx 2.5 \times 10^{-4}i$, which is obtained from an approximation to the eigenvalues $\lambda_{\vecL} = i\sqrt{\ol{\Phi}n(n+1)/a^2} $ of $\vecL = \vecL_{\vecG}$ considering relevant physical parameters and an f-plane approximation (see \cite{caldas_al:2023} for details). The scaling of $\tilde{\xi}_{\vecL}$ can be seen as a result of different choices of time step size $\Dt$ and/or wavenumbers: for instance, the plot considering $\xiL = 20000 \tilde{\xi}_{\vecL}$ in Figure \ref{fig:stability_schemes_1D} corresponds to the stability region of the maximum wavenumber in a spectral resolution $M=256$ and a reasonable time step size of approximately $\Dt = 80\text{s}$. The plots reveal notably an enhancement of the stability of the semi-Lagrangian exponential schemes when compared to their Eulerian counterparts, which suffer from a small intersection of the stability regions of each $\kappa s$, indicating that the schemes have a lack of stability common to a large range of combinations between wavenumbers and spatial displacement. The semi-Lagrangian exponential scheme with second-order discretization of the nonlinear term (\SE{1}{2}) also presents a considerably larger stability region compared to SL-SI-SETTLS \xadded{and \SE{1}{1}, whose stability regions coincide, despite of having different stability functions}.

\subsection{Stability of \SE{2}{1} and \SE{2}{2}}
\label{subsec:stability_discussion}

\indent Using the same procedure described above to assess stability properties of the proposed schemes \SE{2}{1} and \SE{2}{2} (Section \ref{sec:second_order_SLETDRK}) does not provide insights on potential stability improvements or decreases induced by the splitting of the matrix exponential. In particular, this approach is not able to capture the fact that the matrix exponential and the interpolation do not commute \cite{peixoto_schreiber:2019}. Indeed, using $u(x,t_n) = A^n e^{i\kappa x}$ in \eqref{eq:SL-ETD1RK-new} and \eqref{eq:pde_SL}, analogously to \eqref{eq:A_SL_ETD1RK}, provides

\begin{equation}
    \label{eq:stability_SE21_SE22}
    \varphi_0\left(\frac{\xiL}{2}\right) e^{-i \kappa s} \varphi_0\left(\frac{\xiL}{2}\right) = e^{\xiL/2} e^{-i \kappa s} e^{\xiL/2} = e^{\xiL} e^{-i \kappa s} = \varphi_0(\xiL) e^{-i \kappa s}
\end{equation}

\noindent such that $A_{\text{\SE{2}{1}}} \equiv A_{\text{\SE{1}{1}}}$, and, similarly, $A_{\text{\SE{2}{2}}} \equiv A_{\text{\SE{1}{2}}}$. However, different stability properties are verified in the numerical simulations presented in Section \ref{sec:numerical_simulations}.

\indent Next, we develop in detail how the differences between the two sets of schemes arise. Consider the ODE \eqref{eq:pde_SL} with only the linear term:

\begin{equation}
    \label{eq:ODE_linear}
    \dert{u} = \lambdaL u
\end{equation}

\indent We first consider the first-order discretization used in \eqref{eq:SL_ETD1RK}. Since the problem is solved in the spectral space, we may write

\begin{equation}
    \label{eq:spectral_SLETD1RK_linear}
    \hat{u}_k^{n+1} = \varphi_0({\xiLk}) [\hat{u}_k^n]_*, \qquad k = -M, \dots, M
\end{equation}

\noindent where $u(t,x) \approx \sum_{k=-M}^M \hat{u}_k^{n} e^{i \kappa x}$ is the truncated Fourier decomposition of $u$ and $ \xiLk := \Dt{\lambdaLk}$ may depend on $k$. Note, however, that the interpolation is performed on the grid. Therefore, using the direct and inverse discrete Fourier transforms

\begin{equation}
    \label{eq:fourier_transform}
    \hat{u}_k^n = \mc{F}(u^n) = \frac{1}{P} \sum_{j = 0}^{P-1} u(t_n, x_j) e^{-i \kappa x_j}, \qquad u(t_n, x_j) = \mc{F}^{-1}(\hat{u}^n) = \sum_{k = -M}^{M} \hat{u}^n_{k} e^{i \kappa x_j}
\end{equation}

\noindent with $\kappa = 2 \pi k /L$, L the length of the periodic spatial domain and $P+1$ the number of equidistant grid points $x_i = i\Dx$, $i = 0,\dots,P$, $\Dx = L/P$, we may rewrite \eqref{eq:spectral_SLETD1RK_linear} more rigorously as

\begingroup
\allowdisplaybreaks
% \begin{equation}
    \begin{align}
    \hat{u}^{n+1}_k & = \varphi_0(\xiLk) \mc{F} ( [ \mc{F}^{-1} (\hat{u}^n) ]_*) =  \varphi_0(\xiLk) \mc{F} ( [ u^n_j]_*) = \varphi_0(\xiLk) \frac{1}{P}\sum_{j = 0}^{P-1} u(t_n, x_{d_j}) e^{-i \kappa x_j} \notag \\
    & = \varphi_0(\xiLk) \frac{1}{P}\sum_{j = 0}^{P-1} u(t_n, x_{j - p_j}) e^{-i \kappa x_j} = \varphi_0(\xiLk) \frac{1}{P}\sum_{j = 0}^{P-1} u(t_n, x_{r_j}) e^{-i \kappa x_{p_j+r_j}} \notag \\
    & = \varphi_0(\xiLk) \frac{1}{P}\sum_{j = 0}^{P-1} u(t_n, x_{r_j}) e^{-i \kappa (p_j + r_j) \Dx}, \label{eq:A_SLETD1RK_linear}
    \end{align}
% \end{equation}
\endgroup

\noindent where $x_{d_j}$ is the approximate departure point at $t_n$ of the Lagrangian trajectory arriving at $x_j$ at $t_{n+1}$\added[id=R2-new]{, $r_j := j-p_j$}, and, for simplification, $x_{d_j} = x_{j - p_j}$, $p_j \in \integers$, is supposed to belong to the grid (otherwise $u(t_n,x_{d_j})$ may be written as a linear combination (\ie interpolation) of the solution computed at grid points). If the velocity field is constant in space, then $\forall j, p_j = p$, thus

\begin{equation}
    \label{eq:A_SE21_constant_velocity}
    \begin{aligned}
    \hat{u}^{n+1}_k & = \varphi_0(\xiLk) e^{-i\kappa p \Dx} \frac{1}{P}\sum_{j = 0}^{P-1} u(t_n, x_{r_j}) e^{-i \kappa x_{r_j}} \\
    & = \varphi_0(\xiLk) e^{-i\kappa p \Dx} \frac{1}{P}\sum_{r = 0}^{P-1} u(t_n, x_{r}) e^{-i \kappa x_{r}} = \varphi_0(\xiLk) \hat{u}^{n}_k e^{-i \kappa p \Dx} = \varphi_0(\xiLk) \hat{u}^{n}_k e^{-i \kappa s} ,
    \end{aligned}
\end{equation}

\noindent with $s := x_j - x_{d_j} = x_j - x_{j-p} = p\Dx$. Therefore, for each $k$, the stability function $|\hat{u}^{n+1}_k| / |\hat{u}^{n}_k| = \varphi_0(\xiLk) e^{-i \kappa s}$ is the same as derived above (see Eq. \eqref{eq:A_SL_ETD1RK}), \xadded{which corresponds, in fact, the exact solution for wavenumber $k$.} However, if the velocity field is not constant, we are not able to factor out $e^{-i \kappa p_j \Dx}$ from the summation, and the stability function may be different, \xadded{and, in particular, the exact solution for wavenumber $k$ is no longer obtained; this suggests an intricate relation between spectral transforms, matrix exponentiation, and spatial interpolation.}

\indent Consider now the proposed modification to the discretization of the linear term. In the spectral space, we solve

\begin{equation}
    \label{eq:spectral_SLETD1RK_new_linear}
    \hat{u}_k^{n+1} = \varphi_0\left(\frac{\xiLk}{2}\right) \mc{F} ([ \mc{F}^{-1}(\hat{w})^n]_*), \qquad k = -M, \dots, M
\end{equation}

\noindent where 

\begin{equation}
    \label{eq:fourier_w}
    \hat{w}^n_k := \varphi_0\left(\frac{\xiLk}{2}\right)\hat{u}^n_k.
\end{equation}

\indent Conducting the same development as above,

\begin{equation}
    \label{eq:A_SLETD1RK_new_linear_vcte}
    \begin{aligned}
        \hat{u}_k^{n+1} = \varphi_0\left(\frac{\xiLk}{2}\right) \frac{1}{P}\sum_{j = 0}^{P-1} w(t_n, x_{r_j}) e^{-i \kappa (p_j + r_j) \Dx} = \varphi_0\left(\frac{\xiLk}{2}\right) \hat{w}^{n}_k e^{-i \kappa s} = \varphi_0(\xiLk) \hat{u}^{n}_k e^{-i \kappa s},
    \end{aligned}
\end{equation}

\noindent where the second equality holds if the velocity field is constant. We conclude that, if $v$ is constant, the original and the proposed discretizations of the linear term yield the same stability function. However, if the velocity field is not constant, we have

\begingroup
\allowdisplaybreaks
% \begin{equation}
    \begin{align}
        \hat{u}_k^{n+1} & = \varphi_0\left(\frac{{\xiLk}}{2}\right) \frac{1}{P}\sum_{j = 0}^{P-1} w(t_n, x_{r_j}) e^{-i \kappa (p_j + r_j) \Dx} \notag\\
        & =  \varphi_0\left(\frac{{\xiLk}}{2}\right) \frac{1}{P}\sum_{j = 0}^{P-1} \left[ \sum_{l = -M}^M \hat{w}^n_l e^{i \kappa_l x_{r_j}}  \right] e^{-i \kappa (p_j + r_j) \Dx} \notag\\
        & =  \varphi_0\left(\frac{{\xiLk}}{2}\right) \frac{1}{P}\sum_{j = 0}^{P-1} \left[ \sum_{l = -M}^M \varphi_0\left(\frac{ {\xi_{\vecL_l}}}{2}\right) \hat{u}^n_l e^{i \kappa_l x_{r_j}}  \right] e^{-i \kappa (p_j + r_j) \Dx} \notag\\
        & \neq \varphi_0({\xiLk}) \frac{1}{P}\sum_{j = 0}^{P-1} u(t_n, x_{r_j}) e^{-i \kappa (p_j + r_j) \Dx} \label{eq:A_SLETD1RK_new_linear}
    \end{align}
% \end{equation}
\endgroup
\noindent where $\kappa_l := 2\pi l /L$. Note that the dependency of ${\xi_{\vecL_l}}$ on $l$ prevents us from factoring out the exponential term from the inner summation, such that the inverse Fourier transform into the brackets does not yield $u(t_n, x_{r_j})$. As a consequence, different stability functions are provided by the two variants of the semi-Lagrangian exponential scheme.

\indent Therefore, by explicitly writing the spectral transforms involved in the integration of the semi-Lagrangian exponential schemes, it becomes clear why different discretizations of the linear term lead to different stability functions. We remark, however, that these differences are not necessarily linked to the fact that the equations are solved using a spectral method, but rather to the lack of commutation between the matrix exponential and the interpolation. For illustration, consider the equation

\begin{equation}
    \label{eq:linear_illustration}
    \dert{\vecu} = \lambdaL \vecu
\end{equation}

\noindent defined on the (periodic) physical space, with $\vecu = \vecthreeT{u_1}{\dots}{u_P} \in \reals^{P}$ and $\lambdaL = \text{diag}(\lambda_1, \dots, \lambda_P) \in \reals^{P \times P}$ a diagonal matrix with at least one element $\lambda_i$ different from the others. Consider also that the interpolation procedure consists simply of a one grid point shift to the right. Therefore, the integration of one time step of \eqref{eq:linear_illustration} using the two variants of the semi-Lagrangian exponential scheme yields respectively

\begin{equation}
    \label{eq:SE21_physical_grid}
    \vecu^{n+1} = \varphi_0(\Dt \lambdaL) [\vecu^n]_* = \matthree{\varphi_0(\Dt \lambda_1)}{}{0}{}{\ddots}{}{0}{}{\varphi_0(\Dt \lambda_P)}
    \left(
    \begin{array}{c}
         u^n_{P}  \\
         u^n_1 \\
         \vdots\\
         u^n_{P-1}
    \end{array}
    \right) =
    \left(
    \begin{array}{c}
         \varphi_0(\Dt \lambda_1)u^n_{P}  \\
         \varphi_0(\Dt \lambda_2)u^n_1 \\
         \vdots\\
         \varphi_0(\Dt \lambda_P)u^n_{P-1}
    \end{array}
    \right)
\end{equation}

\noindent and

\begingroup
\allowdisplaybreaks
% \begin{equation}
    \begin{align}
    \vecu^{n+1} & = \varphi_0\left(\frac{\Dt \lambdaL}{2}\right) \left[\varphi_0\left(\frac{\Dt \lambdaL}{2}\right)\vecu^n\right]_* \notag\\
    & = \matthree{\varphi_0(\Dt \lambda_1/2)}{}{0}{}{\ddots}{}{0}{}{\varphi_0(\Dt \lambda_P/2)}
    \left[\left(
    \begin{array}{c}
         \varphi_0(\Dt \lambda_1/2)u^n_1  \\
         \varphi_0(\Dt \lambda_2/2)u^n_2 \\
         \vdots\\
         \varphi_0(\Dt \lambda_P/2)u^n_P
    \end{array}
    \right)\right]_* \notag\\
    & = \matthree{\varphi_0(\Dt \lambda_1/2)}{}{0}{}{\ddots}{}{0}{}{\varphi_0(\Dt \lambda_P/2)}
    \left[\left(
    \begin{array}{c}
         \varphi_0(\Dt \lambda_P/2)u^n_P \\
         \varphi_0(\Dt \lambda_1/2)u^n_1 \\
         \vdots\\
         \varphi_0(\Dt \lambda_{P-1}/2)u^n_{P-1}
    \end{array}
    \right)\right]_* \notag\\
    & =
    \left(
    \begin{array}{c}
         \varphi_0(\Dt \lambda_1/2)\varphi_0(\Dt \lambda_P/2)u^n_{P}  \\
         \varphi_0(\Dt \lambda_2/2)\varphi_0(\Dt \lambda_1/2)u^n_1 \\
         \vdots\\
         \varphi_0(\Dt \lambda_P/2)\varphi_0(\Dt \lambda_{P-1}/2))u^n_{P-1}
    \end{array}
    \right) \label{eq:SE22_physical_grid}
    \end{align}
% \end{equation}
\endgroup

\noindent with the two schemes being equivalent only if $\lambdaL$ is a multiple of the identity matrix.

% Also note that, in this example, different stability functions are obtained when the velocity field is constant, since the stability is defined using the amplitudes of the solution in the physical space. In the development made above (eq. \eqref{eq:A_SLETD1RK_linear} and \eqref{eq:A_SLETD1RK_new_linear_vcte}), in which the amplitudes of the spectral coefficients define the stability region, identical stability regions are obtained when $v$ is constant, since it leads to the Fourier transform of a shifted function, which corresponds to the multiplication of the spectral coefficients by a constant complex exponential factor $e^{-i \kappa s}$.

\indent Note that, except for the case in which the velocity field is constant, determining the stability regions of \SE{*}{*} analytically, using \eqref{eq:A_SLETD1RK_linear} and \eqref{eq:A_SLETD1RK_new_linear}, is not feasible. In particular, we are not able to obtain an expression in the form $\hat{u}^{n+1}_k = A(\xiL, \xiN) \hat{u}^{n}_k $ for each $k$ separately, and $\hat{u}^{n+1}_k$ depends \emph{a priori} on all wavenumber modes. It means that we are no longer in the framework of a linear stability analysis, even if the problem \eqref{eq:pde_SL} itself is linear, with this \emph{nonlinearity} being a consequence of the lack of commutation between the exponential and interpolation operators. We choose therefore to assess the stability properties of the schemes through numerical simulations, as presented in Section \ref{subsec:stability_simulations}.

\section{Numerical tests}
\label{sec:numerical_simulations}

\indent In this section, we perform a set of numerical simulations of the SWE on the rotating sphere to evaluate the proposed modifications of the semi-Lagrangian exponential methods and to compare them with their original versions and with other relevant numerical schemes in this study, namely SL-SI-SETTLS and the Eulerian exponential integration methods. We begin by verifying the numerical orders of convergence using benchmark test cases with increasing complexity; then, we conduct an empirical stability study; finally, we consider a relatively challenging test to evaluate the methods qualitatively and in terms of computational cost, including a study of the application of hyperviscosity approaches.

\subsection{Numerical verification of the convergence order considering the advection equation}

\indent The goal of this section is to provide \replaced[id=R2]{an initial}{a first} numerical validation of the truncation error analyses developed in Sections \ref{sec:source_first_order} and \ref{subsec:LTE_second_order}, considering the governing equations in the Lagrangian framework and containing only the linear term on the right-hand side. 

\begin{figure}[!htbp]
    \begin{subfigure}{.475\linewidth}
        \centering
        \includegraphics[scale=0.45]{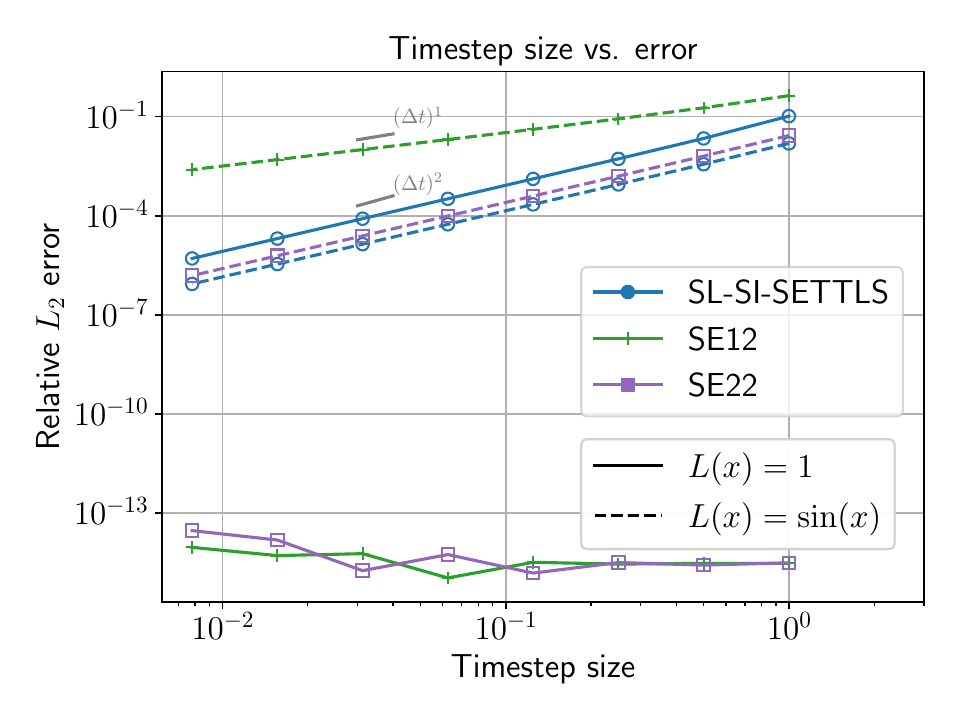}
        \caption{Scalar case\label{fig:errors_linear_scalar}}
    \end{subfigure}
    \begin{subfigure}{.475\linewidth}
        \centering
        \includegraphics[scale=0.45]{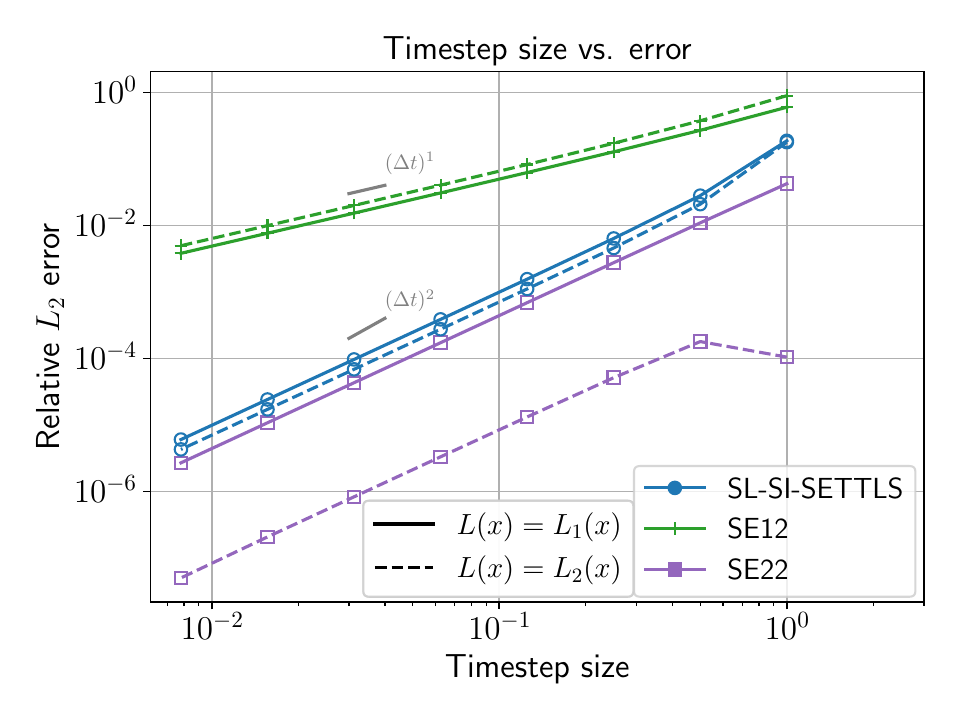}
        \caption{Vectorial case ($N=2$)\label{fig:errors_linear_vectorial}}
    \end{subfigure}
    \caption{\xadded{Numerical verification of the convergence order of the semi-Lagrangian schemes considered in this work (relative $L_{2}$ error \emph{vs} time step size), considering governing equations with only advective and linear terms.}}
    \label{fig:errors_linear}
\end{figure}

\indent We consider first the scalar, one-dimensional equation \eqref{eq:PDE}, in the periodic domain $\Omega = [0,10]$, integrated from $t = 0$ to $t = 10$, and a Gaussian curve centered at $x=5$ as initial solution. We consider time step sizes ranging between $\Dt = 2^{-5}$ and $\Dt = 1$, and the spatial domain is discretized considering $P + 1 = 2049$ equally spaced points. We consider the first- and second-order schemes \SE{1}{2} and \SE{2}{2} (without the nonlinear terms), and, for comparison, we also integrate using SL-SI-SETTLS. A reference solution is computed using a 4th-order explicit Runge Kutta method with time step size $\Dt = 2^{-5}/10$ and the same spatial resolution $P + 1 = 2049$, with the spatial derivatives being computed using the Fourier spectral method with spectral resolution $M = 1024$. Figure \ref{fig:errors_linear_scalar} presents the relative $L_{2}$ errors in two cases, respectively with constant $L(x) \equiv 1$, and variable $L(x) = \sin(x)$. In the former, both semi-Lagrangian exponential methods compute the exact solution (up to the error due to the spectral resolution in the reference simulation), since the integration factor \eqref{eq:integration_factor} is exactly $P_n(t_{n+1}) = \varphi_0(\Dt L)$. In the latter case, in which $L$ depends on space (thus on time, along each Lagrangian trajectory), the first- and second-order accuracies derived analytically are observed.

\indent In Figure \ref{fig:errors_linear_vectorial}, we present the same results in the vectorial case \eqref{eq:base_equation_slexp_linear}, with $N = 2$ and one spatial dimension. The same discretization parameters as in the scalar examples are considered. We perform simulations for two choices of $\vecL(x)$, respectively

\begin{equation}
    \label{eq:choices_L}
    \vecL(x) = \vecL_1(x) := \mattwo{\sin(x)}{\cos(x)}{\cos(x)}{\sin(x)}, \qquad \vecL(x) = \vecL_2(x) := \mattwo{\sin(x)}{\sin(x)}{\sin(x)}{\cos(x)}
\end{equation}

\noindent with the property that $\vecL_1(x)$ and $\vecL_1(y)$ commutes for every $x,y$, which does not hold for $\vecL_2$. In both cases, \SE{1}{2} provides first-order accuracy, due to the approximation of the integrals considering the space-dependent linear operators to be constant along the trajectories. On the other hand, \SE{2}{2} is second-order accurate in both cases, including the one in which the linear operator does not commute along the trajectory. The simulations presented in the following paragraphs, considering the SWE on the rotating sphere, for which the commutation assumption does not hold, confirm these results.

\subsection{Numerical verification of the convergence order \xadded{considering the SWE on the rotating sphere}}

\indent We consider here three test cases to verify the orders of convergence determined analytically. \added[id=R2-new]{To provide a more detailed insight into the stability properties of the time integration schemes, we present, for each simulation, the CFL numbers related to gravity wave propagation and advection. We consider the characteristic speed $c = \ol{\Phi}^{1/2} + |\ol{u}|$, where the square root of the mean geopotential, $\ol{\Phi}^{1/2}=\sqrt{g\ol{h}}$, defines the mean external gravity wave speed, and $\ol{u}$ is the maximum advective velocity magnitude defined by the initial conditions. A reference CFL number is then defined as $c\Dt/\Dx$, where $\Dx$ is the longitudinal grid size along the equator.}

\indent \added[id=R2-new]{The test cases considered here are:}

\begin{itemize}
    \item The Williamson's nonlinear test case 2 \cite{williamson_al:1992}, consisting of a geostrophically balanced steady solution,  \added[id=R2-new]{with mean geopotential $\ol{\Phi} = 2.94 \times 10^4 \text{m}^2/\text{s}^2$, mean advective velocity $|\ol{u}| \approx 39 \text{m}/\text{s}$ and characteristic speed $c \approx 210 \text{m}/\text{s}$.} The authors recommend error evaluations after 5 days of integration; here, we execute simulations for 7 days.
    \item The geostrophic balance test with bathymetry proposed by \cite{peixoto_al:2017}, which is a modification of Williamson's test case 2. In this test case, the fluid depth is constant ($h = h_0$) and the bathymetry profile is constructed to ensure geostrophic balance. By choosing a relatively small $h_0$, nonlinear effects are not negligible, allowing us to evaluate the influence of the discretization of the nonlinear terms on the convergence order. We consider $h_0 = 100\text{m}$ and $h_0 = 1\text{m}$\added[id=R2-new]{, with mean geopotential $\ol{\Phi} = gh_0$, mean advective velocity $|\ol{u}|$ and characteristic speed $c$ respectively $(\ol{\Phi},|\ol{u}|,c) \approx (981,39,70)$ and $(\ol{\Phi},|\ol{u}|,c) \approx (9.81,39,42)$ (in $\text{m}^2/\text{s}^2$, $\text{m}/\text{s}$ and $\text{m}/\text{s}$).} The simulations are also executed for 7 days.
    \item The unstable jet test case proposed by \cite{galewski_al:2004}, which is a popular and challenging benchmark in atmospheric modeling. In this test, with null bathymetry and whose initial solution is a stationary zonal jet, a small Gaussian bump perturbation on the geopotential field leads to the formation of characteristic vortices and important nonlinear processes that \added[id=R2-new]{transfers energy along the spectra, and allows enstrophy cascade} to the small-scale components of the solution. \added[id=R2-new]{The mean fluid depth and geopotential are respectively $\ol{h} = 10000\text{m}$ and $\ol{\Phi} \approx 98100$ $\text{m}^2/\text{s}^2$, and the advective velocity is $|\ol{u}| = 80\text{m}/\text{s}$, resulting in the characteristic speed $c \approx 393 \text{m}/\text{s}$.} A similar test case has been used for numerical validation of \SE{1}{1} and \SE{1}{2} on the plane by \cite{peixoto_schreiber:2019} and, as done in that work, errors are measured after one day of simulation, before the formation of the vortices. We also conduct qualitative analyses considering longer times of simulation of this test case.
\end{itemize}

\begin{figure}[!htbp]
    \begin{subfigure}{.495\linewidth}
        \centering
        \includegraphics[scale=.4]{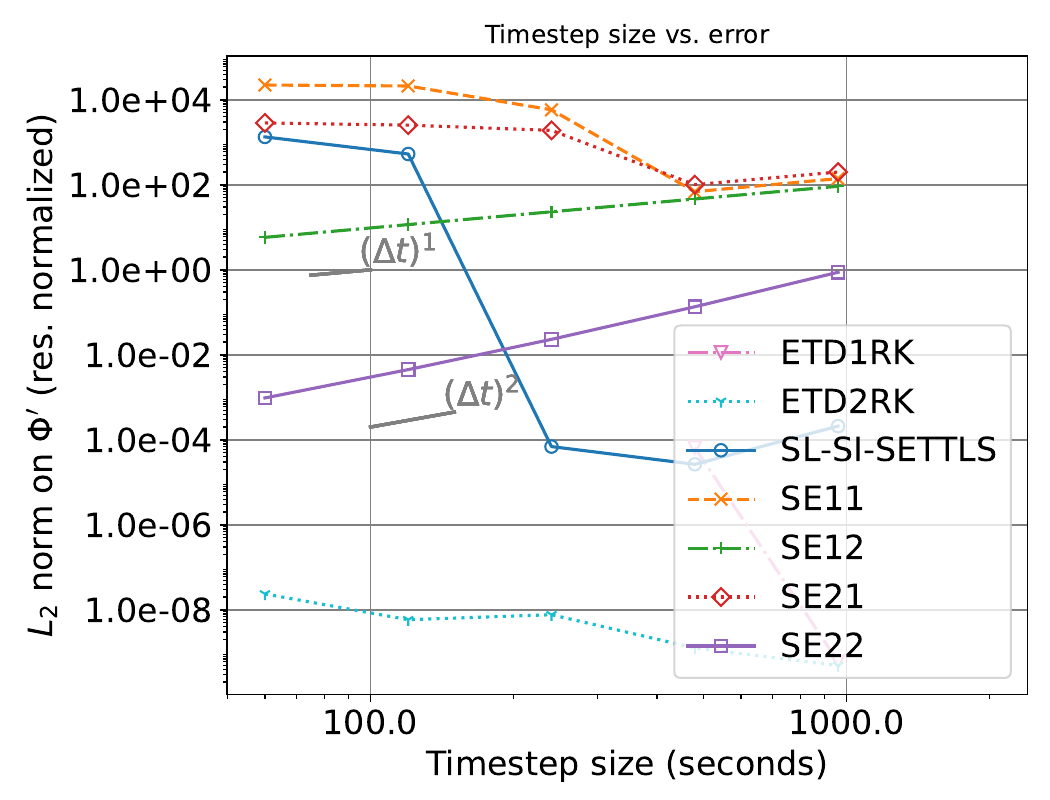}
        \caption{Geostrophic balance}
        \label{fig:error_serial_convergence_geostrophic_balance}
    \end{subfigure}
    \begin{subfigure}{.495\linewidth}
        \centering
        \includegraphics[scale=.4]{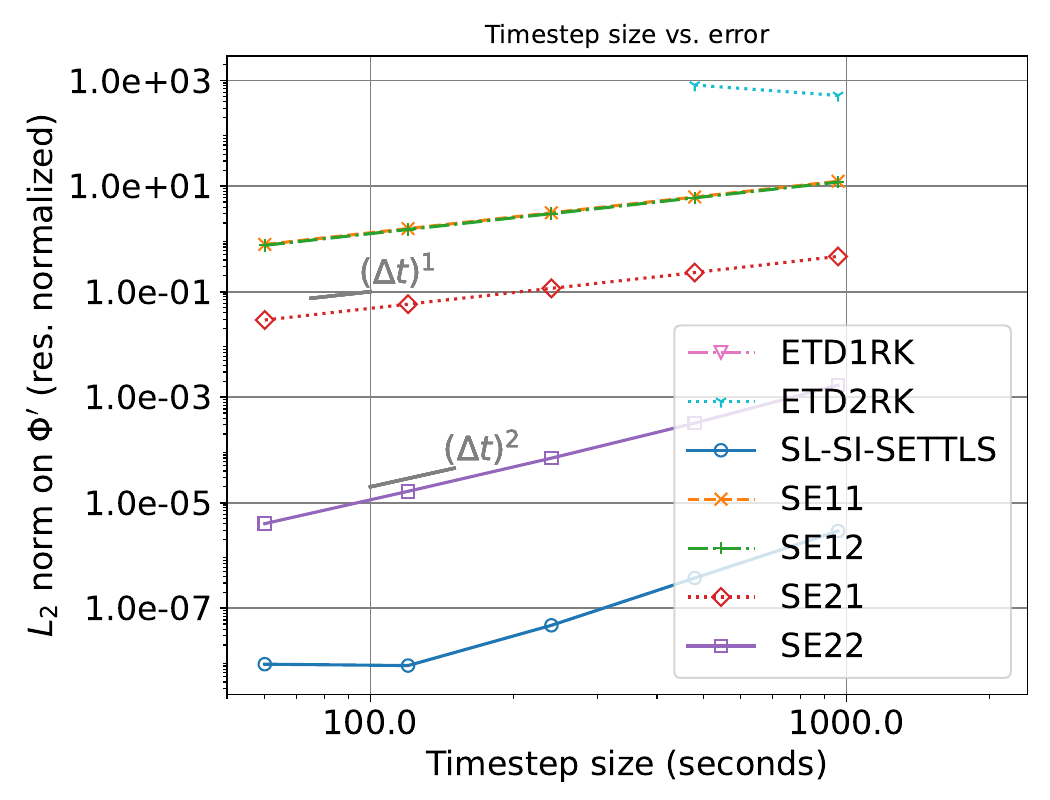}
        \caption{Geostr. balance w/ topography; $h_0 = 100\text{m}$}        \label{fig:error_serial_convergence_topography_100m}
    \end{subfigure}
    \begin{subfigure}{.495\linewidth}
        \centering
        \includegraphics[scale=.4]{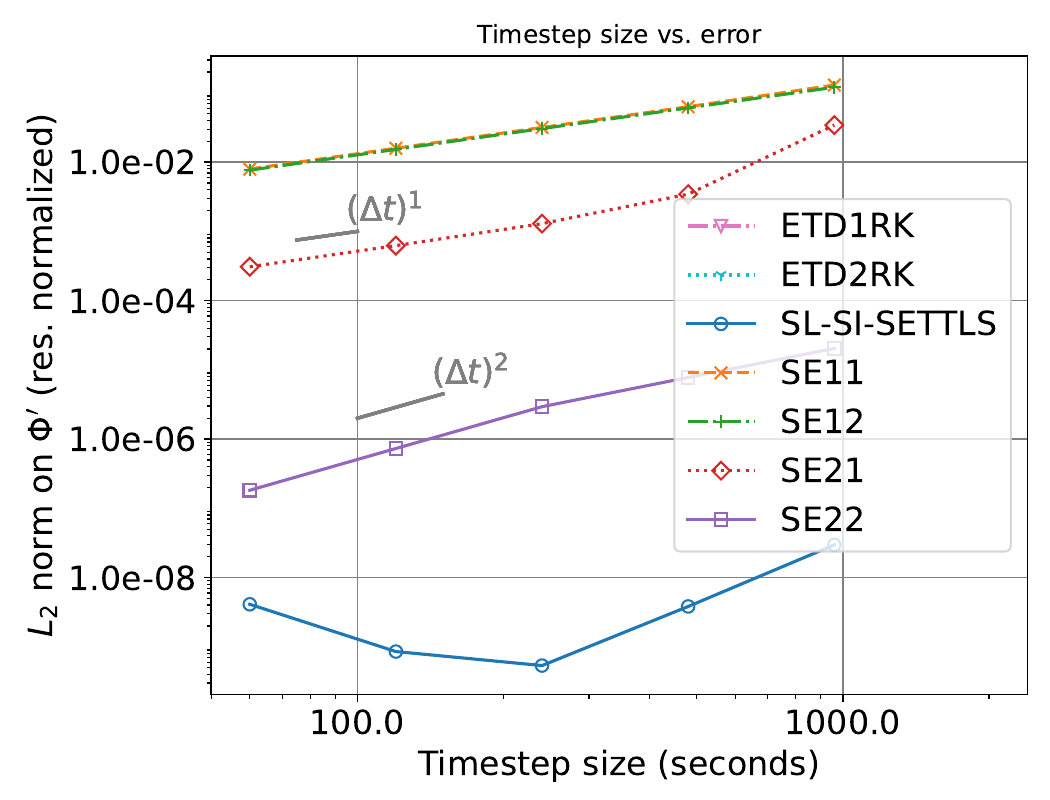}
        \caption{Geostr. balance w/ topography; $h_0 = 1\text{m}$}\label{fig:error_serial_convergence_topography_1m}
    \end{subfigure}
    \begin{subfigure}{.495\linewidth}
        \centering
        \includegraphics[scale=.4]{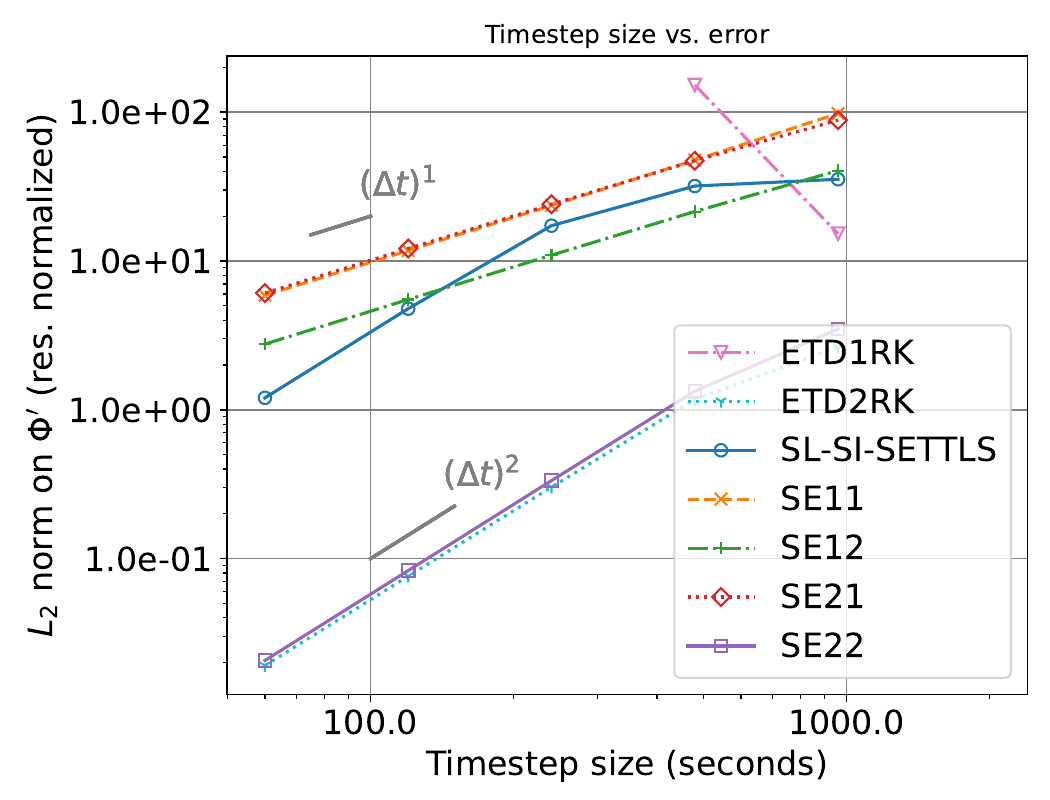}
        \caption{Unstable jet}\label{fig:error_serial_convergence_unstable_jet}
    \end{subfigure}
    \caption{Numerical verification of the convergence order of the schemes considered in this work (normalized $\xLtwo$ norm of the error on the geopotential field perturbation \emph{vs} the time step size). Similar results are obtained with the $\xLinfty$ norm.}
    \label{fig:error_serial_convergence}
\end{figure}

\indent In all test cases and for all schemes, we perform simulations by keeping the ratio between the temporal and spectral resolutions constant, with $\Dt \in \{60,120,240,480,960\}$ (s) and $M \in \{512,256,128,64,32\}$, \added[id=R2-new]{such that the CFL numbers are kept constant for each test case (namely, $0.48$, $0.16$, $0.10$ and $0.91$, respectively for the Williamson's test case 2, the geostrophically balanced test case with bathymetry $(h_0 = 100\text{m}$ and $h_0 = 1\text{m})$, and the unstable jet test case).} The errors in each simulation are computed \wrt a reference solution integrated using a fourth-order explicit Runge-Kutta scheme, with a time step four times smaller and the same spectral resolution (thus allowing us to analyze exclusively the errors due to the temporal discretization without imposing too severe stability constraints by considering a fixed high spectral resolution). None of the simulations use artificial viscosity, \ie $\nu = 0$. Figure \ref{fig:error_serial_convergence} reveals that, among the semi-Lagrangian exponential methods, the second-order accuracy is achieved only by the proposed \SE{2}{2} scheme, with all the other semi-Lagrangian exponential schemes being of first order. This is observed in all test cases, including those with relatively important nonlinear effects (geostrophic balance with topography and small $h_0$). It shows that both our proposed modification of the linear term and the second-order discretization of the nonlinear one proposed by \cite{peixoto_schreiber:2019} are required to ensure global second-order accuracy.

\indent Figure \ref{fig:error_serial_convergence} also confirms our hypotheses on the stability and accuracy of semi-Lagrangian exponential methods. These schemes, notably those with a second-order discretization of the nonlinear term, have improved stability properties compared to the Eulerian ones; indeed, ETD2RK and mainly ETD1RK are unstable in most of the simulations. The only simulations in which ETD2RK is stable are those with less important nonlinear effects, meaning that a weaker \replaced[id=R2-new]{enstrophy}{energy} cascade takes place, namely the geostrophically balanced steady-state solution (Figure \ref{fig:error_serial_convergence_geostrophic_balance}) and the short-run unstable jet simulation (Figure \ref{fig:error_serial_convergence_unstable_jet}). In both simulations, ETD2RK and \SE{2}{2} outperform the accuracy of the also second-order accurate SL-SI-SETTLS (which develops instabilities in the geostrophic balance test case), reflecting the accurate integration of the linear term by the exponential schemes. \xadded{Possibly, the observed instabilities in the geostrophic balance test case (which arise after 3-4 days of simulation) are due to the amplification of small departures from the steady-state solution, which may be caused by the time-splitting approach for integrating equations including linear and nonlinear terms; a detailed study on these amplifications is conducted in Section \ref{subsec:stability_simulations}.} On the other hand, SL-SI-SETTLS is more accurate in the geostrophic balance test case with topography (Figures \ref{fig:error_serial_convergence_topography_100m} and \ref{fig:error_serial_convergence_topography_1m}), in which the nonlinear effects are more relevant due to the small mean depth $h_0$. In this case, although still convergent with the expected accuracy, the semi-Lagrangian exponential schemes produce larger errors. We remark, however, that the observed second-order accuracy of \SE{2}{2} in this test case confirms the discussion in Section \ref{sec:source_first_order} in the sense that the source of the first-order accuracy of \SE{1}{2} is indeed the discretization of the linear term, and no modification of the discretization of the nonlinear term is required for improving its order of accuracy. Finally, note that some simulations reach an error plateau around $10^{-8}-10^{-9}$, which is a consequence of our formulation considering the geopotential perturbation, $\Phi = \ol{\Phi} + \Phi'$, with $\ol{\Phi}$ being several orders of magnitude larger than $\Phi'$; therefore, we consider these solutions to have converged.

\begin{figure}
    \centering
    \includegraphics[scale=.4]{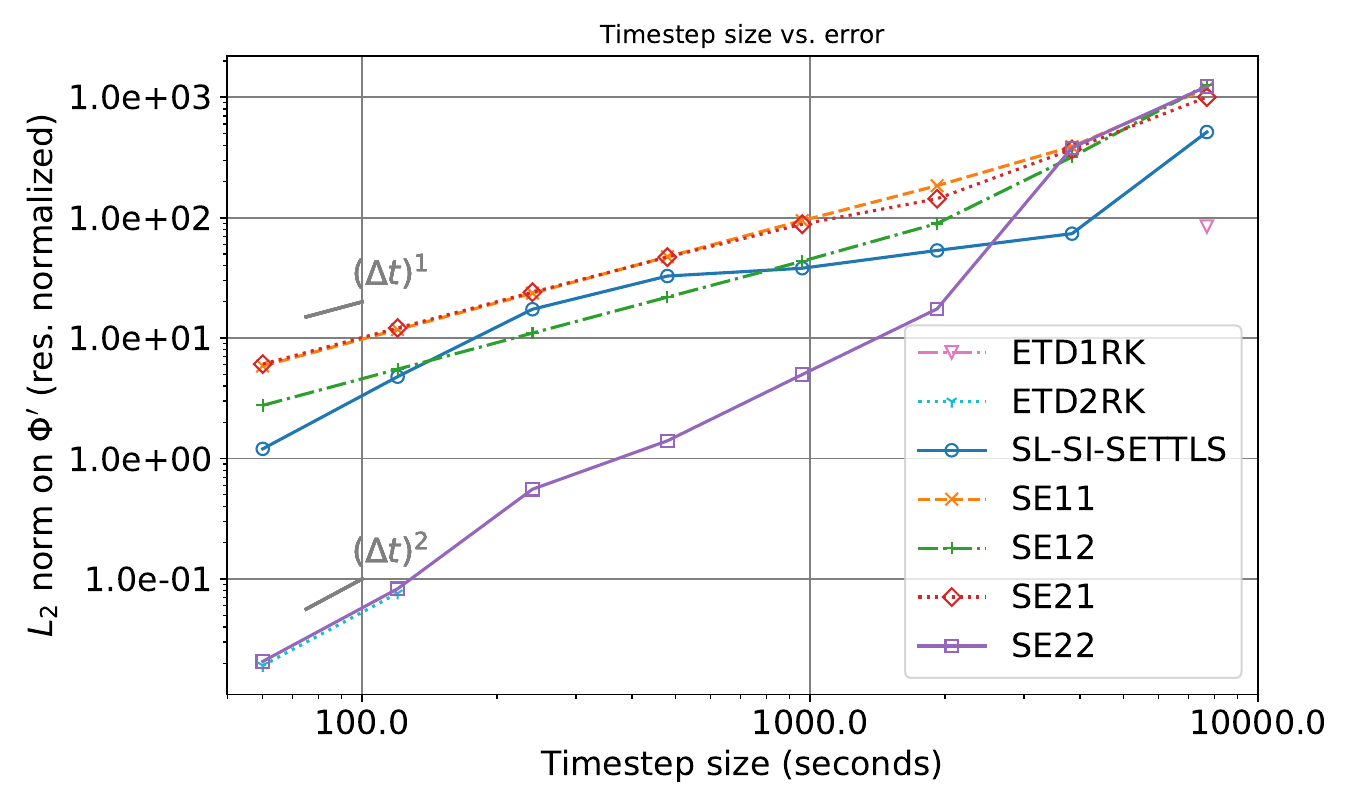}
    \caption{Numerical verification of the convergence order of the schemes considered in this work (normalized $\xLtwo$ norm of the error on the geopotential field perturbation \emph{vs} the time step size). All simulations consider the same spectral resolution $(M=256)$; the reference solution is computed using a fourth-order explicit Runge-Kutta scheme, with $M = 512$ and $\Dt = 2 \text{s}$.}
    \label{fig:error_serial_convergence_fixedM}
\end{figure}

\indent Since our focus in these first results was on the convergence properties of the schemes, we ensured stability in almost all simulations by considering a fixed CFL number. However, we now highlight the advantages of the proposed \SE{2}{2} \wrt the other schemes. We repeat the simulations of the unstable jet test case, but now with a fixed spectral resolution ($M=256$) and time step sizes chosen in $[60, 8192]$ (s)\added[id=R2-new]{, with approximate smallest and largest CFL numbers $c\Dt/\Dx$ of $0.45$ and $62$, respectively.} The reference solution is computed using a fourth-order explicit Runge-Kutta scheme, with $M = 512$ and $\Dt = 2 \text{s}$\added[id=R2-new]{, resulting in an approximate CFL number of 0.03.} The results are presented in Figure \ref{fig:error_serial_convergence_fixedM}. It is clear that the proposed second-order semi-Lagrangian exponential scheme \SE{2}{2} \quotes{extends} the error curve of the second-order Eulerian exponential method ETD2RK, \ie it has a largely improved stability range, maintaining approximately the same accuracy; indeed, the second-order accuracy of \SE{2}{2} is observed with time step sizes as large as approximately $2000\text{s}$\added[id=R2-new]{, corresponding to CFL numbers as large as 15.} (the small peak in the \SE{2}{2} curve will be discussed in Sections \ref{subsec:stability_simulations} and \ref{subsec:qualitative_comparison}), whereas the Eulerian scheme is unstable for $\Dt > 120\text{s}$. \added[id=R2-new]{It indicates a better ability of the proposed method to stably propagate fast gravity waves. }We also observe that, in this test case, \SE{2}{2} is much more accurate than SL-SI-SETTLS and the first-order semi-Lagrangian exponential schemes.

\subsection{Simulation-based stability analysis}
\label{subsec:stability_simulations}

\indent As discussed in Section \ref{subsec:serial_stability}, the lack of commutation between the exponential and interpolation operators does not allow us to compare the stability properties of the semi-Lagrangian exponential schemes through a classical linear stability analysis. Therefore, we conduct an empirical simulation-based stability study, as proposed by \cite{peixoto_al:2017}, in which simulations of the governing equations (here, the SWE on the rotating sphere), considering a steady-state solution, are performed in order to estimate the largest eigenvalue of the Jacobian of the approximate time evolution operator, which allows us to estimate the largest growth rate of perturbations to the steady state. This procedure, based on the power method for the estimation of the largest eigenvalue of a matrix, is briefly recalled below. We refer the reader to Appendix B of \cite{peixoto_al:2017} for details and parameters to be used. Let

\begin{equation}
    \label{eq:power_method}
    \vecU^{n+1} = \vecG(\vecU^{n})
\end{equation}

\noindent represent the evolution of a given state vector $\vecU$ between times $t_n$ and $t_{n+1}$ using a given time discretization method. Let $\ol{\vecU}$ be a steady state, such that $\vecR^{n} := \vecU^n - \ol{\vecU}$ is the perturbation at time $n$. Through linearization of \eqref{eq:power_method}, the perturbation satisfies, to first order

\begin{equation}
    \label{eq:power_method2}
    \vecR^{n+1} = \vecG'(\ol{\vecU})\vecR^n
\end{equation}

\indent As in the power method, the procedure is to iterate over \eqref{eq:power_method2}. If this method converges, then $\vecR^n$ converges to the eigenvector associated with the largest (in absolute value) eigenvalue $\ol{\lambda}$ of $\vecG'$. This converged dominant eigenvalue can also be computed, as well as the growth rate $\ol{\nu}:= \log{\ol{\lambda}}/\Dt$ and the e-folding time $\ol{\tau}$, defined as the time in which $\vecR$ increases by a factor equal to $e$. We have

\begin{equation}
    \label{eq:e-folding_time}
    e = \frac{\normG{\vecR^{n}}}{\normG{\vecR^0}} = |\ol{\lambda}|^n = e^{\ol{\nu} n \Dt} \implies \ol{\tau} = n \Dt = \frac{1}{\ol{\nu}} = \frac{\Dt}{\log \ol{\lambda}}
\end{equation}

\indent The e-folding time is a relevant parameter, since, even if the scheme is unstable (\ie $|\ol{\lambda}| > 1$), a large e-folding time may indicate that instabilities will remain reasonably small within simulation times used in practical applications, \eg in numerical weather prediction, such that the time integration scheme can be successfully used in such contexts.

\begin{figure}[!htbp]
    \begin{subfigure}{.45\linewidth}
        \centering
        \includegraphics[scale=.4]{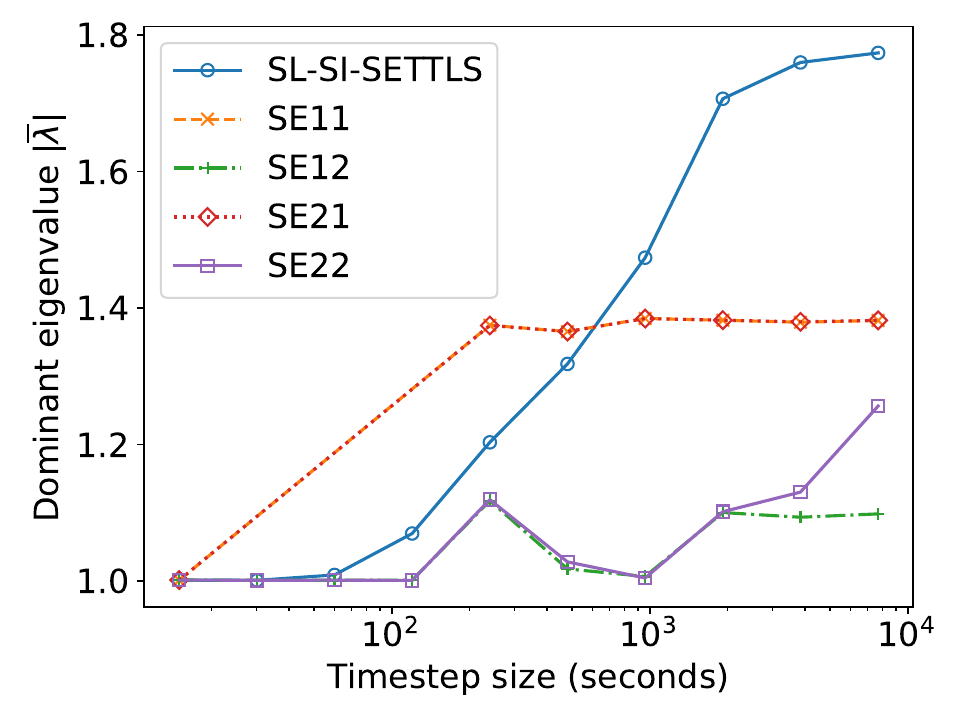}
    \end{subfigure}
    \begin{subfigure}{.45\linewidth}
        \centering
        \includegraphics[scale=.4]{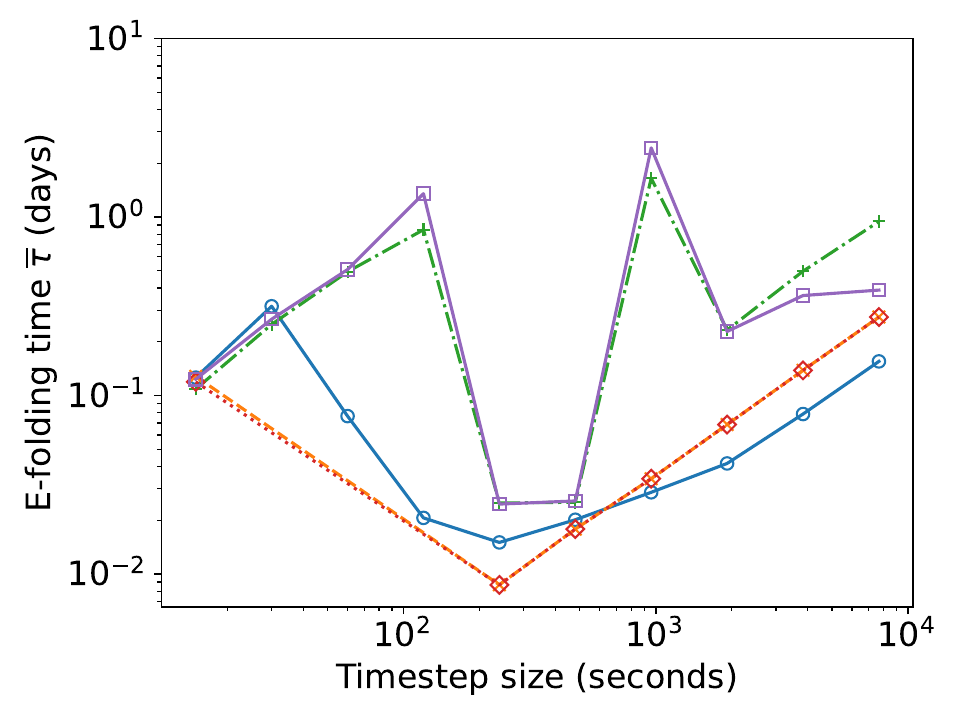}
    \end{subfigure}
    \caption{Simulation-based stability analysis considering the geostrophic balanced test case: dominant eigenvalues (left) and corresponding e-folding times (right) as a function of the time step size and the time integration scheme. The same legend holds for both plots. \added[id=R2-new]{The smallest reported eigenvalues are approximately $1.001$.}}
    \label{fig:simulation_stability}
\end{figure}

\indent We conduct this stability study by considering the geostrophic balance test case \added[id=R2-new]{(Williamson's test case 2)}, to which we introduce small initial perturbations to trigger instabilities. The spectral resolution is set to $M = 512$, and we compute $\ol{\lambda}$ and $\ol{\tau}$ for various time step sizes and for the various time integration schemes considered here, \added[id=R2-new]{with CFL numbers $c\Dt/\Dx$ ranging from $0.12$ to $62$, approximately.} Results are presented in Figure \ref{fig:simulation_stability}, which only includes simulations for which the iterative procedure described above converges. \added[id=R2-new]{The smallest reported eigenvalues are approximately 1.001.} We do not observe, in this test case, clear differences in the stability behavior depending on the discretization of the linear term in the semi-Lagrangian exponential schemes, with very similar results between \SE{1}{1} and \SE{2}{1}, and also between \SE{1}{2} and \SE{2}{2}, \ie the stability seems to be mostly determined by the discretization of the nonlinear term. The \SE{1}{2} and \SE{2}{2} schemes clearly outperform the other methods, including SL-SI-SETTLS, for which the dominant eigenvalue rapidly increases when the time step sizes get larger. Also, unlikely the other schemes, \SE{1}{2} and \SE{2}{2} present e-folding times of more than one day for some time step sizes as large as approximately $\Dt = 1000\text{s}$\added[id=R2-new]{ (approximate CFL number of 8.1)}, indicating that simulations can be conducted for relatively large durations with instabilities remaining quite controlled. We recall that no artificial viscosity is considered in this study. Finally, we highlight the unnatural peak on the dominant eigenvalue (corresponding to the abrupt reduction of the e-folding time) of \SE{1}{2} and \SE{2}{2} for the tested time step sizes $\Dt = 240 \text{s}$ and $\Dt = 480\text{s}$. As described in Section \ref{subsec:qualitative_comparison}, this behavior is also observed in the simulation of the unstable jet test case, but the reasons are unknown; it may be caused, \eg by possible bad behaviors of the exponential functions when computed at time step sizes close to the mentioned values. As identified by \cite{crouseilles_al:2020} through linear stability analyses, exponential integration methods usually have quite erratic stability behaviors, including asymmetric stability regions and discontinuities of the amplification factor as a function of the time step size. 

\subsection{Qualitative comparison}
\label{subsec:qualitative_comparison}

\indent We consider the unstable jet test case to conduct a detailed qualitative comparison of the solutions provided by \SE{1}{1}, \SE{1}{2}, their proposed modifications and SL-SI-SETTLS. This test case is standard in atmospheric modeling research due to its complex dynamics after some days of simulation, and guidelines for evaluating numerical methods applied to it have been proposed by \cite{scott_al:2015}. For this study, we consider a fixed spectral resolution $M = 512$ and various time step sizes in the same range $[60,960] $ (s) considered in the convergence analysis, \added[id=R2-new]{with respective CFL numbers between 0.91 and 14, approximately}. A reference solution is obtained from a simulation using a fourth-order explicit Runge-Kutta method, with $M=1024$ and $\Dt = 2\text{s}$, \added[id=R2-new]{corresponding to an approximate CFL number of 0.06.} Results using the Eulerian exponential schemes are not considered due to their highly unstable behavior for all tested time step sizes. 

\begin{figure}[!htbp]
    \begin{subfigure}{.475\linewidth}
        \centering
        \includegraphics[scale=.35]{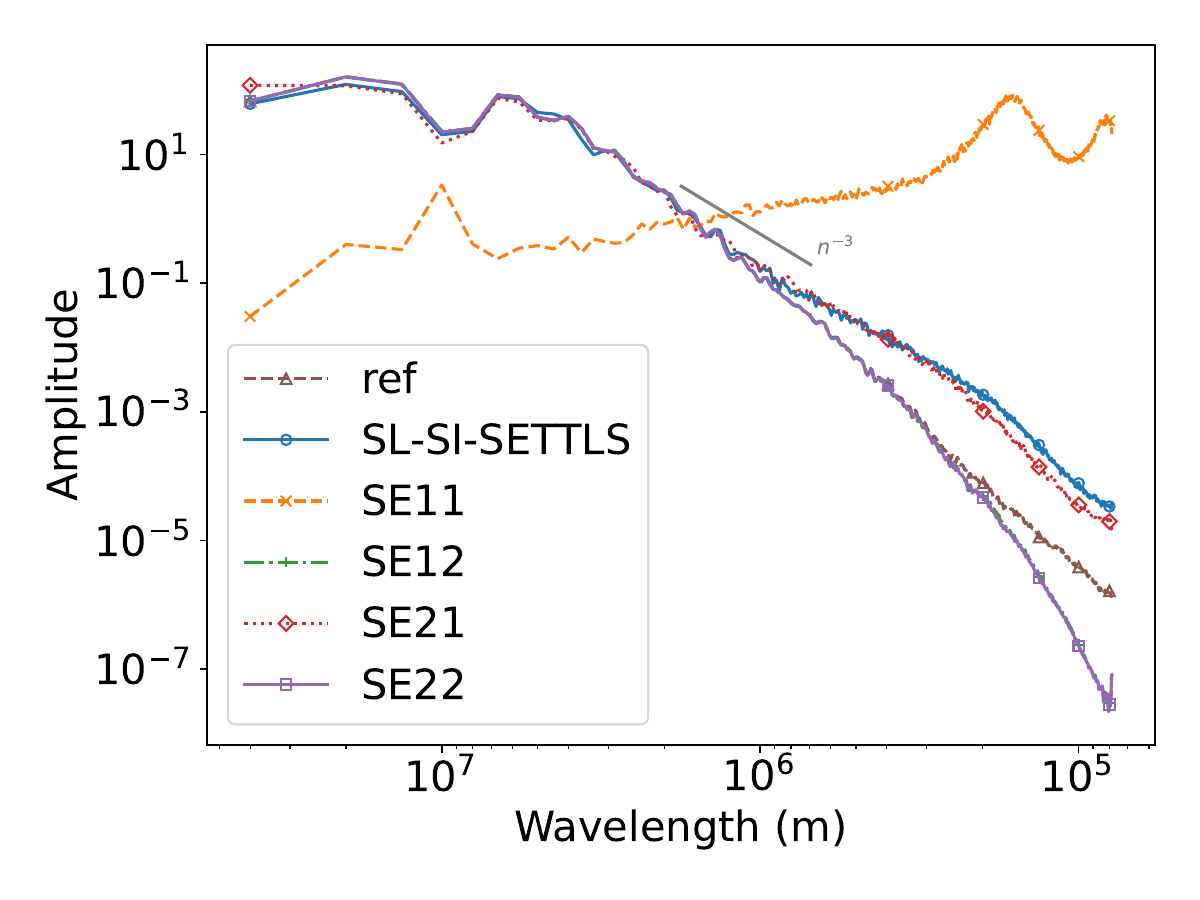}
        \caption{$\Dt = 120 \text{s}$, 6 days}\label{fig:spectrum_serial_120}
    \end{subfigure}
    \begin{subfigure}{.475\linewidth}
        \centering
        \includegraphics[scale=.35]{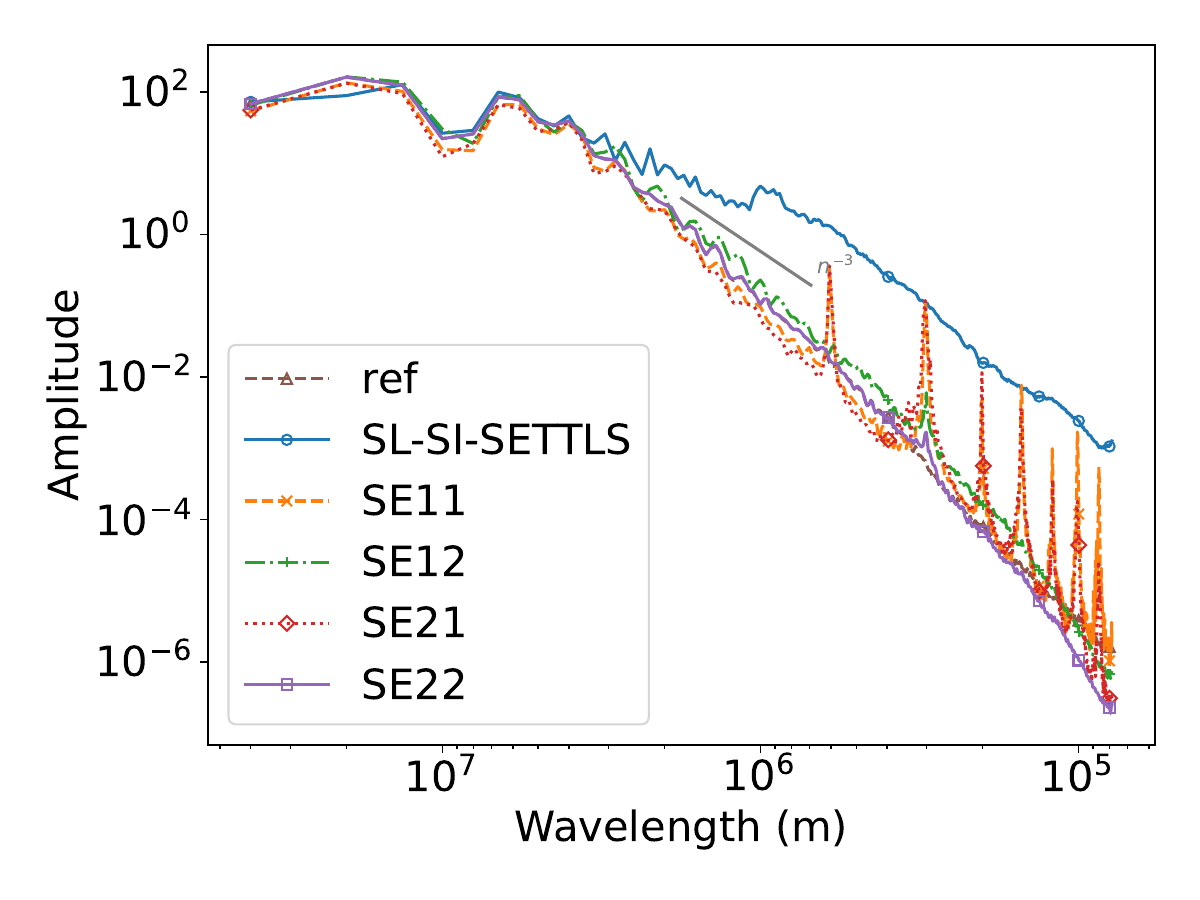}
        \caption{$\Dt = 960 \text{s}$, 6 days}\label{fig:spectrum_serial_960}
    \end{subfigure}
    \begin{subfigure}{\linewidth}
        \centering
        \includegraphics[scale=.35]{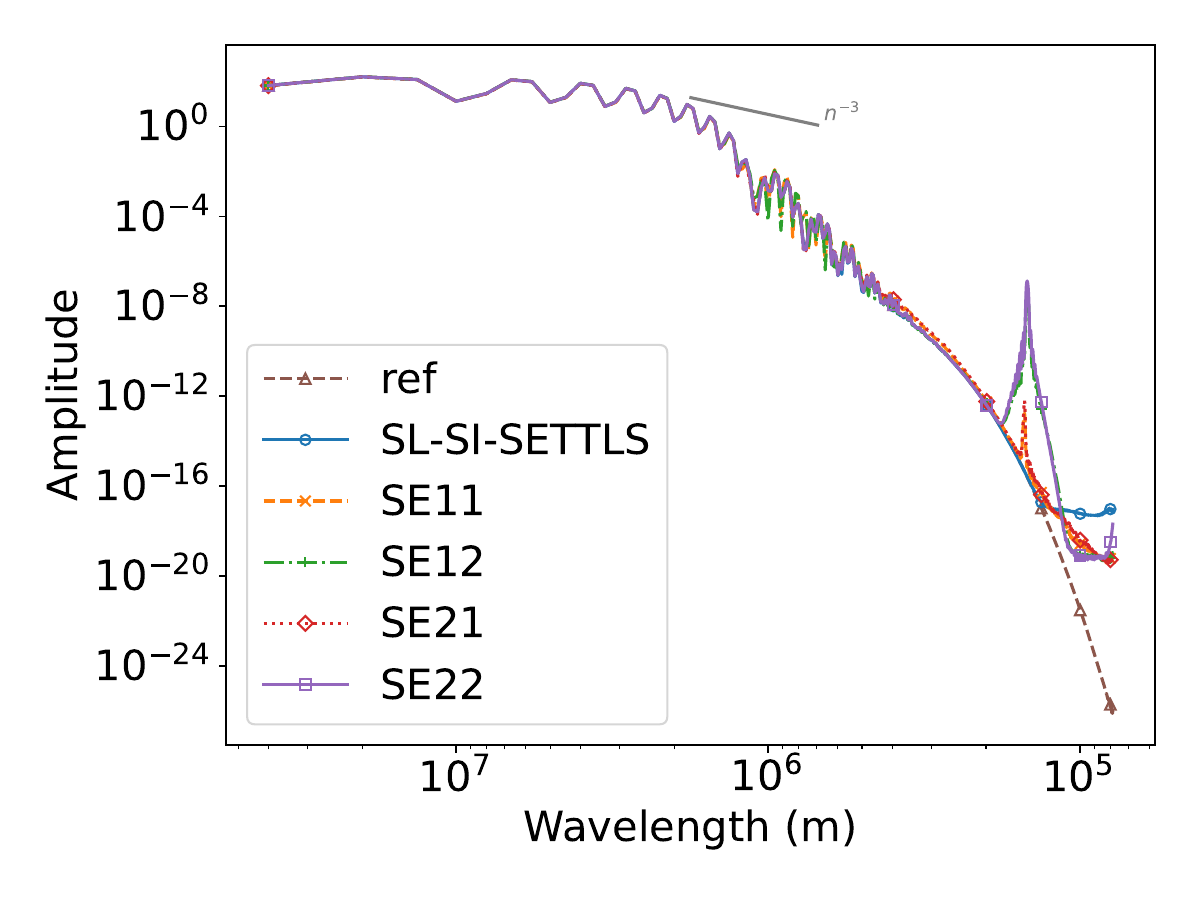}
        \caption{$\Dt = 240 \text{s}$, 1 day}\label{fig:spectrum_serial_240}
    \end{subfigure}
    \caption{Kinetic energy spectra, truncated at resolution $M = 512$, for the unstable jet test case.}
    \label{fig:spectrum_serial}
\end{figure}

\indent We first present, in Figures \ref{fig:spectrum_serial_120} and \ref{fig:spectrum_serial_960}, the kinetic energy spectra produced after six days of simulation by the four semi-Lagrangian exponential schemes, as well as SL-SI-SETTLS, with $\Dt = 120 \text{s}$ and $\Dt = 960 \text{s}$, and compared to the reference solution. All spectra, including the reference one, present upturned tails, indicating the accumulation of energy on the highest wavenumbers due to the \replaced[id=R2-new]{enstrophy}{energy} cascade, which may lead to instabilities in longer executions if there is not enough dissipation in the numerical scheme \cite{skamarock:2004}; thus, an appropriate application of artificial viscosity or hyperviscosity could be considered in order to mitigate unstable behaviors.

\indent In the case $\Dt = 120\text{s}$, the \SE{1}{1} scheme presents a clear unstable behavior, with almost the entire spectrum being strongly overamplified. This behavior is partially controlled with the proposed modified discretization of the linear term, and \SE{2}{1} presents a smaller overamplification of intermediate wavenumbers; concerning \SE{1}{2} and \SE{2}{2}, their spectra, which almost coincide visually, accurately reproduce the reference one in small and intermediate wavenumbers, but are damped out in the largest ones, indicating a stronger diffusion due to second-order discretization of the nonlinear term; finally, we observe an unstable behavior of SL-SI-SETTLS at the end of the wavenumber spectrum, being outperformed by the semi-Lagrangian exponential schemes using a second-order discretization of the nonlinear term.

\indent When the larger time step size $\Dt = 960\text{s}$ is considered, all time integration schemes present overamplifications of at least a portion of the wavenumber spectrum, notably SL-SI-SETTLS. The most stable simulation is provided by our proposed second-order scheme \SE{2}{2}, whose spectrum presents only a small peak at medium wavenumbers, being able to well reproduce the reference one at the large scales and remaining below it at the fine ones. Note that \SE{1}{2}, which provides nearly identical results to \SE{2}{2} when $\Dt = 120\text{s}$, has now a much more pronounced unstable behavior, with overamplification of almost the entire spectrum.

% \begin{figure}[!htbp]
%     \centering
%     \includegraphics[scale=.75]{alternative_SLETD2RK/figures/spectrum_kinetic_energy_t00000518400.00000000.pdf}
%     \caption{Kinetic energy spectrum after six days of simulation.}
%     \label{fig:spectrum_serial}
% \end{figure}

\indent We remark that, among the tested time step sizes, the four semi-Lagrangian exponential schemes present very early unexpected unstable behaviors when $\Dt$ is chosen around $240$s and $480$s, which has been identified in the stability analysis presented above (Figure \ref{fig:simulation_stability}). Figure \ref{fig:spectrum_serial_240} illustrates the kinetic energy spectrum after a single day of simulation (\ie long before the vortices start to appear), in which we observe unnatural peaks in the spectrum of the solutions of the semi-Lagrangian exponential methods. This behavior is not observed, at least not this early, for reasonably larger time step, \eg $\Dt = 960 \text{s}$ as seen above. The reasons for this are not yet understood. A possible explanation is that one or more of the exponential functions $\varphi_k$, $\psi_k$ may be not too well behaved when computed at $\Dt \vecL$, with $\vecL$ defined by the parameters of the SWE on the rotating sphere and $\Dt$ around the values mentioned above. Therefore, the use of artificial viscosity or hyperviscosity would be required to ensure further stability. This is studied in Section \ref{subsec:viscosity}.

\begin{figure}[!htbp]
    \begin{subfigure}{.3\linewidth}
        \centering
        \includegraphics[scale=.225]{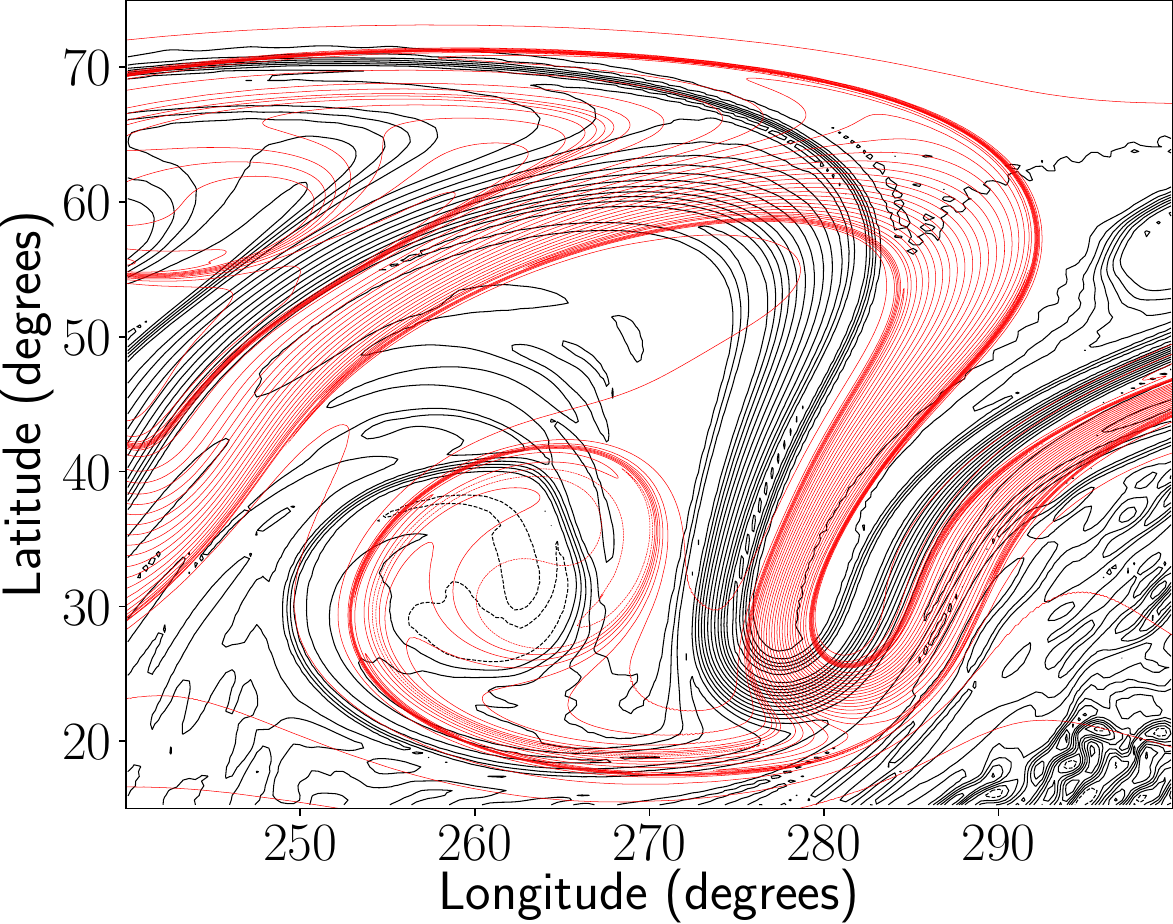}
        \caption{SL-SI-SETTLS, $\Dt = 120\text{s}$}
    \end{subfigure}
    \begin{subfigure}{.3\linewidth}
        \centering
        \includegraphics[scale=.225]{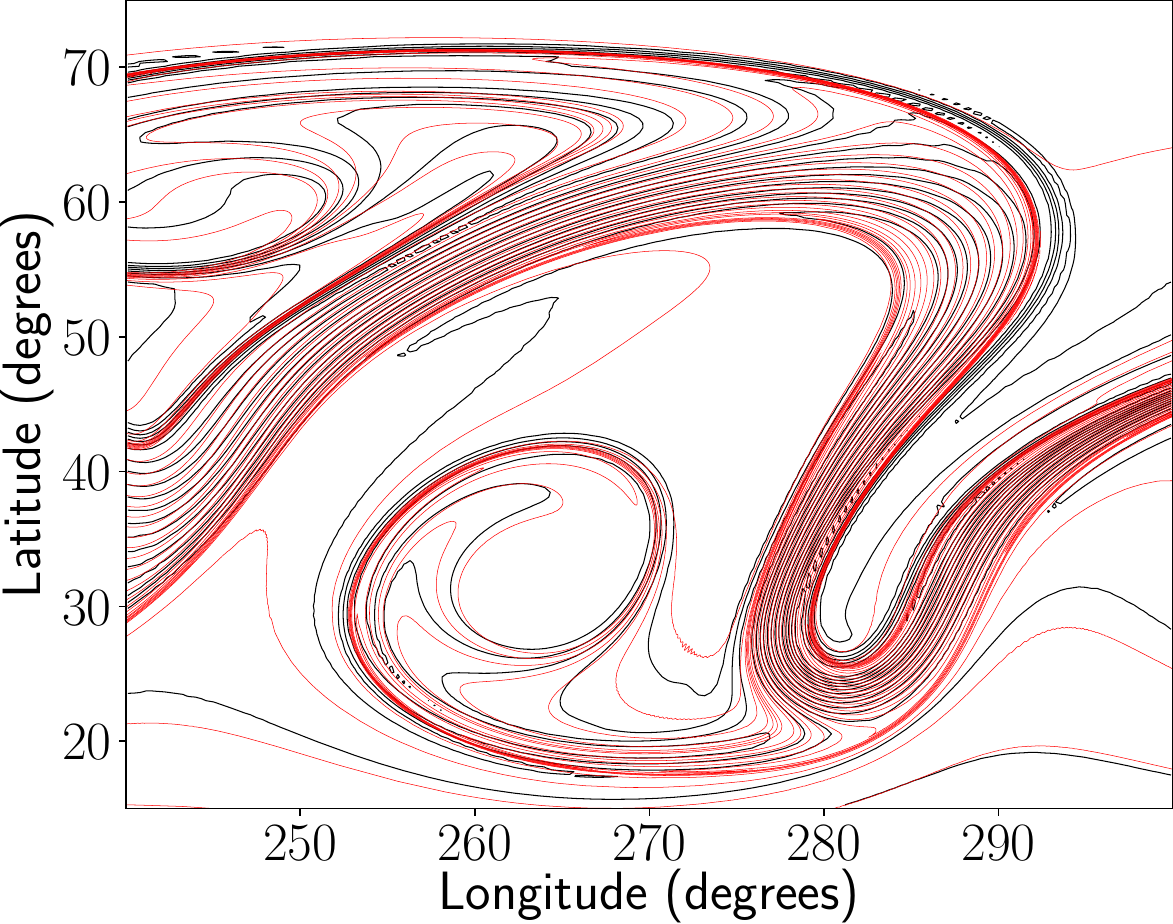}
        \caption{\SE{1}{2}, $\Dt = 120\text{s}$}
    \end{subfigure}
    \begin{subfigure}{.3\linewidth}
        \centering
        \includegraphics[scale=.225]{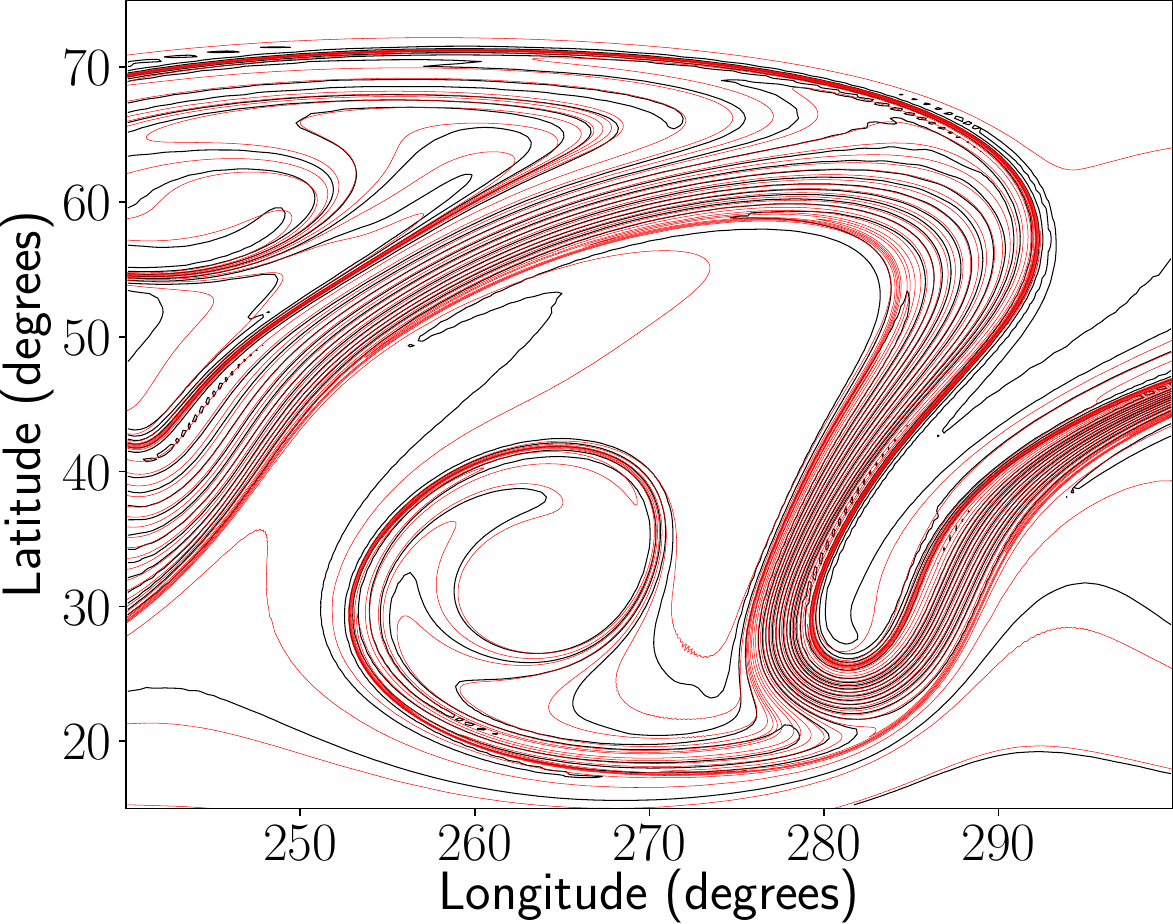}
        \caption{\SE{2}{2}, $\Dt = 120\text{s}$}
    \end{subfigure}
    \begin{subfigure}{.3\linewidth}
        \centering
        \includegraphics[scale=.225]{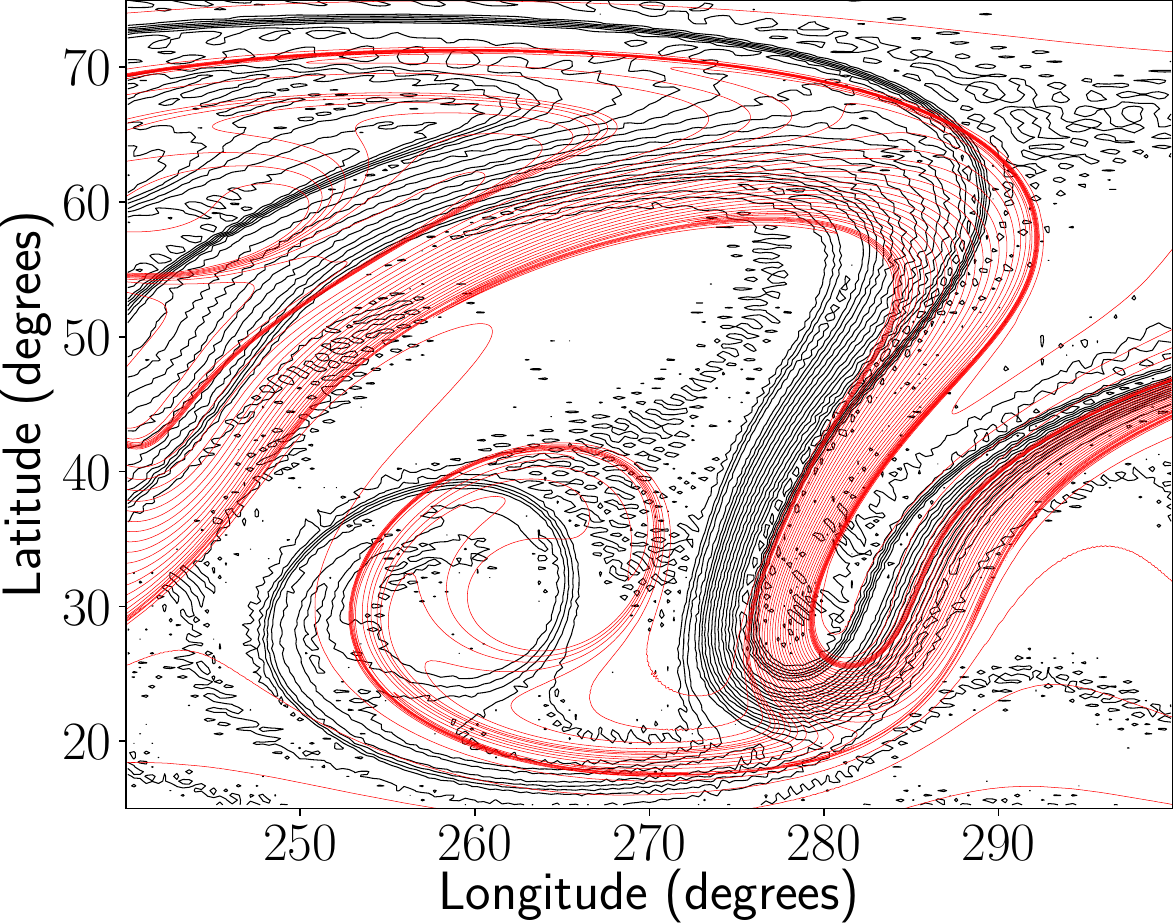}
        \caption{SL-SI-SETTLS, $\Dt = 960\text{s}$}
    \end{subfigure}
    \begin{subfigure}{.3\linewidth}
        \centering
        \includegraphics[scale=.225]{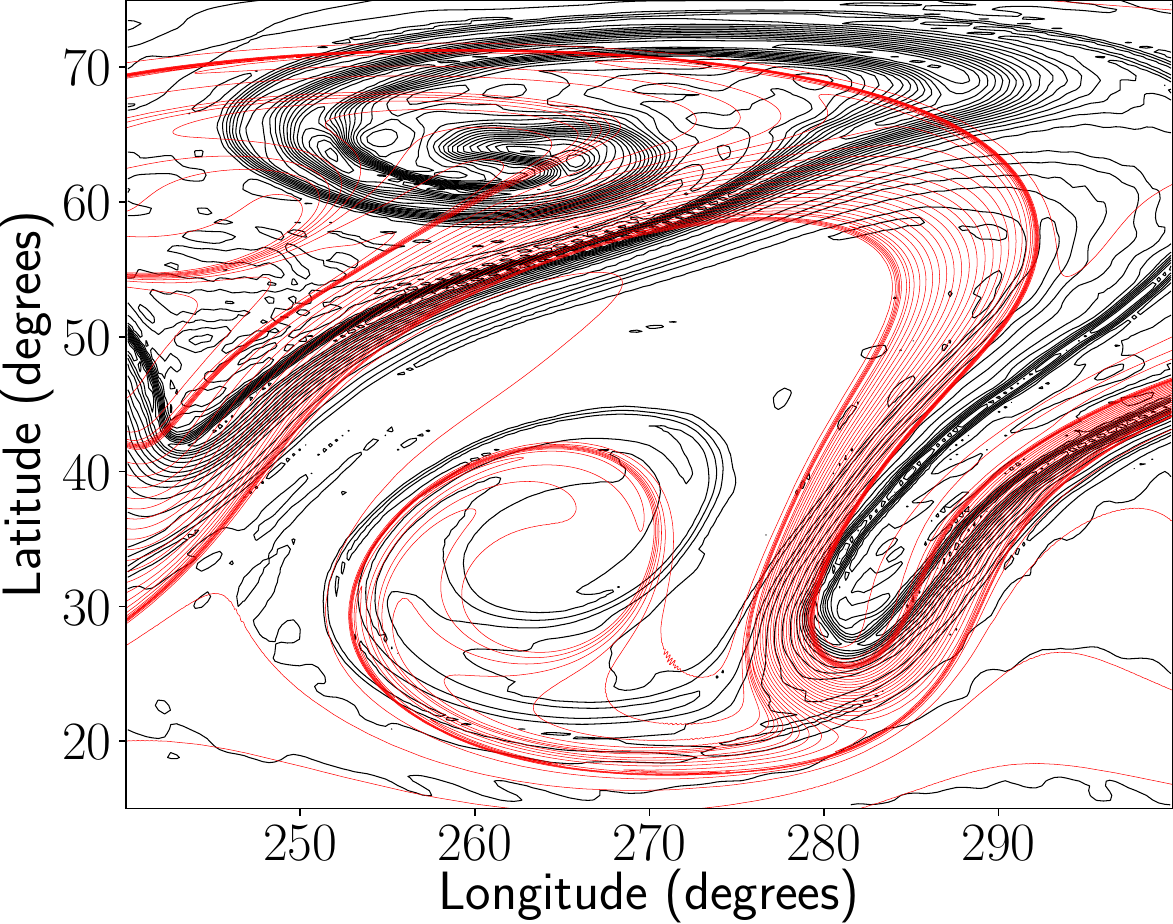}
        \caption{\SE{1}{2}, $\Dt = 960\text{s}$}
    \end{subfigure}
    \begin{subfigure}{.3\linewidth}
        \centering
        \includegraphics[scale=.225]{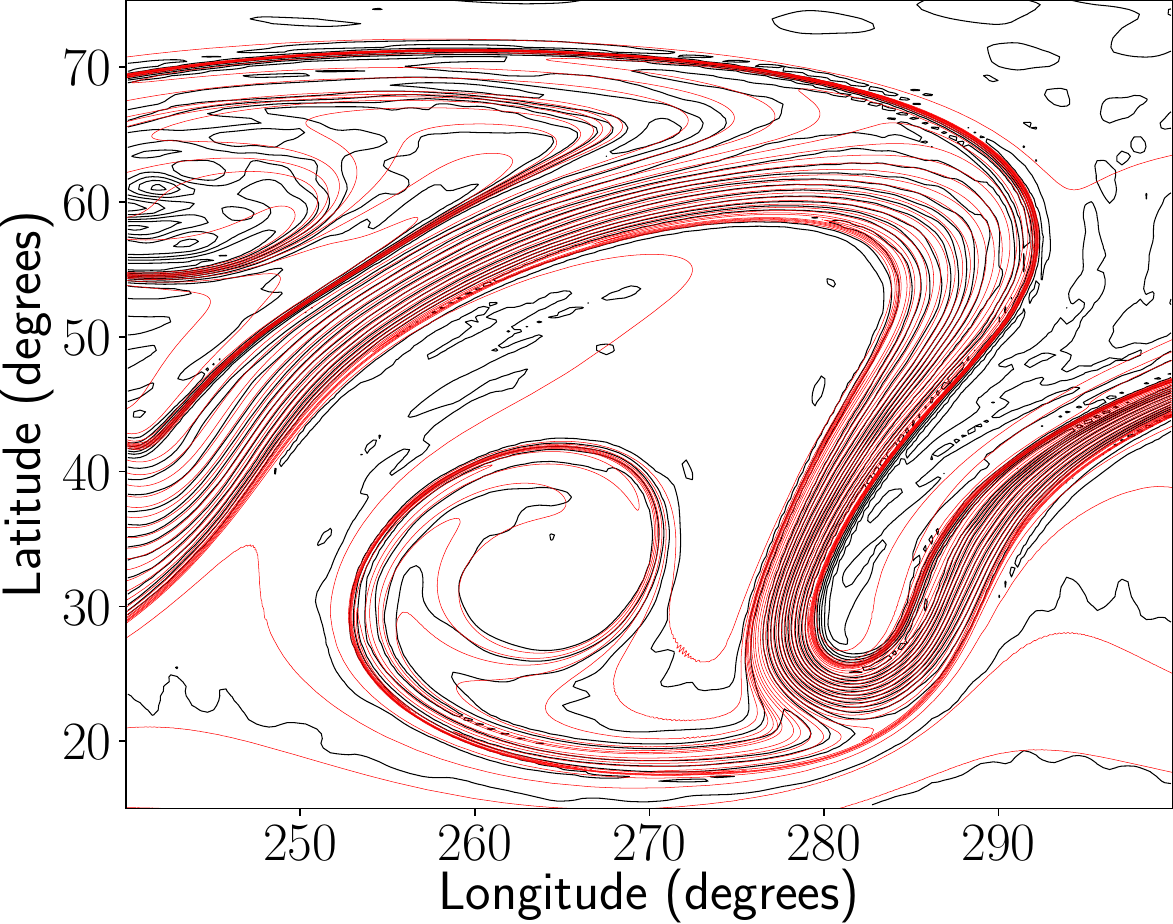}
        \caption{\SE{2}{2}, $\Dt = 960\text{s}$}
    \end{subfigure}
    \caption{Contour lines of the potential vorticity field after six days of simulation of the unstable jet test case, using selected time integration methods (black lines), compared to the reference solution provided by \cite{scott_al:2015} (red lines). Contour levels are separated by $0.2\Omega/\ol{h}$.}
    \label{fig:contour}
\end{figure}

\indent We now compare the obtained results using the guidelines provided by \cite{scott_al:2015} for the study of the unstable jet test case. First, we present, in Figure \ref{fig:contour}, the contour lines of the potential vorticity field of each simulation after six days on the region $[240^{\circ},300^{\circ}]\times[15^{\circ},75^{\circ}]$, compared to the reference solution provided by \cite{scott_al:2015}, which is obtained considering a spectral resolution $M = 2730$. The potential vorticity is defined by $\xi_p := \xi_a/h$, where $\xi_a := 2 \Omega \sin \theta + \xi = f + \xi$ is the absolute vorticity. For the sake of conciseness, among the semi-Lagrangian exponential methods, only results produced by \SE{1}{2} and \SE{2}{2} are presented. With $\Dt = 120\text{s}$, small-scale oscillating errors are observed in the solution produced by SL-SI-SETTLS, which is also slightly displaced \wrt the reference one; both \SE{1}{2} and \SE{2}{2} produce solutions close to the reference contour, with a better representation of some small-scale details by the latter scheme. With $\Dt = 960\text{s}$, small-scale and displacement errors are still present in the SL-SI-SETTLS solution; concerning \SE{2}{2}, it still represents well the general profile and position of the reference solution, indicating smaller dispersion errors compared to the other schemes, but clear oscillations reveal the formation of unstable behaviors; finally, the solution of \SE{1}{2} is strongly deteriorated both in shape and position.

\begin{figure}[!htbp]
    \centering
    \includegraphics[scale=.475]{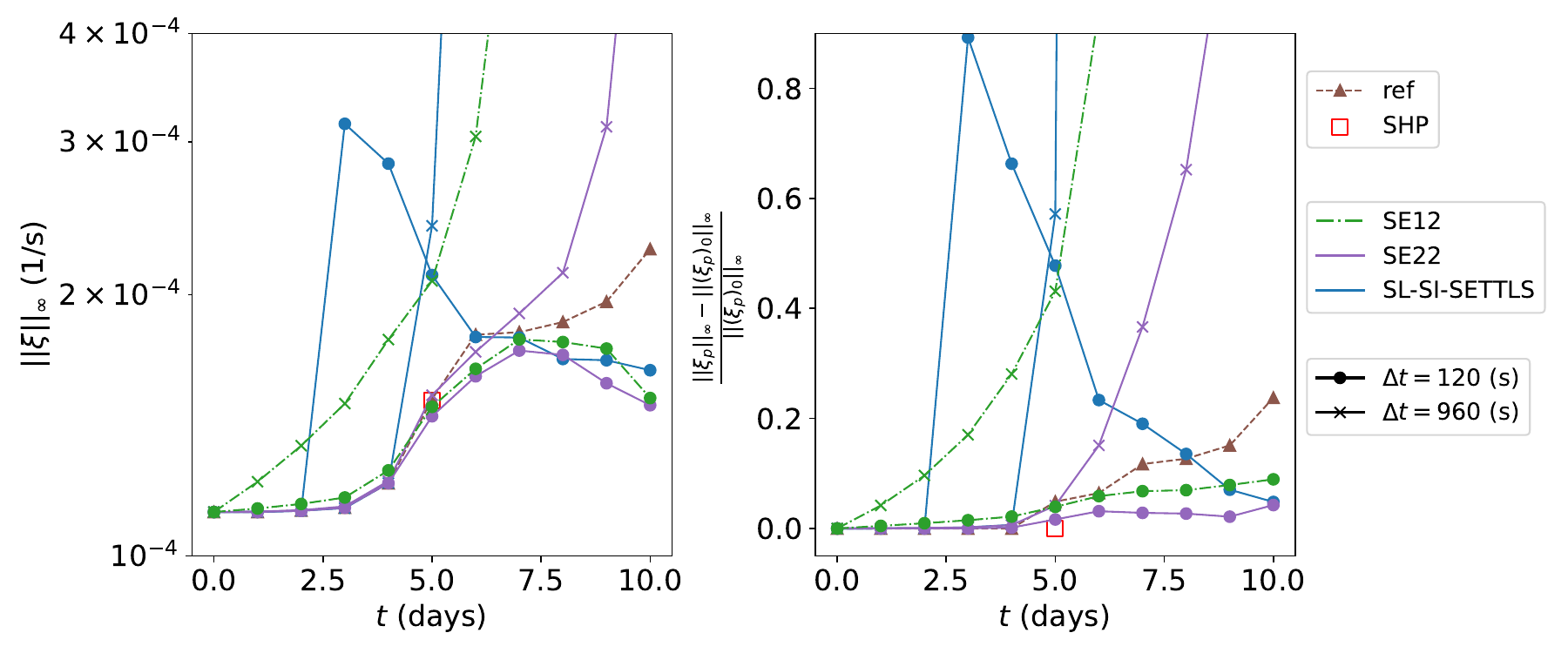}
    \caption{Evolution of diagnostic quantities in the unstable jet test case, measured at every day of simulation. \quotes{ref} indicates the reference solution (4th order explicit Runge-Kutta, with $M=1024$ and $\Dt = 2\text{s}$) and \quotes{SHP} indicates the reference values provided by Scott, Harris and Polvani \cite{scott_al:2015} at day 5. Left: maximum absolute value of the vorticity field; right: deviation of the maximum absolute value of the potential vorticity field.}
    \label{fig:diagnostics}
\end{figure}

\indent \cite{scott_al:2015} also propose a set of metrics for a quantitative evaluation of the convergence of numerical solutions of the unstable jet test case. Among them, the maximum absolute value of the vorticity ($\norm{\infty}{\xi}$) and the relative error of the maximum potential vorticity ($(\norm{\infty}{\xi_p} - \norm{\infty}{(\xi_p)_0}) / \norm{\infty}{(\xi_p)_0} $, where $(\xi_p)_0$ is the initial potential vorticity), are identified to be proper convergence diagnostics. These quantities at every day of simulation are plotted in Figure \ref{fig:diagnostics} for the solutions computed with SL-SI-SETTLS, \SE{1}{2} and \SE{2}{2}, being compared to the reference solution, as well as to reference values provided by \cite{scott_al:2015} at day 5. With the small time step $\Dt = 120\text{s}$, the semi-Lagrangian exponential schemes outperform SL-SI-SETTLS, and they even provide better approximations to the reference value of the potential vorticity deviation compared to our reference solution, which, as discussed above, presents a slightly unstable behavior at the very end of the wavenumber spectrum. Under the larger time step $\Dt = 960\text{s}$, the instabilities of \SE{1}{2} are clear, with a fast divergence from the reference values from the first day of simulation; SL-SI-SETTLS and \SE{2}{2} have a much later divergence, mainly the latter, which stays relatively close to our reference solution until day 6 or 7.

\subsection{Evaluation of computational cost}

\indent We compare, in Table \ref{tab:times_unstable_jet}, the wall-clock times for the integration of the simulations whose kinetic energy spectra are presented in Figures \ref{fig:spectrum_serial_120} and \ref{fig:spectrum_serial_960}. Each simulation is executed with threaded spatial parallelization in 16 physical cores of Intel Xeon Gold 6130 \@2.10GHz in the GRICAD cluster from the University of Grenoble Alpes. As expected from the complexity analysis depicted in Table \ref{tab:complexity}, the proposed modification of the discretization of the linear term strongly increases the computing time of the semi-Lagrangian exponential schemes, by factors ranging approximately between 40\% and 65\%. However, the proposed \SE{2}{2} scheme, which provided the most stable and accurate results, presents a comparable cost with the well-established SL-SI-SETTLS. Moreover, the evaluation of the exponential functions is not fully optimized in our implementation, such that further wall clock time improvements could be obtained.

\begin{table}[!htbp]
    \centering
    \begin{tabular}{|c|c|c|c|c|c|}
        \cline{2-6}
          \multicolumn{1}{c|}{} & SL-SI-SETTLS & \SE{1}{1} & \SE{2}{1} & \SE{1}{2} & \SE{2}{2} \\
         \hline
         $\Dt = 120\text{s}$ & 1378 & 757  & 1051  & 1244 & 1538 \\
         \hline
         $\Dt = 960\text{s}$ & 178 & 97 & 132 & 149 & 190 \\
         \hline
    \end{tabular}
    \caption{Approximate wall-clock times (in seconds) for the integration of the unstable jet test case (10 days of simulation) by selected schemes.}
    \label{tab:times_unstable_jet}
\end{table}

\subsection{Application of artificial viscosity and hyperviscosity}
\label{subsec:viscosity}

\indent No artificial viscosity or hyperviscosity was applied in the numerical simulations presented above. Although relatively stable simulations were indeed obtained, mainly with the proposed \SE{2}{2} scheme, we identified some indications that viscosity approaches should be used for further stability, \eg the lack of convergence in the geostrophic balance test case (Figure \ref{fig:error_serial_convergence_geostrophic_balance}), the upturned tails at the end of the wavenumber spectra (Figure \ref{fig:spectrum_serial}) and the localized peaks in the spectra of semi-Lagrangian exponential methods when specific time step sizes are used (Figures \ref{fig:simulation_stability} and \ref{fig:spectrum_serial_240}). The goal of this section is to verify if artificial viscosity and hyperviscosity approaches are able to mitigate these issues and to identify a compromise between stability and accuracy when choosing the viscosity parameters, namely the viscosity order $q$ and the viscosity coefficient $\nu$. In practice, artificial viscosity and/or numerical diffusion inherent to the discretization schemes are required in operational global atmospheric circulation models \cite{thuburn:2008}. Ideally, one would use higher-order viscosity approaches ($q > 2$), in order to damp only the largest wavenumbers (the ones in which instabilities may be triggered first due to the \replaced[id=R2-new]{enstrophy}{energy} cascade induced by the nonlinearity of the problem) and preserve the largest ones, thus maintaining relatively high accuracy. Fourth-, sixth- and eight-order hyperviscosity approaches are usually considered in models using spectral discretization \cite{jablonowski_williamson:2011}; for instance, the spectral models IFS-ECMWF \cite{ECMWF:2003} and GFS-NCEP/NOAA \cite{NOAA_NCEP:2016} use respectively fourth- and eight-order viscosities. The determination of viscosity coefficients depends on the spectral resolution and the damping timescale, \ie the time necessary for damping a given wavenumber by a given fraction, which varies between units and tenths of hours. We refer the reader to Appendix A of \cite{caldas_al:2023} for details and to the references therein for reports of viscosity parameters usually applied in spectral methods for atmospheric circulation models.

\begin{figure}[!htbp]
    \begin{subfigure}{.475\linewidth}
        \centering
        \includegraphics[scale=.35]{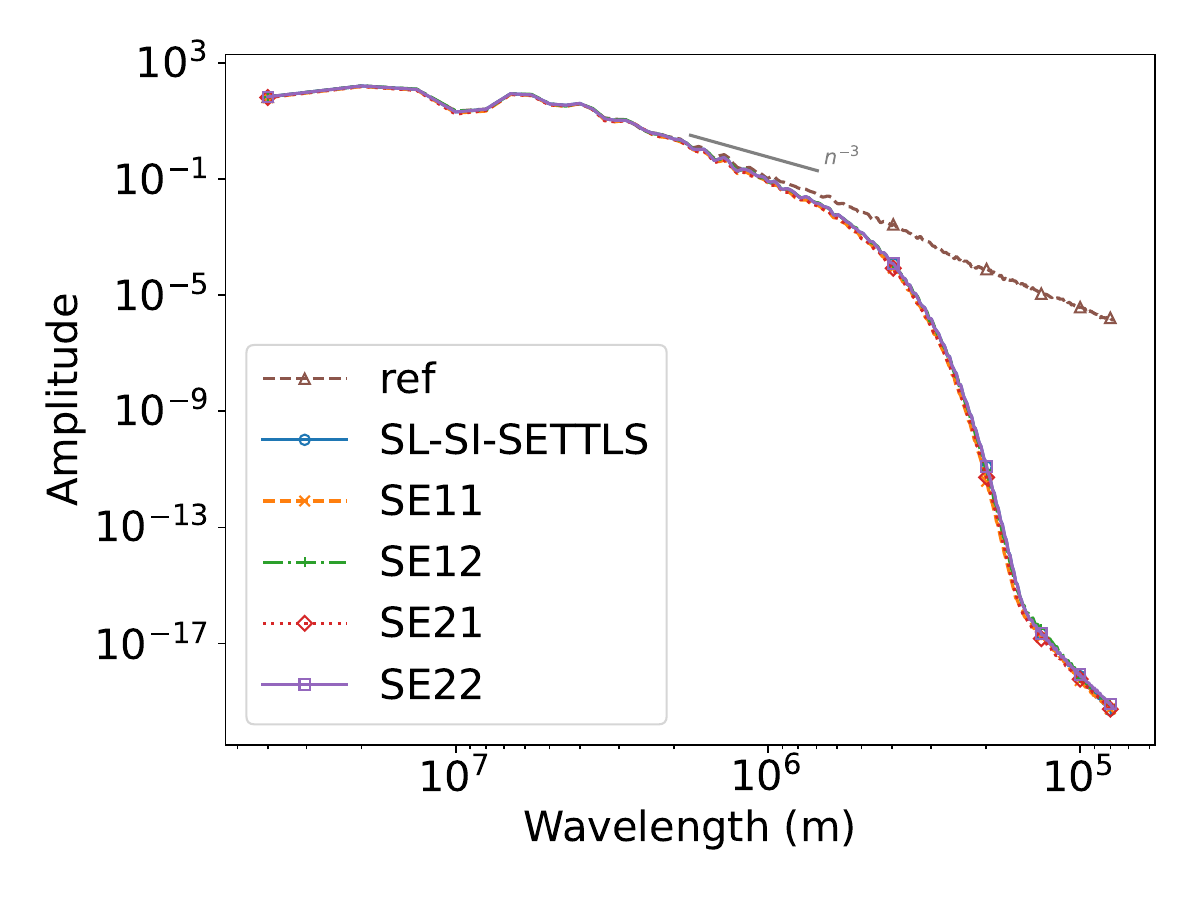}
        \caption{$(q,\nu) = (4, 10^{15})$}
    \end{subfigure}
    \begin{subfigure}{.475\linewidth}
        \centering
        \includegraphics[scale=.35]{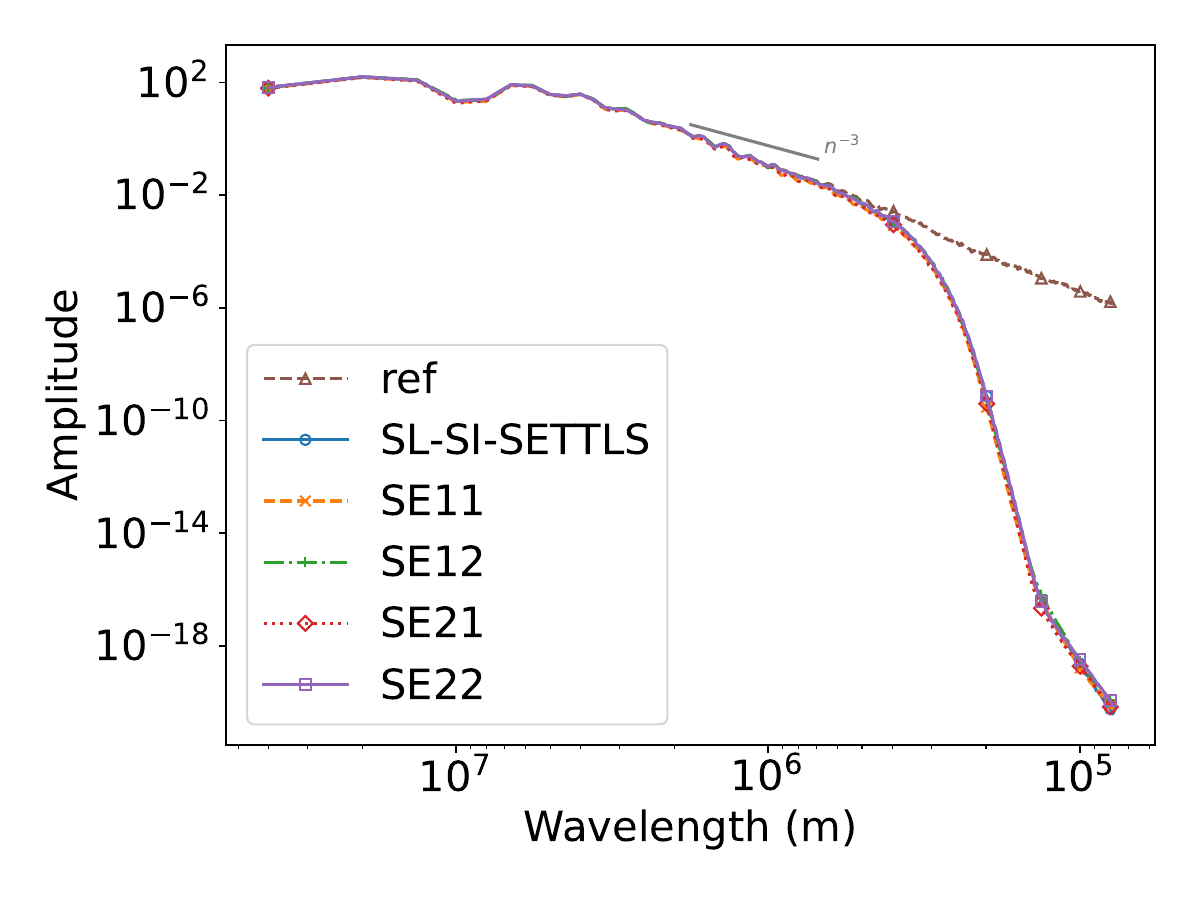}
        \caption{$(q,\nu) = (6, 10^{24})$}
    \end{subfigure}
    \caption{Kinetic energy spectra after six days, with $\Dt = 240\text{s}$ (and $\Dt = 2\text{s}$ for the reference solution), truncated at resolution $M = 512$, for the unstable jet case with artificial viscosity.}
    \label{fig:spectrum_viscosity}
\end{figure}

\begin{figure}[!htbp]
    \begin{subfigure}{.45\linewidth}
        \centering
        \includegraphics[scale = .275]{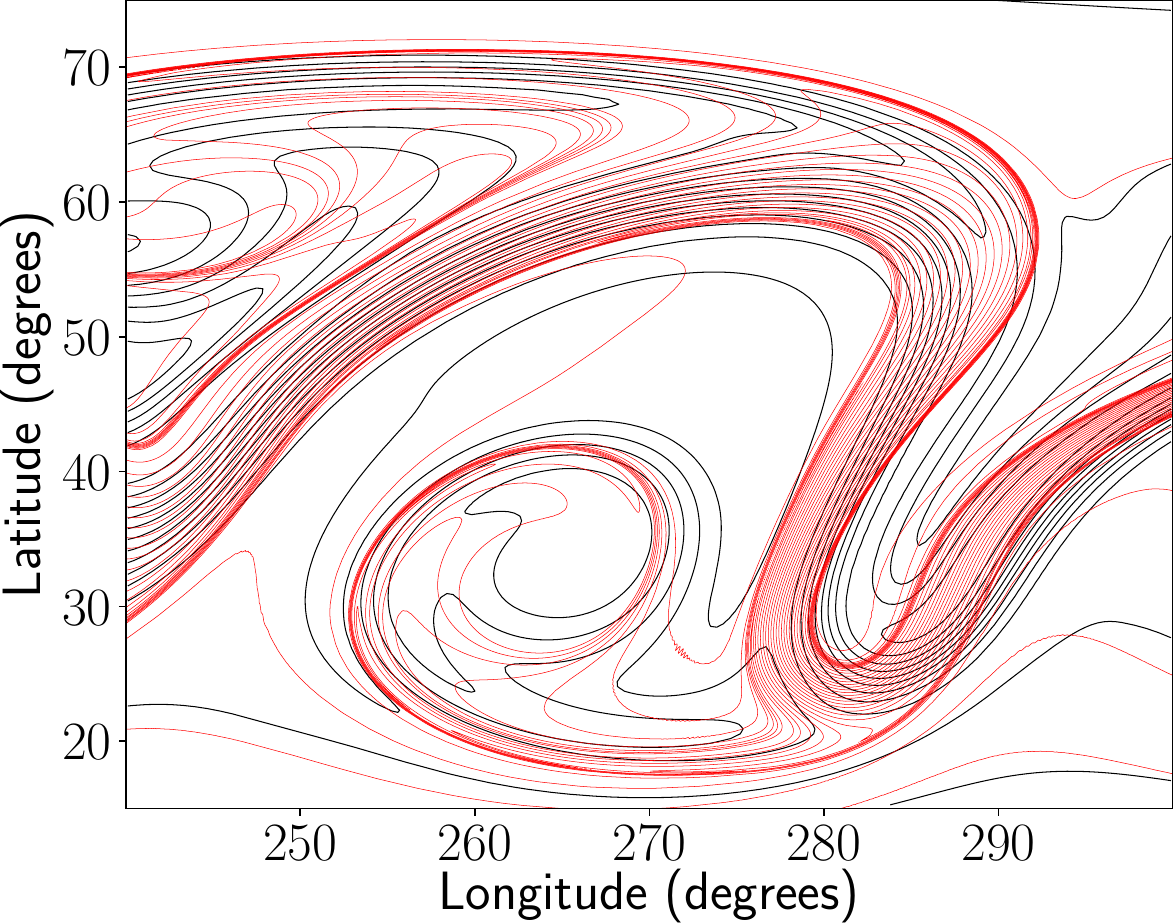}
        \caption{$(q,\nu) = (4, 10^{15})$}
    \end{subfigure}
    \begin{subfigure}{.45\linewidth}
        \centering
        \includegraphics[scale=.275]{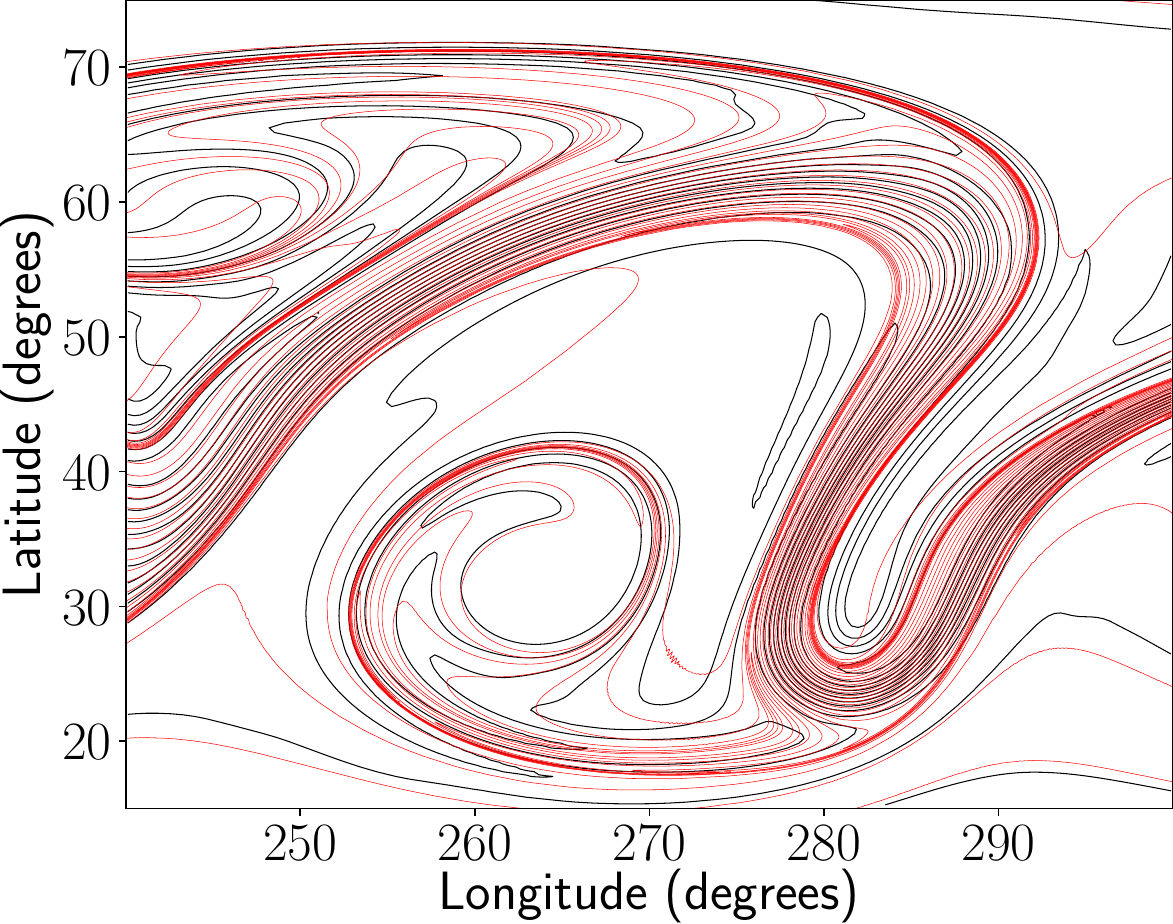}
        \caption{$(q,\nu) = (6, 10^{24})$}
    \end{subfigure}
    \caption{Contour lines of the potential vorticity field after six days of the unstable jet test case with artificial viscosity, using the \SE{2}{2} scheme and $\Dt = 240\text{s}$ (black lines), compared to the reference solution provided by \cite{scott_al:2015} (red lines). Contour levels are separated by $0.2\Omega/\ol{h}$.}
    \label{fig:contour_lines_viscosity}
\end{figure}

\indent In all simulations presented here, the viscosity term $\vecL_{\nu}$ is solved at the end of each time step using a first-order backward Euler scheme. We consider the unstable jet test case with $\Dt = 240 \text{s}$, for which clear instabilities were observed in the integration of the semi-Lagrangian exponential schemes, with $(q,\nu) = (4, 10^{15})$ and  $(q,\nu) = (6, 10^{24})$. Smaller viscosity coefficients were not able to fully prevent the formation of peaks in the kinetic energy spectra. As presented in Figure \ref{fig:spectrum_viscosity}, both viscosity configurations ensure stability at the cost of important damping of the highest wavenumber modes (and stronger damping of intermediate modes with the fourth-order viscosity), with very close results for all tested time integration methods. As illustrated in Figure \ref{fig:contour_lines_viscosity} for the \SE{2}{2} scheme (qualitatively similar results are obtained using \SE{1}{2} or SL-SI-SETTLS), the sixth-order viscosity produces a better approximation to the reference solution provided by \cite{scott_al:2015}.

\indent By comparing these results with the simulations without viscosity (Figure \ref{fig:contour}), we observe that \SE{2}{2} with $\Dt = 120\text{s}$, $\nu = 0$ is qualitatively superior to the simulation with $\Dt = 240\text{s}$ and $(q,\nu) = (6,10^{24})$, but the latter outperforms the simulation with $\Dt = 960\text{s}$, $\nu = 0$, which presents several small-scale oscillations despite a relatively stable simulation compared to other time integration schemes. On the other hand, SL-SI-SETTLS and \SE{1}{2} clearly benefit from the use of artificial viscosity, with considerably more stable and accurate results compared to the inviscid case with both $\Dt = 120\text{s}$ and $\Dt = 960\text{s}$ (in the case of SL-SI-SETTLS) or $\Dt = 960\text{s}$ (in the case of \SE{1}{2}).

\section{Conclusion}
\label{sec:conclusion}

\indent Semi-Lagrangian exponential methods, proposed by \cite{peixoto_schreiber:2019}, have a potential for application in atmospheric circulation models, since they have enhanced stability compared to Eulerian exponential schemes and, contrary to traditional semi-Lagrangian methods, are able to integrate the linear terms very accurately; therefore, this class of time stepping scheme is more suitable for the integration of advection-dominated models containing stiff linear terms. The work presented here aimed to extend \cite{peixoto_schreiber:2019} in several aspects: first, we identified and explained the accuracy limitations of their methods, showing that the combination of exponential and interpolation operators adopted in the discretization of the linear term\xadded{, which corresponds to a first-order discretization of the integration factor,} produces a first-order leading term in the truncation error. Second, inspired by this initial study, we proposed a new method, with an alternative discretization of the linear term, leading to an effective second-order accuracy. Also, we conducted a detailed stability study comparing several Eulerian and semi-Lagrangian methods, both through a linear stability analysis and an empirical simulation-based stability study, which was required because the standard stability analysis is not able to provide detailed insights in the stability of semi-Lagrangian exponential schemes. Finally, while \cite{peixoto_schreiber:2019} developed their work in the context of the SWE on the biperiodic plane, we applied the schemes to the integration of the SWE on the rotating sphere, thus providing results that are more closely related to real-world operational atmospheric circulation models.

\indent Promising results were obtained using the proposed second-order semi-Lagrangian exponential scheme. Its theoretical convergence order was confirmed in several benchmark test cases, and, considering the especially challenging unstable jet simulation, it proved to be more stable and accurate than Eulerian exponential schemes, other semi-Lagrangian exponential schemes and the well-established second-order SL-SI-SETTLS, even without the use of artificial viscosity approaches, which is a consequence of the combined accurate integration of the stiff linear terms and the semi-Lagrangian treatment of the advection. The latter approach largely extends the range of stable time step sizes compared to a second-order Eulerian exponential method, maintaining approximately the same accuracy, while the former leads to larger accuracy, compared to SL-SI-SETTLS, in simulations dominated by linear processes; however, our proposed scheme is outperformed by SL-SI-SETTLS (if stable) when nonlinear effects are more important. In terms of computational cost, a drawback of the proposed method is its higher complexity compared to the previously formulated first-order semi-Lagrangian exponential schemes, as a consequence of the splitting of the exponential operator applied to the linear term, which is required to achieve second-order accuracy. We highlight, however, that this larger cost can be mitigated by a straightforward computation of matrix exponentials, by combining a spectral discretization with a convenient arrangement of the terms of the governing equations, and by precomputing these exponential terms, which is possible in simulations with a constant time step size. Also, the observed wall-clock times are comparable to those of SL-SI-SETTLS, which is used in operational atmospheric circulation models; therefore, our proposed method seems to be a competitive candidate for operational applications.

\indent Important aspects still need to be considered in future work for further improvement of semi-Lagrangian exponential methods. As illustrated in the numerical simulations, these schemes presented unexpected unstable behaviors in a specific range of time step sizes, which can possibly be caused by badly behaved exponential functions; even if these issues have been easily avoided by using high-order artificial viscosity, with little damage to the accuracy of the numerical solution, a more detailed understanding needs to be achieved. Also, higher-order semi-Lagrangian exponential schemes can be formulated to achieve improved accuracy in simulations dominated by nonlinearities; from the truncation error analyses conducted here, they could probably be obtained in a quite straightforward way following the derivation of higher-order Eulerian exponential methods by \cite{cox_matthews:2002} and using higher-order approximations to the integral involved in the integration factor term\xadded{; however, the hypothesis on the commutation of the linear operator along the Lagrangian trajectories, which was proved here to introduce a second-order error, should be carefully studied.} \xadded{Moreover, the discussion on the derivation of linear stability analyses of semi-Lagrangian exponential schemes combined with spectral methods, conducted in Section \ref{subsec:stability_discussion}, suggests that more intricate studies should be necessary to fully understand the relation between spectral transforms, matrix exponentiation, and spatial interpolation, and the influence on accuracy properties of semi-Lagrangian exponential methods, which would certainly contribute to further developments of this class of methods.} Finally, applications to more complex atmospheric models would be a natural and challenging extension of this work.

\section*{Acknowledgements}

% \indent \textbf{Funding:} This work was supported by the São Paulo Research Foundation (FAPESP) grants 2021/03777-2 and 2021/06176-0, as well as the Brazilian National Council for Scientific and Technological Development (CNPq), Grant 303436/2022-0. This project also received funding from the Federal Ministry of Education and Research and the European High-Performance Computing Joint Undertaking (JU) under grant agreement No 955701, Time-X. The JU receives support from the European Union’s Horizon 2020 research and innovation programme and Belgium, France, Germany, Switzerland.

\indent This work was supported by the São Paulo Research Foundation (FAPESP) grants 2021/03777-2 and 2021/06176-0, as well as the Brazilian National Council for Scientific and Technological Development (CNPq), Grant 303436/2022-0. This project also received funding from the Federal Ministry of Education and Research and the European High-Performance Computing Joint Undertaking (JU) under grant agreement No 955701, Time-X. The JU receives support from the European Union’s Horizon 2020 research and innovation programme and Belgium, France, Germany, Switzerland.

\indent Most of the computations presented in this paper were performed using the GRICAD infrastructure (\url{https://gricad.univ-grenoble-alpes.fr}), which is supported by Grenoble research communities.

\section*{Data availability statement}

\indent The research code, benchmark tests and postprocessing scripts used in this article are available in the repository \url{https://gitlab.inria.fr/sweet/sweet/-/tree/paper_2nd_order_SLETDRK/benchmarks_sphere/paper_2nd_order_sletdrk} \cite{code}.

\bibliographystyle{abbrv}
\bibliography{biblio}

\end{document}